\documentclass[12pt]{article}

\usepackage{amsmath}
\usepackage{amssymb}
\usepackage{amsthm}
\usepackage{amsfonts}
\usepackage{amscd}
\usepackage[dvipsnames]{xcolor}
\usepackage[T1]{fontenc}  
\usepackage[nottoc]{tocbibind}  
\usepackage{enumitem} 
\usepackage{caption}
\usepackage[numbers,sort&compress]{natbib}
\usepackage{hyperref}
\usepackage{cleveref} 
\usepackage{adjustbox}
\usepackage{fancyvrb} 
\usepackage{tikz}
\usepackage{verbatim} 
\usepackage{framed}
\usepackage{geometry}
\usepackage[medium]{titlesec}
\usepackage{etoolbox} 
\usepackage{arydshln} 
\usepackage{tocloft} 

\theoremstyle{plain}
\newtheorem{lemma}{Lemma}
\newtheorem{proposition}[lemma]{Proposition}
\newtheorem{theorem}{Theorem}
\newtheorem{corollary}{Corollary}
\newtheorem{conjecture}{Conjecture}
\newtheorem{theoremext}{Theorem}

\theoremstyle{definition}
\newtheorem{definition}[lemma]{Definition}
\newtheorem{remark}[lemma]{Remark}
\newtheorem{example}[lemma]{Example}

\definecolor{colLinkBlue}{RGB}{23,111,192}
\definecolor{colCiteGreen}{RGB}{8,144,8}

\hypersetup{colorlinks,urlcolor=colCiteGreen,citecolor=colCiteGreen,linkcolor=colLinkBlue}
\makeatletter\patchcmd{\ttlh@hang}{\parindent\z@}{\parindent\z@\leavevmode}{}{}\patchcmd{\ttlh@hang}{\noindent}{}{}{}\makeatother 
\titlespacing*{\section}{0pt}{1mm}{1mm}
\titlespacing*{\subsection}{0pt}{0mm}{0mm}
\titlespacing*{\paragraph}{0pt}{1mm}{1mm}
\newcommand{\myspace}{\setlength{\abovedisplayskip}{1mm}\setlength{\belowdisplayskip}{0mm}}
\usetikzlibrary{positioning}
\usetikzlibrary{arrows.meta}

\newenvironment{Mlist}{\begin{itemize}[topsep=0pt,itemsep=0pt,leftmargin=7mm]}{\end{itemize}}
\newenvironment{claims}{\begin{enumerate}[topsep=0pt,itemsep=0pt,leftmargin=8mm,label=\textnormal{\textbf{(\alph*)}},ref=(\alph*)]}{\end{enumerate}}

\newcommand{\csep}[1]{\setlength{\tabcolsep}{#1}}
\newcommand{\fig}[3]{\includegraphics[height=#1cm, width=#2cm]{img/#3}}
\newcommand{\SEC}[1]{\textsection\ref{sec:#1}}
\newcommand{\ASN}[1]{Assertion~\ref{#1}}
\newcommand{\END}{\hfill $\vartriangleleft$}
\newcommand{\RL}[2]{\Cref{lem:#1}\ref{lem:#1:#2}}
\newcommand{\RLS}[4]{Lemmas~\ref{lem:#1}\ref{lem:#1:#2} and~\ref{lem:#3}\ref{lem:#3:#4}}
\newcommand{\RP}[2]{\Cref{prp:#1}\ref{prp:#1:#2}}
\newcommand{\RPS}[4]{Propositions~\ref{prp:#1}\ref{prp:#1:#2} and~\ref{prp:#3}\ref{prp:#3:#4}}
\newcommand{\df}[1]{{\it #1}}
\newcommand{\Wlog}{without loss of generality }
\newcommand{\resp}{respectively}
\newcommand{\st}{such that }
\newcommand{\wrt}{with respect to }
\newcommand{\Iff}{if and only if }
\renewcommand{\c}{\colon}
\newcommand{\set}[2]{\{#1 ~|~ #2\}}
\newcommand{\dto}{\dasharrow}

\newcommand{\bas}[1]{\langle #1\rangle}
\newcommand{\SingType}{\operatorname{SingType}}
\newcommand{\Sing}{\operatorname{Sing}}
\renewcommand{\P}{{\mathbb{P}}}
\renewcommand{\S}{{\mathbb{S}}}
\newcommand{\C}{{\mathbb{C}}}
\newcommand{\R}{{\mathbb{R}}}
\newcommand{\Z}{{\mathbb{Z}}}
\newcommand{\Q}{{\mathbb{Q}}}
\newcommand{\U}{{\mathbb{U}}}
\newcommand{\E}{{\mathbb{E}}}
\newcommand{\Y}{{\mathbb{Y}}}
\renewcommand{\H}{{\mathbb{H}}}
\newcommand{\bP}{{\mathbf{P}}}
\newcommand{\bS}{{\mathbf{S}}}
\newcommand{\bR}{{\mathbf{R}}}
\newcommand{\bC}{{\mathbf{C}}}
\newcommand{\bO}{{\mathbf{O}}}
\newcommand{\cE}{\mathcal{E}}
\newcommand{\cF}{\mathcal{F}}
\newcommand{\cG}{\mathcal{G}}
\newcommand{\cL}{\mathcal{L}}
\newcommand{\cM}{\mathcal{M}}
\newcommand{\cR}{\mathcal{R}}
\newcommand{\cW}{\mathcal{W}}

\newcommand{\oL}{\overline{L}}
\newcommand{\oR}{\overline{R}}
\newcommand{\oM}{\overline{M}}
\newcommand{\oT}{\overline{T}}
\newcommand{\oF}{\overline{F}}

\renewcommand{\l}{\ell}
\newcommand{\p}{\varepsilon}

\newcommand{\pp}{\mathfrak{p}}
\newcommand{\op}{\overline{\mathfrak{p}}}
\newcommand{\qq}{\mathfrak{q}}
\newcommand{\oq}{\overline{\mathfrak{q}}}
\renewcommand{\aa}{\mathfrak{a}}
\newcommand{\oa}{\overline{\mathfrak{a}}}

\newcommand{\ii}{\mathfrak{i}}
\newcommand{\HC}{\mathbb{H}_{\C}}
\newcommand{\qi}{\mathbf{i}}
\newcommand{\qj}{\mathbf{j}}
\newcommand{\qk}{\mathbf{k}}
\newcommand{\hstar}{\hat{\star}}
\newcommand{\Hstar}{\star_{\mathbb{H}}}
\newcommand{\pstar}{\star_{\mathbb{P}}}
\newcommand{\sstar}{\star_S}

\newcommand{\AN}{\varnothing}
\newcommand{\Ai}{{A_1}}

\newcommand{\uAi}{\underline{A_1}}
\newcommand{\uAii}{\underline{A_2}}
\newcommand{\uAiii}{\underline{A_3}}
\newcommand{\aut}{\operatorname{Aut}}
\newcommand{\ut}{\operatorname{UT}}
\newcommand{\lt}{\operatorname{LT}}
\newcommand{\rt}{\operatorname{RT}}

\newcommand{\diag}{\operatorname{diag}}
\newcommand{\sgn}{\operatorname{sgn}}

\newif\ifhide
\ifhide
\renewenvironment{tikzpicture}{D\comment}{\endcomment}
\renewcommand{\fig}[3]{F}
\fi

\begin{document}
\myspace

\begin{center}
~\vspace{1cm}\\
\LARGE
Translational and great Darboux cyclides
\\[5mm]\large
Niels Lubbes
\\[2mm]\large
\today
\end{center}

\begin{abstract}
A surface that is the pointwise sum of circles in Euclidean space
is either coplanar or contains no more than 2 circles through a general point.
A surface that is the pointwise product of circles in the unit-quaternions
contains either 2, 3, 4, or 5 circles through a general point.
A surface in a unit-sphere of any dimension
that contains 2 great circles through a general point
contains either 4, 5, 6, or infinitely many circles through a general point.
These are some corollaries from our classification of translational and great
Darboux cyclides. We use the combinatorics associated to the set of low degree curves on such surfaces
modulo numerical equivalence.
\\[2mm]
\\[2mm]
%
\end{abstract}

\section{Introduction}
\label{sec:intro}

In this article, we characterize surfaces in $\R^3$ that contain at least two
circles through each point.
Such surfaces are algebraic by \citep[Theorem~2]{conical}
and thus with \df{surface} we shall mean a real irreducible algebraic surface (see \SEC{model}).

Surfaces that are a union of circles in two different ways
have applications in
architecture \cite{arch},
kinematics \cite{nls-darboux,linkage}
and geometric modeling in general \cite{jut,p,kz}.
In particular, the ``Darboux cyclides'' have a long history \cite{dar,kum},
and its various properties are still
a topic of recent research \cite{ive,blum,tak,morph,pot,parabola}.
In order to clarify our main result and its relation to \cite{sko},
we recall some definitions for the non-expert.

An \df{inversion} with respect to a sphere~$O\subset \R^3$
with center~$c$ and radius~$r$ is the
map~$f\c \R^3\setminus\{c\}\to\R^3\setminus\{c\}$
\st $||x-c||\cdot ||f(x)-c||=r^2$ and the vectors $x-c$
and $f(x)-c$ are codirected for all $x\in \R^3\setminus\{c\}$.
Such a map exchanges the interior and exterior of~$O$
and takes generalized circles to generalized
circles, where a \df{generalized circle} is either a circle or a line.
We call two surfaces in $\R^3$ \df{M\"obius equivalent} if
one surface is mapped to the other by a composition of inversions.

Let $\mu\c S^3\dto\R^3$ with $\mu(y):=(y_1,y_2,y_3)/(1-y_4)$
denote the \df{stereographic projection}
from the point $(0,0,0,1)$ on the \df{3-dimensional unit-sphere}~$S^3\subset\R^4$.
The \df{M\"obius degree} of a surface $Z\subset\R^3$
is defined as~$\deg\mu^{-1}(Z)$.
A surface~$Z\subset\R^3$ is called \df{$\lambda$-circled}
if the Zariski closure of $\mu^{-1}(Z)$
contains at least $\lambda\in \Z_{\geq 0}\cup\{\infty\}$ circles through a general point.
If $\lambda\in\Z_{\geq0}$, then we assume that $Z$ is not $(\lambda+1)$-circled.
If $\lambda\geq 2$, then we call $Z$ \df{celestial}.
The M\"obius degree and $\lambda$ are both M\"obius invariants.

We may identify the unit-sphere $S^3\subset\R^4$ with
the unit quaternions and we denote the Hamiltonian product by~$\star$.
We consider the following constructions where $A$ and $B$
are curves in $\R^3$ or $S^3$:
\begin{align*}
A+B&:=\set{a+b \in \R^3}{ a\in A \text{ and } b\in B},
\\
A\star B&:=\set{a\star b \in S^3}{ a\in A \text{ and } b\in B}.
\end{align*}
Suppose that $Z\subset \R^3$ is a surface.
We call $Z$ \df{Bohemian} or \df{Cliffordian} if
there exist generalized circles $A$ and $B$ \st $Z$ is the Zariski closure of $A+B$ and $\mu(A\star B)$, \resp.
A surface that is either Bohemian or Cliffordian is called \df{translational}.
If $A$ and $B$ are great circles \st $A\star B\subset S^3$ is a surface,
then $A\star B$ is called a \df{Clifford torus}.

\begin{figure}[!ht]
\centering
\csep{2mm}
\begin{tabular}{ccccc}
\fig{2.8}{2}{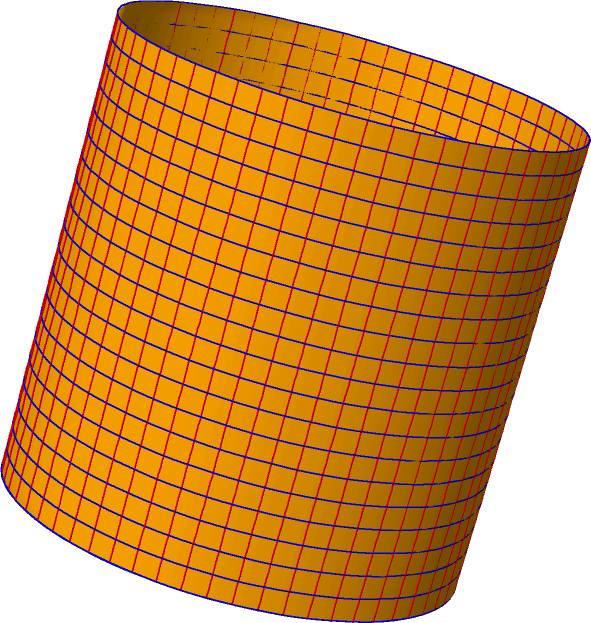} &
\fig{2.8}{2}{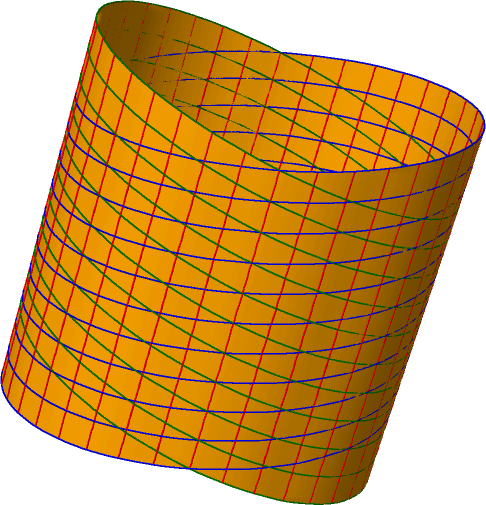} &
\fig{2.8}{2.8}{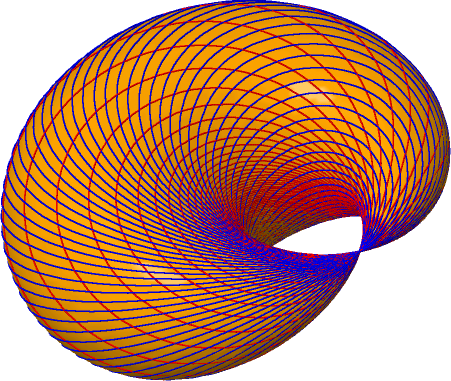} &
\fig{2.8}{2.8}{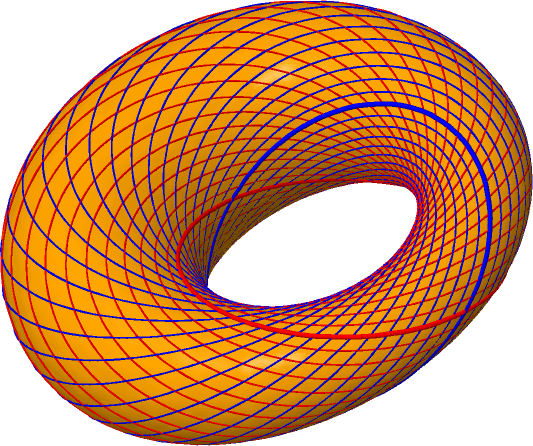} &
\fig{2.8}{2.8}{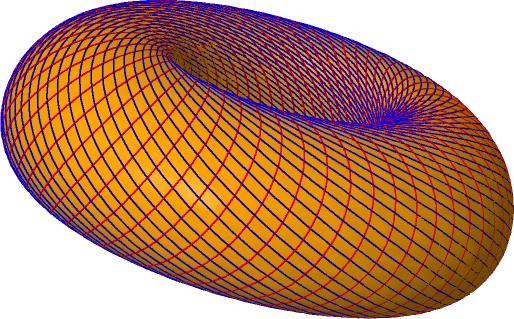}
\\
CY  &
EY  &
CH1 cyclide &
ring cyclide &
Perseus cyclide
\end{tabular}
\caption{Examples of Darboux cyclides.}
\label{fig:4}
\end{figure}

A \df{Darboux cyclide} in~$\R^3$ is a surface of M\"obius degree four.
A \df{$Q$ cyclide} is a Darboux cyclide that is M\"obius equivalent to
a quadric~$Q$.
For example, a \df{CH1~cyclide} is M\"obius equivalent
to a {C}ircular {H}yperboloid of {1} sheet (see \Cref{fig:4}),
where we used the following abbreviations:
\\[2mm]
\begin{tabular}{llll}
E = elliptic/ellipsoid      & P = parabolic/paraboloid   &  O = cone     \\
C = circular                & H = hyperbolic/hyperboloid &  Y = cylinder \\
\end{tabular}
\\[2mm]
A CO cyclide and CY cyclide is also known as a \df{spindle cyclide}
and \df{horn cyclide}, \resp.
A \df{ring cyclide}, \df{Perseus cyclide} or \df{Blum cyclide}
is a Darboux cyclide without real singularities
that is $4$-circled, $5$-circled and $6$-circled, \resp~(see \Cref{fig:4}).
See \Cref{tab:B} for a complete list of names for celestial Darboux cyclides.

It follows from \citep[Main~Theorem~1.1]{sko}
that a celestial surface in~$\R^3$
is either a Darboux cyclide or M\"obius equivalent to
a Bohemian or Cliffordian surface.
The following question arises:
{\it what are the Bohemian and Cliffordian Darboux cyclides?}
We shall provide necessary conditions using the combinatorics
of divisor classes of curves on such surfaces.
We also classify ``great'' celestial Darboux cyclides;
we call a surface $Z\subset \R^3$ \df{great} if its inverse stereographic projection~$\mu^{-1}(Z)$ is covered by great circular arcs.

We will use Theorems~\ref{thm:B} and~\ref{thm:8} and therefore build on~\cite{circle}.

\begin{theorem}
\label{thm:x}
Suppose that $Z\subset\R^3$ is a $\lambda$-circled surface of M\"obius degree~$d$
\st $\lambda\geq 2$ and  $(d,\lambda)\neq(8,2)$.
\begin{claims}

\item\label{thm:x:a}
The surface~$Z$ is Bohemian if and only if $Z$ is either a plane, CY or EY.

\item\label{thm:x:b}
If $Z$ is Cliffordian, then $Z$ is either a
Perseus cyclide,
ring cyclide or
\\CH1~cyclide.
Conversely, if $Z$ is a ring cyclide, then $Z$ is M\"obius equivalent
to a Cliffordian surface.

\item\label{thm:x:c}
The surface $Z$ is M\"obius equivalent to a great celestial surface
if and only if
$Z$ is either a
plane,
sphere,
Blum cyclide,
Perseus~cyclide,
ring cyclide,
EO cyclide or
CO cyclide.
\end{claims}
\end{theorem}

\begin{remark}
\Cref{thm:x}\ref{thm:x:b} should be considered the main result of this article,
since
Theorems~\ref{thm:x}\ref{thm:x:a} and~\ref{thm:x}\ref{thm:x:c} can
be proven using alternative and well-known methods from \cite{cool,pot,tak}
(see \Cref{rmk:boh,rmk:great} for details).
However, we propose that our proof methods for \ref{thm:x:a} and \ref{thm:x:c}
provide additional geometric insight, some of which is
made precise by \Cref{cor:moebius,cor:bohemian,cor:g2,cor:great} below
(see also \Cref{rmk:ring,rmk:greatblum,rmk:greatperseus}).
\END
\end{remark}

We summarized \Cref{thm:x} in \Cref{tab:x}.

\begin{table}[!ht]
\centering
\caption{Overview of $\lambda$-circled surfaces in $\R^3$
of M\"obius degree~$d$
that are either Bohemian, Cliffordian, or great and celestial.}
\label{tab:x}
\csep{5mm}
\begin{tabular}{cccc}
 name           & $d$ & $\lambda$  & possible types               \\\hline
plane/sphere    & $2$ & $\infty$   & Bohemian, great              \\
\hdashline
Blum cyclide    & $4$ & $6$        & great                        \\
Perseus cyclide & $4$ & $5$        & Cliffordian, great           \\
ring cyclide    & $4$ & $4$        & Cliffordian, great           \\
CH1 cyclide     & $4$ & $3$        & Cliffordian                  \\
EY              & $4$ & $3$        & Bohemian                     \\
CY              & $4$ & $2$        & Bohemian                     \\
EO cyclide      & $4$ & $3$        & great                        \\
CO cyclide      & $4$ & $2$        & great                        \\
\hdashline
                & $8$ & $2$        & Bohemian, Cliffordian, great \\
\end{tabular}
\end{table}

\begin{figure}[!ht]
\centering
\csep{6mm}
\begin{tabular}{ccc}
\fig{3}{3}{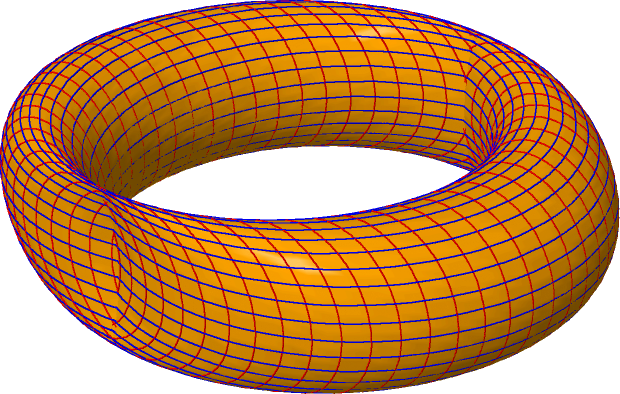} & \fig{3}{4}{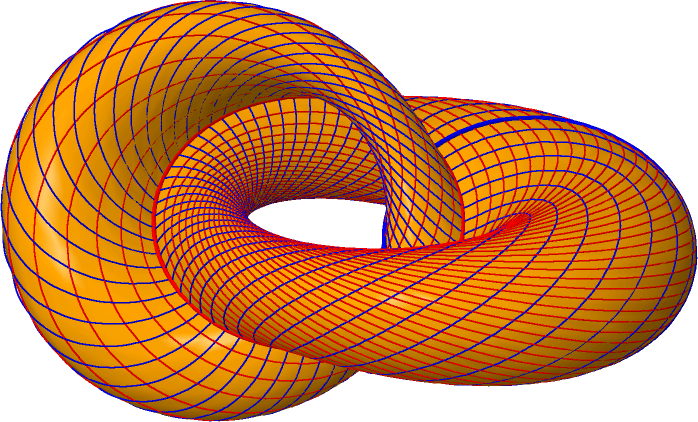} & \fig{3}{4}{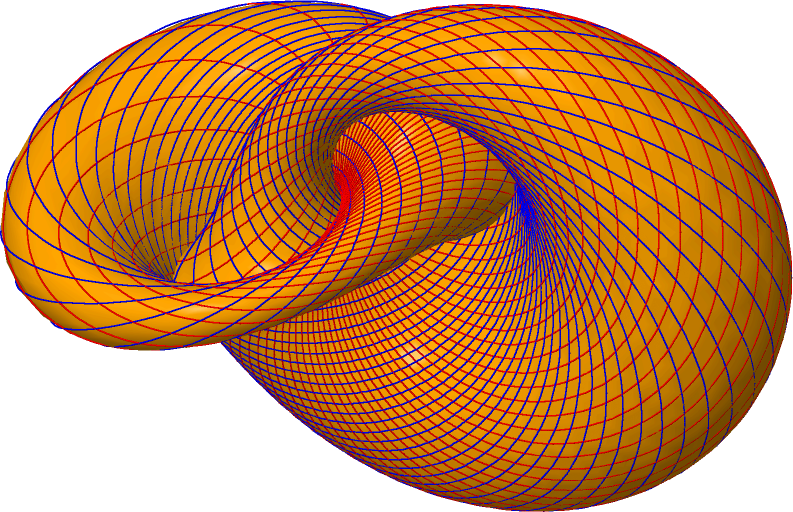}
\\
Bohemian & Cliffordian and great & Cliffordian
\end{tabular}
\caption{$2$-circled surfaces of M\"obius degree $8$ (see \Cref{exm:cert}).}
\label{fig:8}
\end{figure}

See \Cref{exm:EOCO,exm:cert,exm:RPB}
for an example for each row and each possible type.
In particular, we consider for $0\leq i,j\leq 8$ the surface~$Z_{ij}\subset\R^3$,
which is defined as the Zariski closure of a stereographic projection of the surface
\[
\set{C_i(\alpha)\star C_j(\beta)}{0\leq \alpha,\beta<2\pi}\subset S^3,
\]
where the circle parametrizations $C_i(t)$ are defined in \Cref{tab:C}.
We show in \Cref{exm:cert} that $Z_{01}$, $Z_{23}$ and $Z_{45}$
are a ring cyclide, Perseus cyclide and CH1 cyclide, \resp.
The Cliffordian surfaces $Z_{06}$ and $Z_{78}$ are of degree 8 and
illustrated in \Cref{fig:8}.

\begin{table}[!ht]
\caption{Parametrizations of circles in $S^3$ with $0\leq t<2\pi$.
Only $C_0(t)$ and $C_1(t)$ define great circles.}
\label{tab:C}
\footnotesize
\begin{align*}
C_0(t)&:=                       \Bigl(\cos(t),\sin(t),0,0               \Bigr),\\
C_1(t)&:=\tfrac{1}{5}           \Bigl(5\cos(t),4\sin(t),3\sin(t),0      \Bigr),\\
C_2(t)&:=\tfrac{1}{3}           \Bigl(\cos(t),\sin(t),-2,2              \Bigr),\\
C_3(t)&:=\tfrac{1}{3}           \Bigl(2\cos(t),2\sin(t),2,1             \Bigr),\\
C_4(t)&:=\tfrac{1}{3-2\cos(t)}  \Bigl(-2+2\cos(t),2\sin(t),0,1-2\cos(t) \Bigr),\\
C_5(t)&:=\tfrac{1}{3+2\cos(t)}  \Bigl(2+2\cos(t),2\sin(t),0,1+2\cos(t)  \Bigr),\\
C_6(t)&:=\tfrac{1}{17+12\cos(t)}\Bigl(12+8\cos(t),8\sin(t),0,9+12\cos(t)\Bigr),\\
C_7(t)&:=\tfrac{1}{3-2\cos(t)}  \Bigl(2-2\cos(t),-2\sin(t),0,1-2\cos(t) \Bigr),\\
C_8(t)&:=\tfrac{1}{3+2\cos(t)}  \Bigl(2+2\cos(t),2\sin(t),0,1+2\cos(t)  \Bigr).
\end{align*}
\end{table}

From \Cref{thm:x} and its proof we recover the following four corollaries.

\begin{corollary}
\label{cor:moebius}
A Darboux cyclide is not M\"obius equivalent to both
a Bohemian surface and a Cliffordian surface.
\end{corollary}

\begin{corollary}
\label{cor:bohemian}
If $A,B\subset \R^3$ are circles \st
$A+B$ is a non-planar $\lambda$-circled surface
of M\"obius degree $d$, then $(\lambda,d)=(2,8)$.
\end{corollary}

\begin{corollary}
\label{cor:g2}
If $Z\subset S^n$ with $n\geq 3$ is a surface that contains two great circles through
a general point and is not contained in a hyperplane section,
then $n=3$ and its stereographic projection $\mu(Z)$
is either a Blum cyclide, Perseus cyclide, or ring cyclide.
\end{corollary}

See \Cref{exm:RPB} for a Blum cyclide, Perseus cyclide and ring cyclide
that contain two great circles through each point.
See \Cref{fig:4,fig:blum} for renderings of these cyclides.

\begin{corollary}
\label{cor:great}
If $Z\subset \R^3$ is a great ring cyclide, then $Z=\mu(A\star B)$
for some great circles $A,B\subset S^3$.
Great Perseus cyclides are not Cliffordian.
\end{corollary}

Notice that  a ring cyclide is by \Cref{cor:great} M\"obius equivalent
to the stereographic projection of a Clifford torus.

We conjecture the converse of Theorem~\ref{thm:x}\ref{thm:x:b}.

\begin{conjecture}
If $Z\subset \R^3$ is either a
Perseus cyclide or
CH1 cyclide,
then $Z$ is M\"obius equivalent to a Cliffordian surface.
\end{conjecture}

\subsection*{Overview}
In \SEC{model}, we propose a projective model for
M\"obius geometry which is a compactification
of $\R^3$. In \SEC{pencil}, we show that
the intersection of
Bohemian and Cliffordian surfaces
with the boundary of this compactification
consist of complex lines and/or base points of pencils of circles.
In \SEC{N}, we state a classification
of possible incidences between complex lines and base points in Darboux
cyclides.
We use these results in \SEC{c}, \SEC{b} and \SEC{g}
and obtain a list of all possible candidates
for Cliffordian, Bohemian and great Darboux cyclides, \resp.
Moreover, we show that each candidate is realized by some example.
In \SEC{proof}, we conclude the proof for \Cref{thm:x}.

\makeatletter
\renewcommand{\@cftmaketoctitle}{}
\makeatother
\vspace{2mm}
\begingroup
\def\addvspace#1{\vspace{-1mm}}
\tableofcontents
\endgroup

\section{A projective model for M\"obius geometry}
\label{sec:model}

We define a \df{real variety} $X$ to be a complex variety together with
an antiholomorphic involution $\sigma\c X\to X$ called the \df{real structure}
and we denote its real points by
\[
X_\R:=\set{p\in X}{\sigma(p)=p}.
\]
Such varieties can always be defined by polynomials with real coefficients
(see \citep[Section~I.1]{sil} and \citep[Section~6.1]{serre}).
In what follows, points, curves, surfaces and projective spaces~$\P^n$ are real algebraic varieties
and maps between such varieties are compatible with their real structures
unless explicitly stated otherwise.
In particular, curves and surfaces in this article are by default reduced and irreducible.
We assume that the real structure $\sigma\c\P^n\to\P^n$ sends $x$
to $(\overline{x_0}:\ldots:\overline{x_n})$.

Let $f\c X \dto Y\subset\P^n$ be a rational map that is not defined at~$U\subset X$.
By abuse of notation we denote $f(X\setminus U)\subseteq Y$ by $f(X)$.
We call $f$ a \df{morphism} if it is everywhere defined and thus $U=\varnothing$.

We consider the following hyperquadric $\S^3\subset\P^4$ and three different hyperplane sections
$\U,\E,\Y\subset\S^3$:
\begin{Mlist}
\item \df{M\"obius quadric}: $\S^3:=\set{x\in\P^4}{-x_0^2+x_1^2+x_2^2+x_3^2+x_4^2=0}$,
\item \df{Euclidean absolute}: $\U:=\set{x\in \S^3}{x_0-x_4=0}$,
\item \df{elliptic absolute}: $\E:=\set{x\in \S^3}{x_0=0}$, and
\item \df{hyperbolic absolute}: $\Y:=\set{x\in\S^3}{x_4=0}$.
\end{Mlist}

The following operators $\bP$, $\bC$, $\bR$ and $\bS$
are used to switch between different affine and projective models:
\begin{Mlist}
\item If $Z\subset\R^n$, then $\bP(Z)\subset\P^n$ denotes the Zariski closure of~$\iota_n(Z)$,
where the embedding
$\iota_n\c \R^n\hookrightarrow \P^n$ sends $(z_1,\ldots,z_n)$ to $(1:z_1:\ldots:z_n)$.
\item If $Z\subset\R^n$, then $\bC(Z)\subset\C^n$ denotes the Zariski closure of the
embedding of~$Z$ into~$\C^n$ via the standard embedding $\R^n\hookrightarrow\C^n$.
\item If $C\subset\P^n$, then $\bR(C)\subset\R^n$ is defined as $\iota_n^{-1}(\set{x\in C_\R}{x_0\neq 0})$,
where $\iota_n^{-1}\c\P_\R^n\dto\R^n$ sends $(x_0:\ldots:x_n)$ to $(x_1,\ldots,x_n)/x_0$.
\item If $Z\subset\R^3$, then $\bS(Z)\subset\S^3$ is defined as $\bP(\mu^{-1}(Z))$,
where $\mu\c S^3\dto\R^3$ is the stereographic projection at \SEC{intro}.
\end{Mlist}
If $\bR(\{a\})=\{b\}$ for $a\in\P^n$, then we write $\bR(a)=b$ instead.

\begin{remark}
\label{rmk:P}
We observe that, $\bP(S^3)=\S^3$, $\bR(\S^3)=S^3$, $\bS(\R^3)=\S^3$, $\bR(\E)=\E_\R=\varnothing$
and $\bR(\Y)=S^2$.
The quadrics $\S^3$, $\E$ and $\Y$ are smooth,
and
$\U$ is a quadratic cone with vertex in~$\U_\R=\{(1:0:0:0:1)\}$.
\END
\end{remark}

We consider the following linear projections from $\S^3$ to $\P^3$:
\begin{Mlist}
\item \df{stereographic projection}~$\pi\c\S^3\dto\P^3$, $\pi(x):=(x_0-x_4:x_1:x_2:x_3)$,
\item \df{central projection}~$\tau\c\S^3\to\P^3$, $\tau(x):=(x_1:x_2:x_3:x_4)$, and
\item \df{vertical projection}~$\nu\c\S^3\to\P^3$, $\nu(x):=(x_0:x_1:x_2:x_3)$.
\end{Mlist}

\begin{remark}
\label{rmk:proj}
The stereographic projection~$\pi$ corresponds via $\bP$ with $\mu\c S^3\dto \R^3$
as defined in \SEC{intro}.
The projection center of~$\pi$ lies in~$\U_\R$.
Moreover, $\pi$ defines a biregular isomorphism $\S^3\setminus\U\cong\pi(\S^3\setminus\U)$
and $\pi(\U)=\set{x\in\P^3}{x_0=x_1^2+x_2^2+x_3^2=0}$ is an irreducible conic in~$\P^3$ without real points.
The central and vertical projections define 2:1 morphisms with
ramification locus $\E$ and $\Y$, \resp.
The branching loci $\tau(\E)$ and $\nu(\Y)$ are quadrics in~$\P^3$.
The central projection $\tau\c\S^3\to\P^3$ corresponds
via $\bR$ to a 2:1 linear map~$S^3\to \R^3$ whose fibers are antipodal points.
For intuition, we remark that a linear projection~$S^1\to \R$
of the unit circle~$S^1\subset\R^2$ is via $\bP$
a 1-dimensional analogue of $\pi$, $\tau$ and $\nu$,
if the center lies either on~$S^1$, in the interior of~$S^1$, or in the exterior of~$S^1$, \resp.
\END
\end{remark}

The following two complex maps will be used for defining ``translations'' of $\S^3$:
\begin{Mlist}
\item
$\zeta_b\c \P^4\to \P^4$ with $b=(b_1,b_2,b_3)\in\C^3$
is the linear transformation corresponding to the
following $5\times 5$ matrix, where
$\Delta:=\tfrac{1}{2}(b_1^2 + b_2^2 + b_3^2)$:
\vspace{1mm}
\[
\left(
\begin{smallmatrix}
1+\Delta & b_1 & b_2 & b_3 & -\Delta \\
b_1      & 1   & 0   & 0   & -b_1\\
b_2      & 0   & 1   & 0   & -b_2\\
b_3      & 0   & 0   & 1   & -b_3\\
\Delta   & b_1 & b_2 & b_3 & 1-\Delta
\end{smallmatrix}
\right).
\]
\item
$\_\hstar\_\c\S^3\times\S^3\dto \S^3$ is the rational map defined by
\begin{align*}
(x,y)\mapsto (x_0y_0:& ~x_1y_1-x_2y_2-x_3y_3-x_4y_4:x_1y_2+x_2y_1+x_3y_4-x_4y_3:\\
                     & ~x_1y_3-x_2y_4+x_3y_1+x_4y_2:x_1y_4+x_2y_3-x_3y_2+x_4y_1).
\end{align*}
\end{Mlist}

We consider the following complex transformations of~$\S^3$,
where $\aut_\C\P^4$ denotes the complex projective transformations of~$\P^4$
and $H\in\{\U,\E,\Y\}$:
\begin{align*}
\aut_\C\S^3:=&\set{\varphi\in\aut_\C\P^4}{\varphi(\S^3)=\S^3},\\
\aut_H\S^3:=&\set{\varphi\in\aut_\C\S^3}{\varphi(H)=H},\\
\ut\S^3:=&\set{\zeta_b\c\P^4\to\P^4}{b\in\C^3},\\
\lt\S^3:=&\set{\varphi\c\S^3\dto\S^3}{\varphi(x)=p\,\hstar\,x,~ p\in\S^3\setminus\E}, \text{ and}\\
\rt\S^3:=&\set{\varphi\c\S^3\dto\S^3}{\varphi(x)=x\,\hstar\,p,~ p\in\S^3\setminus\E}.
\end{align*}
The \df{M\"obius transformations} are defined as
\[
\aut\S^3:=\set{\varphi\in\aut_\C\S^3}{\varphi\circ\sigma=\sigma\circ\varphi},
\]
where $\sigma\c\P^4\to\P^4$ denotes the real structure.
The
\df{Euclidean transformations},
\df{elliptic transformations} and
\df{hyperbolic transformations} of $\S^3$,
are defined as
\[
\aut_\U\S^3\cap\aut\S^3,\quad
\aut_\E\S^3\cap\aut\S^3\quad\text{and}\quad
\aut_\Y\S^3\cap\aut\S^3,
\quad\text{\resp.}
\]
The \df{Euclidean translations}, \df{left Clifford translations} and
\df{right Clifford translations} are defined
as
\[
\ut\S^3\cap\aut\S^3,\quad
\lt\S^3\cap\aut\S^3\quad\text{and}\quad
\rt\S^3\cap\aut\S^3,
\quad\text{\resp.}
\]
The \df{left generator} and \df{right generator}
that pass through $p\in\E$ are defined as
\[
\cL_p:=\set{q\,\hstar\,p}{q\in\S^3\setminus\E}
\quad\text{and}\quad
\cR_p:=\set{p\,\hstar\,q}{q\in\S^3\setminus\E},
\quad\text{\resp.}
\]
We shall refer to the complex lines in~$\U$ as \df{generators}.

The following proposition is classical and concerns translations in elliptic geometry
(see \citep[\textsection7.9 and 7.93]{coxeter}).
Our proof is based on \cite[Proposition~1]{gen}.
Recall from \SEC{intro} that $\_\star\_\c S^3\times S^3\to S^3$ denotes the Hamiltonian product for
the unit quaternions.

\begin{proposition}
\label{prp:E}
~
\begin{claims}
\item\label{prp:E:a}
$\bR(x\,\hstar\,y)=\bR(x)\star\bR(y)$ for all $x,y\in\S^3_\R$.
\item\label{prp:E:b}
$\lt\S^3,\rt\S^3\subset\aut_\E\S^3$.
\item\label{prp:E:c}
For all $p\in\E$,
the generators $\cL_p$ and $\cR_p$ are the two complex lines in $\E$ containing~$p$.
\item\label{prp:E:d}
For all $\varphi\in\lt\S^3$ and $p\in\E$, we have $\varphi(\cL_p)=\cL_p$.
\\For all $\varphi\in\rt\S^3$ and $p\in\E$, we have $\varphi(\cR_p)=\cR_p$.
\end{claims}
\end{proposition}

\begin{proof}
We start by introducing some terminology, which is only needed in this proof.
The algebra of \df{quaternions} consist of the
vector space $\H:=\bas{1,\qi,\qj,\qk}_\R$
together with the associative product $\_\star_{\H}\_\c \H\times \H\to \H$
that is defined by
\[
\qi\star_{\H}\qi=\qj\star_{\H}\qj=\qk\star_{\H}\qk=\qi\star_{\H}\qj\star_{\H}\qk=-1
\]
with $1\in \H$ being the multiplicative unit.
The algebra of \df{complex quaternions} is defined as $\HC:=\bas{1,\qi,\qj,\qk}_\C$
and $\ii\in\C$ denotes the imaginary unit.
The product induced by $\Hstar$ is denoted by
$\_\bullet\_\c\HC\times \HC\to \HC$.

The \df{conjugate} of a quaternion or complex quaternion $h=h_1+h_2\,\qi+h_3\,\qj+h_4\,\qk$
is defined as $h^*:=h_1-h_2\,\qi-h_3\,\qj-h_4\,\qk$.
A direct calculation shows that
\[
h\bullet h^*=h_1^2+h_2^2+h_3^2+h_4^2.
\]
We observe that $S^3=\set{h\in \H}{h\Hstar h^*=1}\subset\R^4=\H$,
and thus the Hamiltonian product $\_\star\_\c S^3\times S^3\to S^3$ is induced by $\Hstar$.
Similarly,
\[
\bC(S^3)=\set{h\in\HC}{h\bullet h^*=1}\subset\C^4=\HC.
\]
A direct calculation shows that the product
$\_\star_{S}\_\c\bC(S^3)\times\bC(S^3)\to\bC(S^3)$ induced by $\bullet$
extends to the rational map~$\_\hstar\_\c \S^3\times \S^3\dto \S^3$ defined before.
This implies that \ASN{prp:E:a} holds.

If we identify $\P^3$ with the projectivized vector space~$\P(\bas{1,\qi,\qj,\qk}_\C)$,
then $\bullet$
extends to the following rational map $\_\pstar\_\c\P^3\times\P^3\dto \P^3$,
where $x=(x_1:\ldots:x_4)$ and $y=(y_1:\ldots:y_4)$:
\begin{align*}
(x,y)\mapsto(&x_1y_1-x_2y_2-x_3y_3-x_4y_4:x_1y_2+x_2y_1+x_3y_4-x_4y_3:\\
             &x_1y_3-x_2y_4+x_3y_1+x_4y_2:x_1y_4+x_2y_3-x_3y_2+x_4y_1).
\end{align*}
We define $\iota\c \bC(S^3)\hookrightarrow\S^3$
and $\kappa\c \HC\dto \P^3$ as follows, where
the complex quaternion~$h=h_1+h_2\,\qi+h_3\,\qj+h_4\,\qk\in\HC$ is non-zero:
\[
\iota(h):=(1:h_1:h_2:h_3:h_4)
\qquad\text{and}\qquad
\kappa(h):=(h_1:h_2:h_3:h_4).
\]
Notice that $(\tau\circ\iota)(h)=\kappa(h)$,
where $\tau\c\S^3\to\P^3$ denotes the central projection.

Suppose that $q\in \bC(S^3)$
and that either
\begin{Mlist}
\item
$\psi_1(x):=q\,\sstar\,x$,\quad
$\psi_2(x):=\iota(q)\,\hstar\,x$,\quad
$\psi_3(x):=\kappa(q)\,\pstar\,x$,\quad
$\psi_4(x):=q\,\bullet\,x$, or
\item
$\psi_1(x):=x\,\sstar\,q$,\quad
$\psi_2(x):=x\,\hstar\,\iota(q)$,\quad
$\psi_3(x):=x\,\pstar\,\kappa(q)$,\quad
$\psi_4(x):=x\,\bullet\,q$.
\end{Mlist}
In both cases the diagram in Table~\ref{tab:E} commutes
as a direct consequence of the definitions.

\begin{table}[!ht]
\caption{See the proof of Proposition~\ref{prp:E}.}
\label{tab:E}
\centering
\begin{tikzpicture}[node distance=10mm and 20mm, auto]
\node (A) {$\bC(S^3)$};
\node (B) [right= of A] {$\S^3$};
\node (C) [right= of B] {$\P^3$};
\node (D) [right= of C] {$\HC$};
\node (a) [below= of A] {$\bC(S^3)$};
\node (b) [right= of a] {$\S^3$};
\node (c) [right= of b] {$\P^3$};
\node (d) [right= of c] {$\HC$};

\draw[arrows={Hooks[right]->}] (A) to node {$\iota$} (B);
\draw[->] (B) to node {$\tau$} (C);
\draw[dashed,<-] (C) to node {$\kappa$} (D);

\draw[arrows={Hooks[right]->}] (a) to node[swap] {$\iota$} (b);
\draw[->] (b) to node[swap] {$\tau$}(c);
\draw[dashed,<-] (c) to node[swap] {$\kappa$} (d);

\draw[->] (A) to node[swap] {$\psi_1$} (a);
\draw[->] (B) to node[swap] {$\psi_2$} (b);
\draw[->] (C) to node       {$\psi_3$} (c);
\draw[->] (D) to node       {$\psi_4$} (d);
\end{tikzpicture}
\end{table}

We are now ready to prove the remaining Assertions~\ref{prp:E:b}, \ref{prp:E:c} and \ref{prp:E:d}.

\ref{prp:E:b}
Since $q\,\sstar\,q^*=q^*\,\sstar\,q=1$,
we find that $\psi_1^{-1}(x)$ is equal to either $q^*\,\sstar\,x$ or $x\,\sstar\,q^*$ for all $x\in\bC(S^3)$.
It follows that $\psi_2\in\aut_\C\S^3$.
Moreover, $\iota(\bC(S^3))=\S^3\setminus\E$, which implies that $\psi_2(\E)=\E$.
Therefore, $\lt\S^3,\rt\S^3\subset\aut_\E\S^3$
and thus we concluded the proof of \ASN{prp:E:b}.

We set $E:=\set{h\in\HC}{h\,\bullet\,h^*=0}$
and for all $\alpha\in E$ we define
\[
L_\alpha:=\set{h\in\HC}{h\,\bullet\,\alpha=0}
\quad\text{and}\quad
R_\alpha:=\set{h\in\HC}{\alpha\,\bullet\,h=0}.
\]

\ref{prp:E:c}
Suppose that $\alpha\in E\setminus\{0\}$ and $\beta:=\alpha^*$.
The map $\HC\to\HC$ that sends $h$ to $h\,\bullet\,\alpha$ is linear
with respect to the underlying vector space~$\langle 1,\qi,\qj,\qk \rangle_\C$
and has kernel~$L_\alpha$.
Let $V_\beta:=\langle \beta,~\qi\bullet\beta,~\qj\bullet\beta,~\qk\bullet\beta\rangle_\C$.
By assumption,
$
\alpha\bullet\alpha^*=
\beta\bullet\beta^*=
\beta\bullet\alpha=0$ and thus $V_\beta\subseteq L_\alpha$.
Now suppose by contradiction that $\dim V_\beta<2$.
In this case,
\[
\beta=c_1\bullet\qi\bullet\beta=c_2\bullet\qj\bullet\beta=c_3\bullet\qk\bullet\beta
\]
for some $c_1,c_2,c_3\in\C$. Since
\[
\qi=\qj\,\bullet\,\qk=-\qk\,\bullet\,\qj,~~
\qj=\qk\,\bullet\,\qi=-\qi\,\bullet\,\qk,~~
\qk=\qi\,\bullet\,\qj=-\qj\,\bullet\,\qi,~~
\qi^2=\qj^2=\qk^2=-1,
\]
we find that there exists $\beta_0,\beta_1,\beta_2,\beta_3\in\C$ \st
\[
\begin{array}{rlrlrlrlr}
                \beta &=&          \beta_0 &+&         \beta_1\bullet\qi &+&         \beta_2\bullet\qj &+&         \beta_3\bullet\qk, \\
c_1\bullet\qi\bullet\beta &=& -c_1\bullet\beta_1 &+& c_1\bullet\beta_0\bullet\qi &-& c_1\bullet\beta_3\bullet\qj &+& c_1\bullet\beta_2\bullet\qk, \\
c_2\bullet\qj\bullet\beta &=& -c_2\bullet\beta_2 &+& c_2\bullet\beta_3\bullet\qi &+& c_2\bullet\beta_0\bullet\qj &-& c_2\bullet\beta_1\bullet\qk, \\
c_3\bullet\qk\bullet\beta &=& -c_3\bullet\beta_3 &-& c_3\bullet\beta_2\bullet\qi &+& c_3\bullet\beta_1\bullet\qj &+& c_3\bullet\beta_0\bullet\qk.
\end{array}
\]

By comparing the coefficients, we see that
$\beta_0=-c_i\bullet\beta_i$,
$\beta_i=c_i\bullet\beta_0$
and thus
$\beta_0\neq 0$ and
$c_i^2=-1$
for all $1\leq i\leq 3$.
This implies that $\beta_1=\beta_2=\beta_3=\pm\ii\bullet\beta_0$
and thus
\[
\beta\bullet\beta^*=\beta_0^2+\beta_1^2+\beta_2^2+\beta_3^2=-2\bullet\beta_0^2=0.
\]
We arrived at a contradiction as $\beta_0\neq 0$.
We established that
\[
\dim L_\alpha\geq \dim V_\beta\geq 2.
\]
Thus, $\kappa(L_\alpha)$
and $\kappa(R_\alpha)$ are the two complex lines in the projective quadric~$\kappa(E)$
that contain the complex point~$\kappa(\alpha^*)\in \kappa(L_\alpha)\cap\kappa(R_\alpha)$.
Since $\tau(\E)=\kappa(E)$ as a direct consequence of the definitions,
we may assume \Wlog that $\kappa(\alpha^*)=\tau(p)$.
It follows that
$\tau(\cL_p)=\kappa(L_\alpha)$
and $\tau(\cR_p)=\kappa(R_\alpha)$.
The central projection~$\tau$ is a 2:1 morphism with ramification locus~$\E$
and thus \ASN{prp:E:c} is true.

\ref{prp:E:d}
We may assume \Wlog that $\psi_2=\varphi$ and $\alpha\in E$ \st $\kappa(\alpha^*)=\tau(p)$.
Thus, $\psi_2(x)=\iota(q)\,\hstar\,x$ for some $q\in \bC(S^3)$.
If $\beta\in L_\alpha$, then $h\,\bullet\,\beta\in L_\alpha$ for all $h\in\HC$.
This implies that $\psi_4(L_\alpha)=L_\alpha$.
As \Cref{tab:E} commutes, we deduce that $\psi_2(\cL_p)=\varphi(\cL_p)=\cL_p$ as was to be shown.
The prove of the second statement
is analogous and thus we concluded the proof.
\end{proof}

The following proposition shows that the
Euclidean translations of $\S^3$ correspond to Euclidean translations of $\R^3$
and leave the generators of $\U$ invariant.

\begin{proposition}
\label{prp:U}
\hspace{1em}
\begin{claims}
\item\label{prp:U:a}
$\bR(\pi\circ\zeta_v(x))=\bR(\pi(x))+v$ for all $x\in\S^3_\R\setminus\U$ and $v\in\R^3$.
\item\label{prp:U:b}
$\ut\S^3\subset\aut_\U\S^3$.
\item\label{prp:U:c}
If $L\subset\U$ is a generator, then $\varphi(L)=L$ for all $\varphi\in\ut\S^3$.
\end{claims}
\end{proposition}

\begin{proof}
\ref{prp:U:a}
It follows from a
straightforward
\href{https://github.com/niels-lubbes/cyclides#euclidean-translations-for-s3-in-section-2}{calculation}
(see \citep[\texttt{cyclides}]{cyclides})
that
\begin{equation}
\label{eqn:U}
(\pi\circ\zeta_b)(x)=(x_{04}:x_1 + x_{04}\,b_1: x_2 + x_{04}\,b_2: x_3 + x_{04}\,b_3)
\end{equation}
for all $b\in \C^3$, where $x_{04}:=x_0-x_4$.
Thus,
$\bR(\pi\circ\zeta_b(x))=\bR(\pi(x))+b$ for all $b\in\R^3$ and $x\in\S^3_\R$ \st $x_{04}\neq 0$.

\ref{prp:U:b}
Suppose that $\varphi\in\ut\S^3$ so that $\varphi=\zeta_b$ for some $b\in\C^3$.
Let $M$ be the $5\times 5$ matrix associated to $\varphi$
and let $J$ be the diagonal matrix with $(-1,1,1,1,1)$ on its diagonal.
We verify that $M^\top\cdot J\cdot M=c\cdot J$ for some $c\in\C$
and thus $\varphi\in\aut_\C\S^3$. Since $\varphi(\U)=\U$, it follows that $\varphi\in\aut_\U\S^3$.
See \citep[\texttt{cyclides}]{cyclides}
for an automatic \href{https://github.com/niels-lubbes/cyclides#euclidean-translations-for-s3-in-section-2}{verification}.

\ref{prp:U:c}
By assumption, $\varphi=\zeta_b$ for some $b\in\C^3$.
Let $\iota(z):=(1:z_1:z_2:z_3)$,
\[
\psi_2(x):=(x_0:x_1+x_0\,b_1:x_2+x_0\,b_2:x_3+x_0\,b_3),\quad
\psi_3(z):=(z_1+b_1,z_2+b_2,z_3+b_3).
\]
It follows from \Cref{eqn:U} that the diagram in \Cref{tab:U} commutes.

\begin{table}[!ht]
\caption{See the proof of Proposition~\ref{prp:U}.}
\label{tab:U}
\centering
\begin{tikzpicture}[node distance=10mm and 20mm, auto]
\node (A) {$\S^3$};
\node (B) [right= of A] {$\P^3$};
\node (C) [right= of B] {$\C^3$};
\node (a) [below= of A] {$\S^3$};
\node (b) [right= of a] {$\P^3$};
\node (c) [right= of b] {$\C^3$};

\draw[dashed,->] (A) to node {$\pi$} (B);
\draw[arrows={<-Hooks[left]}] (B) to node {$\iota$} (C);

\draw[dashed,->] (a) to node[swap] {$\pi$} (b);
\draw[arrows={<-Hooks[left]}] (b) to node[swap] {$\iota$}(c);

\draw[->] (A) to node[swap] {$\varphi$} (a);
\draw[->] (B) to node[swap] {$\psi_2$} (b);
\draw[->] (C) to node       {$\psi_3$} (c);
\end{tikzpicture}
\end{table}

Let $H\subset \C^3$ be a complex plane.
Since $\psi_3(H)$ is parallel to $H$, we deduce
that the Zariski closures of the images~$\iota(H)$ and $(\iota\circ\psi_3)(H)$ in~$\P^3$
intersect at a complex line~$K$ at infinity.
Recall from \Cref{rmk:proj} that $\pi(\U)$
is a conic at infinity and thus $1\leq|K\cap \pi(\U)|\leq 2$ by B\'ezout's theorem.
We find that $\psi_2(K\cap \pi(\U))=K\cap \pi(\U)$.
Since $H$ was chosen arbitrary, we deduce that $\psi_2(\qq)=\qq$
for all $\qq\in\pi(\U)$.
We have $\pi(L)=\qq$ for some $\qq\in\pi(\U)$,
and thus $\varphi(L)=L$ as asserted.
\end{proof}

\begin{remark}
\label{rmk:gens}
Notice that $\E_\R=\varnothing$ and that the complex conjugate of a
left or right generator in~$\E$ is again left and right, \resp.
The generators in $\U$ are all concurrent.
Complex conjugate lines in~$\Y$ intersect in a real point.
A hyperplane section of~$\S^3$ is M\"obius equivalent to
either $\U$, $\E$ or $\Y$. Since $\U$ is unlike $\E$ and $\Y$ a tangent hyperplane section,
we could interpret Euclidean geometry as a limit case of both the hyperbolic
and elliptic geometries.
\END
\end{remark}

\begin{definition}
\label{def:model}
If $Z\subset\R^3$ is a Darboux cyclide, then
we shall call its M\"obius model~$\bS(Z)$ in $\S^3$
also a \df{Darboux cyclide}.
Similarly for CH1 cyclide, EY cyclide, Blum cyclide and so on,
and for attributes such as $\lambda$-circled, celestial, Bohemian, Cliffordian and great.
We call a conic $C\subset \S^3$ a \df{(great/small) circle}
if $\bR(C)$ is a (great/small) circle in $S^3\subset \R^4$.
\END
\end{definition}

\section{Associated pencils and absolutes}
\label{sec:pencil}

In \SEC{model}, we considered elliptic and Euclidean geometries as subgroups of
the M\"obius transformations that preserve some fixed hyperplane section
of the M\"obius quadric.
In this section, we characterize the intersection of
Cliffordian and Bohemian surfaces with such hyperplane sections.
In particular, we analyze how circles meet the elliptic or Euclidean absolute
as they move in their respective pencils.

A \df{pencil} on a surface $X\subset\S^3$ is defined as
an irreducible real hypersurface
\[
F\subset X\times \P^1
\]
\st the 1st and 2nd projections $\pi_1\c X\times \P^1\dto X$
and $\pi_2\c X\times \P^1\dto \P^1$ are dominant.
If $i\in \P^1$ is reached by $\pi_2$,
then the Zariski closure of $\pi_1(F\cap X\times\{i\})\subset X$
is called a \df{member} of $F$ and denoted by $F_i$.
We call a complex point $p\in X$ a \df{base point} of $F$,
if $p\in F_i$ for all~$i\in\pi_2(F)$.
We call $F$ a \df{pencil of conics} if
$F_i$ is a complex irreducible conic for almost all~$i\in \P^1$.
We call $F$ a \df{pencil of circles} if it is a pencil of conics \st
$F_i$ is a circle for infinitely many~$i\in \P^1_\R$.

\begin{example}
Suppose that
\[
F:=\set{(x_0 : x_1 : x_2 : x_3 ; i_0 : i_1 )}{i_0\,x_3 = i_1\,x_0,~ -x_0^2 + x_1^2 + x_2^2 + x_3^2=0}.
\]
Thus $F\subset \S^2\times \P^1$ is defined by the latitudinal circles on a sphere,
where $F_i$ is a circle if $i\in\P^1_\R$ \st $-1\leq i_1/i_0\leq 1$
and otherwise $F_i$ is a complex conic with at most one real point.
By definition, $F$ is a pencil of circles.
\END
\end{example}

\begin{remark}
If $F\subset X\times \P^1$ is a pencil of conics
on a celestial Darboux cyclide, then it follows from \citep[Theorem~9]{conical} that
$F$ is the Zariski closure of the graph of a rational map $f\c X\dto\P^1$ whose fibers are complex conics.
This implies that the first projection~$\pi_1$ is birational
and the second projection~$\pi_2$ is surjective.
The complex points where the map~$f$ is not defined correspond to the base points.
\END
\end{remark}

The following lemma is used in \Cref{lem:cc}.
We present an elementary proof based on B\'ezout's theorem to make our result accessible to a wider audience.

\begin{lemma}
\label{lem:L}
Suppose that $X\subset\S^3$ is a Darboux cyclide
and $F\subset X\times\P^1$ is a pencil of circles.
\begin{claims}

\item\label{lem:L:a}
The member $F_i$ is a complex conic
that is not contained in a complex line for all~$i\in \P^1$.

\item\label{lem:L:b}
If $L\subset X$ is a complex line \st $F_i\cap L\neq\varnothing$ for almost
all $i\in\P^1$, then $|\set{i\in\P^1}{p\in F_i}|=1$
for all complex $p\in L$.

\item\label{lem:L:c}
If $R\subset X$ is a complex line \st $F_i\cap R=\varnothing$ for almost
all $i\in\P^1$, then there exists $k\in \P^1$ \st $R\subset F_k$ is a component.

\end{claims}
\end{lemma}

\begin{proof}
Let $q:=(1:0:0:0:1)$ denote the center of the stereographic projection~$\pi\c\S^3\dto\P^3$
and recall from \Cref{rmk:proj} that $\pi$
defines a biregular isomorphism between $X\setminus \U$ and $\pi(X\setminus\U)$.
We may assume up to M\"obius equivalence that $q\in X_\R$ is a general point in~$X$.

{\bf Claim 1.} {\it
The surface $\pi(X)$ is of degree three and does not contain a complex line through a general point.
Moreover, $\pi(F_i)$ is a complex irreducible conic for almost all~$i\in\P^1$.}
\\
Since $q$ is general, it is smooth in~$X$ so that that $\deg\pi(X)=3$.
Moreover, $q$ is not a base point of~$F$ or any other pencil of complex conics.
A surface in $\S^3$ that contains infinitely many points
and a complex line through a general point,
must contain two complex conjugate lines through this point and thus be a sphere.
As the inverse stereographic projection of a complex line is either a complex
line or a complex irreducible conic that passes through $q$,
we deduce that Claim~1 holds true.

{\bf Claim 2.} {\it
There exists a circle $C\subset X$ containing~$q$ and a set~$\{H_i\}_{i\in \P^1}$ of complex hyperplane sections of~$\S^3$ \st
$H_i\cap X=F_i\cup C$ for all $i\in \P^1$.}
\\
Let $I:=\set{i\in\P^1}{\pi(F_i) \text{ is not contained in a complex line}}$
so that $|\P^1\setminus I|<\infty$ by Claim~1.
Let $P_i\subset\P^3$ denote the complex plane containing~$\pi(F_i)$ for all $i\in I$.
Recall that $\deg \pi(X)=3$ by Claim~1
and thus there exists by B\'ezout's theorem for all $i\in I$ a complex line~$M_i\subset\pi(X)$
\st $P_i\cap \pi(X)=\pi(F_i)\cup M_i$.
Since $\pi(X)$ contains by Claim~1 no continuous family of complex lines
and $F$ is a continuous family of complex curves,
we deduce that $M:=M_i=M_j$ for all $i,j\in I$.
Moreover, since $\pi(F_i)$ and the plane~$P_i$ are for all~$i\in I\cap \P^1_\R$ real and coplanar with the complex line~$M$,
it follows that $M$ must be real as well.
Let $C\subset X$ be the circle \st $\pi(C)=M$ so that $q\in C$.
Let $\{H_i\}_{i\in\P^1}$ be the complex hyperplane sections of $\S^3$ that contain the circle~$C$.
We find that $\pi(H_i)=P_i$ for all $i\in I$, which implies that $H_i\cap X=F_i\cup C$ for all $i\in I$.
Since the pencil $F$ is defined as the zero set of algebraic equations and thus Zariski closed,
we deduce that $\pi_2(F)=\P^1$ and $H_i\cap X=F_i\cup C$ for all $i\in\P^1\setminus I$.
This concludes the proof of Claim~2.

We are now ready to prove Assertions~\ref{lem:L:a}, \ref{lem:L:b} and~\ref{lem:L:c}.

\ref{lem:L:a}
Since $\deg X=4$, it follows from Claim~2 and B\'ezout's theorem that $\deg(F_i\cup C)=4$ for all~$i\in\P^1$
when counted with multiplicities. We remark that
if the hyperplane section~$H_i$ is tangent to $X$ along~$C$,
then $H_i\cap X=F_i=C$ has intersection multiplicity two.
In any case, $F_i$ is for all $i\in \P^1$ a complex curve of degree at most two.
Now suppose by contradiction that there exists $c\in\P^1$ \st
$\deg F_c=1$.
Let $\xi\c X\dto \P^3$ be the complex linear projection from a general complex point in~$F_c$.
Since $\S^3_\R$ does not contain lines, the complex conjugate line~$\oF_c$ is contained in~$X$ as well.
Notice that $F_c,\oF_c\nsubseteq\Sing X$, because otherwise
$\xi(X)$ would be a reducible complex quadric \st $\xi(\oF_c)\subset\Sing\xi(X)$.
Thus the center of $\xi$ lies in $X\setminus\Sing X$ so that $\deg\xi(X)=3$.
Since $X\cap H_c=F_c\cup C$ is a complex hyperplane section and
the complex line~$F_c$ is projected to the complex point~$\xi(F_c)\in\xi(C)$,
we find that $\xi(C\cup F_c)\subset\xi(X)$ is a complex hyperplane section consisting of the
complex irreducible conic~$\xi(C)$.
We arrived at a contradiction with B\'ezout's theorem, since a complex hyperplane section of~$\xi(X)$
must be of odd degree $\deg\xi(X)=3$ when counted with multiplicities.
This concludes the proof of \ASN{lem:L:a}.

\ref{lem:L:b}
Let $j\in \P^1$ be general so that $F_j\cap L\neq\varnothing$.
Since $H_j\cap X=F_j\cup C$ by Claim~2, it follows from B\'ezout's theorem
that $|H_j\cap L|=|F_j\cap L|=1$.
As $F$ has no base point on~$L\cap C$ by generality of $q$,
we have $F_j\cap L\nsubseteq C$ which implies that $L\cap C=\varnothing$.
Hence, for all $p\in L$ there exists a unique $i\in \P^1$ \st
$H_i$ contains~$p$.
We conclude the proof of \ASN{lem:L:b}
as $H_i\cap X=F_i\cup C$ and $L\cap C=\varnothing$, and thus~$p\in F_i$.

\ref{lem:L:c}
Let $j\in \P^1$ be general so that $F_j\cap R=\varnothing$.
Because $|H_j\cap R|=1$ by B\'ezout's theorem and $X\cap H_j=F_j\cup C$ by Claim~2,
we find that $|R\cap C|=1$.
Hence, $\pi(R\cup C)$ spans a complex plane containing the line~$\pi(C)$,
which implies that there exists $k\in\P^1$ \st $R\cup C\subset H_k$.
It follows from Claim~2 and \ASN{lem:L:a} that $R\cup C\subset F_k\cup C$ so that $R\subset F_k$.
\end{proof}

Suppose that $A,B\subset S^3$ are circles \st $A\star B\subset S^3$ is a surface
and observe that $\bP(A)\cong \P^1$.
The \df{left associated pencil} of the surface~$A\star B$
is defined as the pencil of circles
\[
F\subset\bP(A\star B)\times\bP(A)
\]
\st $F_a=\varphi_a(\bP(B))$ for all $a\in\bP(A)\setminus\E$,
where $\varphi_a\in\lt\S^3$ sends $x$ to $a\,\hstar\,x$.
It follows from \RP{E}{a} that for all $a\in\bP(A)_\R$, we have
\[
\bR(F_a)=\set{\bR(a)\,\star\,b}{b\in B}.
\]
We remark that $\lt\S^3\subset\aut_\E\S^3$ by \RP{E}{b}
and for all $a\in\bP(A)\setminus\E$, we have
$\set{a\,\hstar\,b}{b\in\bP(B)\setminus\E}\subseteq F_a$.
The member $F_a$ may be reducible for $a\in \bP(A)\cap \E$.

Similarly, the \df{right associated pencil} of~$A\star B$
is defined as the pencil of circles
\[
G\subset\bP(A\star B)\times\bS(B)
\]
\st $G_b=\varphi_b(\bP(A))$ for all $b\in\bP(B)\setminus\E$,
where $\varphi_b\in\rt\S^3$ sends $x$ to $x\,\hstar\,b$.

\begin{lemma}
\label{lem:E}
Suppose that $A\star B$ is a surface
for some circles $A,B\subset S^3$.
\begin{claims}
\item\label{lem:E:a}
If the left associated pencil $F\subset\bP(A\star B)\times\bP(A)$
has no base points on $\E$,
then $\bP(A\star B)$ contains
complex conjugate left generators $L,\oL\subset\E$
\st $|F_a\cap L|=|F_a\cap\oL|=1$ for almost all $a\in\bP(A)$.

\item\label{lem:E:b}
If the right associated pencil $G\subset\bP(A\star B)\times\bP(B)$
has no base points on $\E$,
then $\bP(A\star B)$ contains
complex conjugate right generators $R,\oR\subset\E$
\st $|G_b\cap R|=|G_b\cap\oR|=1$ for almost all $b\in\bP(B)$.
\end{claims}
\end{lemma}

\begin{proof}
\ref{lem:E:a}
For infinitely many points~$i,j\in \bP(A)_\R$ the members $F_i$ and $F_j$ are circles
and these circles are related
by a left Clifford translation.
We know from B\'ezout's theorem
that $F_i$ intersects $\E$ in two complex conjugate points,
since $\E$ is a hyperplane section of~$\S^3$.
It follows from \RP{E}{d} and \Cref{rmk:gens} that
the complex conjugate left generators $L,\oL\subset\E$
that pass through these points are left invariant.
As $F$ is Zariski closed and without base points on~$\E$ we conclude that
$L,\oL\subset\bP(A\star B)$.

The proof for \ASN{lem:E:b} is analoguous.
\end{proof}

Let $A,B\subset \R^3$ be circles and notice that $\bS(A)\cong\P^1$.
The \df{associated pencil} of $A+B$
is defined as the pencil of circles $F\subset \bS(A+B)\times \bS(A)$
such that for almost all~$a\in \bS(A)$, there exists a unique $c\in\bC(A)$
\st $F_a=\zeta_c(\bS(B))$.

As a straightforward consequence of the definitions and \RP{U}{a},
we find that $\zeta_a(b)\in\bS(A+B)$ for all $a\in A$ and $b\in\bS(B)$.
Hence, the circles
$\{C_a\}_{a\in A}$ on the surface $A+B\subset\R^3$ that are defined by $C_a:=\set{a+b}{b\in B}$
correspond via $\bR(\pi(\_))$ to a subset of $\{F_a\}_{a\in\bS(A)}$.

\begin{lemma}
\label{lem:U}
Suppose that $A,B\subset\R^3$ are circles so that $A+B$
is a surface.
If the associated pencil $F\subset \bS(A+B)\times \bS(A)$
of $A+B$ has no base points on $\U$,
then $\bS(A+B)$ contains complex conjugate generators $L,\oL\subset\U$
\st $|F_a\cap L|=|F_a\cap\oL|=1$ for almost all $a\in\bS(A)$.
\end{lemma}

\begin{proof}
Let $i\in \bS(A)$ be general.
Since $F$ has no base points on~$\U$, we may assume
\Wlog that $F_i$ does not meet the vertex of the tangent hyperplane section~$\U$.
Thus it follows from B\'ezout's theorem that
\[
|F_i\cap\U|=|F_i\cap\set{x\in\P^4}{x_0-x_4=0}|=2.
\]
Recall from \RP{U}{c} that the
Euclidean translations of $\S^3$ leave the generators of the
hyperplane section~$\U$ invariant.
Hence, there exist complex conjugate generators
$L,\oL\subset\U$ \st
$|F_a\cap L|=|F_a\cap\oL|=1$ for almost all $a\in\bS(A)$.
As $F$ is Zariski closed and without base points on~$\U$,
we conclude that $L,\oL\subset\bS(A+B)$.
\end{proof}

\section{Divisor classes of curves on Darboux cyclides}
\label{sec:N}

We recall from \cite{circle} the possible sets
of divisor classes of complex low degree curves on Darboux cyclides in~$\S^3$
(recall \Cref{def:model}).
Each entry in this classification translates into a diagram
that visualizes how complex lines, complex isolated singularities
and circles intersect.

A \df{smooth model}~$\bO(X)$ of a surface $X\subset\P^n$ is a nonsingular surface
\st there exists a birational morphism~$\varphi\c\bO(X)\to X$
that does not contract complex $(-1)$-curves.
We refer to $\varphi$ as a \df{desingularization}.

The smooth model $\bO(X)$ is unique up to biregular isomorphisms
and there exists a desingularization~$\varphi\c\bO(X)\to X$ (see \citep[Theorem~2.16]{kolsing}).

The \df{N\'eron-Severi lattice} $N(X)$ is
an additive group  defined by the divisor classes on~$\bO(X)$ up to numerical equivalence.
This group comes with an unimodular intersection product $\cdot$
and a unimodular involution $\sigma_*\c N(X)\to N(X)$ induced by the real structure $\sigma\c X\to X$.
We denote by $\aut N(X)$ the group automorphisms that are compatible
with both $\cdot$ and~$\sigma_*$.

The \df{class}~$[C]$ of a complex curve $C\subset X$ is defined as the divisor class of $\widetilde{C}$ in $N(X)$,
where $\widetilde{C}\subset \bO(X)$ is the union of complex curves
in~$\varphi^{-1}(C)$ that are not contracted to complex points
by the morphism~$\varphi$.

We consider the following subsets of $N(X)$:
\begin{Mlist}
\item $B(X)$ denotes the set of divisor classes of complex irreducible curves $C\subset \bO(X)$
\st $\varphi(C)$ is a complex point in~$X$,
\item $G(X)$ denotes the set of classes of complex irreducible conics in $X$, and
\item $E(X)$ denotes the set of classes of complex lines in $X$.
\end{Mlist}
We call $W\subset B(X)$ a \df{component} if it defines
a maximal connected subgraph of the graph with vertex set $B(X)$ and edge set $\set{\{a,b\}}{a\cdot b>0}$.
The latter subgraph is called the \df{graph of the component}.

We write $c\cdot W\succ 0$ for a subset $W\subset N(X)$, if there exists $w\in W$ \st $c\cdot w>0$.

The \df{singular locus} of $X$ is denoted by $\Sing X$.

The following proposition is an application of intersection theory on surfaces
(see \citep[Section~V.1]{har}).
For its proof we assume some background in algebraic geometry, but its assertions
are meant to be accessible to non-experts.

\begin{proposition}
\label{prp:C}
Suppose that $X\subset\S^3$ is a celestial Darboux cyclide.
Let
\begin{Mlist}
\item $\cW(X)$ be the set of components in~$B(X)$,
\item $\cF(X)$ be the set of pencils of circles on~$X$,
\item $\cG(X):=\set{g\in G(X)}{\sigma_*(g)=g}$, and
\item $\cE(X)$ be the set of complex lines in $X$.
\end{Mlist}
Suppose that
$C, C'\subset X$ are complex lines and/or circles \st $C\neq C'$.
\begin{claims}

\item\label{prp:C:a}
There exists a bijection $\Gamma\c\cW(X)\to\Sing X$
\st for all $W\in\cW(X)$ the following two properties hold:
\begin{Mlist}
\item
$\Gamma(W)\in\Sing X_\R$ \Iff $\sigma_*(W)=W$,
\item
if $[C]\cdot W\succ 0$, then $\Gamma(W)\in C$.
\end{Mlist}
The graph of a component in~$B(X)$ is a Dynkin graph of type $A_1$, $A_2$ or $A_3$.

\item\label{prp:C:b}
The map $\Lambda\c\cF(X)\to\cG(X)$ that sends $F$ to $[F_{(0:1)}]$ is a
bijection that satisfies the following two properties for all $F,G\in\cF(X)$:
\begin{Mlist}
\item $\Lambda(F)^2=0$, and
\item if $\Lambda(F)\cdot\Lambda(G)=2$, then $|F_i\cap G_j|=2$ for almost all~$i,j\in\P^1$.
\end{Mlist}
In particular, $[F_i]=[F_j]$ for all $i,j\in\P^1$.

\item\label{prp:C:c}
The map $\cE(X)\to E(X)$ that sends $L$ to $[L]$ is bijective.

\item\label{prp:C:d}
We have that $C\cap C'\neq \varnothing$
if and only if
either $[C]\cdot [C']\neq 0$, or there exists $W\in \cW(X)$
\st both $[C]\cdot W\succ0$ and $[C']\cdot W\succ0$.

\item\label{prp:C:e}
For all $p\in X$ and $F\in\cF(X)$,
the complex point~$p$ is a base point of $F$
if and only if
there exists $W\in\cW(X)$ \st $\Lambda(F)\cdot W\succ 0$ and $\Gamma(W)=p$.

\end{claims}
\end{proposition}

\begin{proof}
Let $\varphi\c \bO(X)\to X$ be a desingularization.

{\bf Claim~1.} {\it
The smooth model $\bO(X)$ is
a weak del Pezzo surface of degree four
and $X$ is its anticanonical model with at most isolated singularities.}
\\
This claim follows from \citep[Proposition~1]{conical}, where
we followed the terminology at \citep[Definition~8.1.18, Theorems~8.3.2 and~8.6.4]{dol}.

Recall that we defined the class of a complex curve~$U\subset X$ as the divisor class of the
strict transform of~$U$ in the smooth model~$\bO(X)$.
We denote the strict transforms of $C\subset X$ and $C'\subset X$
via~$\varphi$ by $D\subset \bO(X)$ and $D'\subset \bO(X)$, \resp.
If $W\in\cW(X)$, then we denote by~$C_W$ the union of complex curves in~$\bO(X)$ whose divisor class is in~$W$.

{\bf Claim~2.} {\it
The morphism~$\varphi$ restricted to~$\bO(X)\setminus\varphi^{-1}(\Sing X)$ is an isomorphism,
and $p\in \Sing X$ if and only if $p=\varphi(C_W)$ for some component~$W$
whose graph is a Dynkin graph of type $A_1$, $A_2$ or $A_3$.}
\\
It follows from Claim~1 and \citep[Proposition~8.1.10 and Theorem~8.2.27]{dol} that $\varphi$
contracts $(-2)$-curves to isolated singularities
and is an isomorphism outside the $(-2)$-curves.
Thus $B(X)$ consist of the divisor classes of the $(-2)$-curves
and $C_W$ is a union of $(-2)$-curves with classes in $W$.
We know from \citep[Theorem~8.1.11 and Theorem~8.2.28]{dol} that the graph
of the component~$W$ is a Dynkin graph.
By \citep[Corollary~5]{circle} this graph can only be of type $A_1$, $A_2$ or $A_3$.

{\bf Claim~3.} {\it
$D\cap U\neq\varnothing$ if and only if $[D]\cdot [U]>0$
for all $U\in\set{C_W}{W\in\cW(X)}\cup\{D'\}$.}
\\
Curves on rational surfaces are linearly equivalent
\Iff the curves are numerically equivalent.
Thus, it follows from \citep[Theorem~V.1.1 and Proposition~V.1.4]{har}
that $[D]\cdot [U]$ is equal to number of intersections in~$D\cap U$, when counted with multiplicity.

{\bf Claim~4.} {\it
$C\cap C'\setminus\Sing X\neq\varnothing$ if and only if $[C]\cdot [C']>0$.}
\\
Notice that $[D]\cdot[D']=[C]\cdot [C']$ by definition.
If $C\cap C'\setminus\Sing X\neq\varnothing$, then
$D\cap D'\neq\varnothing$ by Claim~2 and thus $[C]\cdot [C']>0$ by Claim~3.
If $[C]\cdot [C']>0$, then $D\cap D'\neq\varnothing$ by Claim~3
and thus $C\cap C'\setminus\Sing X\neq\varnothing$ by Claim~2.

{\bf Claim~5.} {\it
$C\cap C'\cap\Sing X\neq\varnothing$ if and only if
there exists $W\in\cW(X)$ \st
$[C]\cdot W\succ0$ and $[C']\cdot W\succ0$.}
\\
It follows from Claim~2 that
$C\cap C'\cap\Sing X\neq\varnothing$ if and only if
there exists $W\in\cW(X)$ \st
$C\cap C_W\neq\varnothing$ and $C'\cap C_W\neq \varnothing$.
By definition, $[D]\cdot [C_W]>0$ if and only if $[C]\cdot W\succ0$,
and thus Claim~5 follows from Claim~3.

{\bf Claim~6.} {\it
If $C$ and $C'$ meet transversally,
then $|C\cap C'|\geq [C]\cdot [C']$.}
\\
It follows from \citep[Theorem~V.1.1]{har} that $|D\cap D'|=[D]\cdot [D']$
and thus Claim~6 follows from Claim~2.

We are now ready to prove the
Assertions
\ref{prp:C:a},
\ref{prp:C:b},
\ref{prp:C:c},
\ref{prp:C:d} and
\ref{prp:C:e} of \Cref{prp:C}.

\ref{prp:C:a}
We define $\Gamma(W):=\varphi(C_W)$ for all $W\in\cW(X)$
and thus $\Gamma$ is well-defined and bijective by Claim~2.
Since $\sigma(C_W)=C_W$ if and only if $\sigma_*(W)=W$,
the proof for this assertion is concluded by Claims~2 and~5.

\ref{prp:C:b}
Complex curves on rational surfaces are linearly equivalent
\Iff the complex curves are numerically equivalent
and thus $[F_i]=[F_j]$ for all $i,j\in\P^1$.
We may assume \Wlog that $F_{(0:1)}$ is a circle
and thus $[F_i]\in\cG(X)$ for all $i\in\P^1$.
It follows that the map $\Lambda$ is well-defined and injective.
We know from Claim~1 and \citep[Theorem~9]{conical}
that the strict transform of a conic in~$X$
to the smooth model~$\bO(X)$ belongs to a 1-dimensional
base point free complete linear series of curves on~$\bO(X)$.
This implies that any conic in~$X$
is the member of a pencil of conics on~$X$.
Hence, if $[U]\in\cG(X)$ for some irreducible conic~$U\subset X$,
then $[U]^2=0$ and there exists a pencil~$T\subset X\times\P^1$ of conics
that has $U$ as member.
Since the members of~$T$ cover a Zariski open set of~$X_\R$,
we deduce that infinitely many members of~$T$ are circles,
which implies that~$T\in\cF(X)$.
We established that $\Lambda(T)=[U]$ so that $\Lambda$ is surjective.
Moreover, $[U]^2=0$ and thus $\Lambda(F)^2=0$.
The members $F_i$ and $G_j$ meet
transversally for general choice of~$i,j\in\P^1$.
Since complex conics in $\S^3$ intersect in at most two points,
the second property follows from Claim~6.

\ref{prp:C:c}
Complex lines $L\subset X$ do not move in a pencil and
thus are in 1:1~correspondence with their classes~$[L]\in E(X)$.

\ref{prp:C:d}
Direct consequence of Claims~4 and 5.

\ref{prp:C:e}
Let us additionally assume that $C$ and $C'$ are members of $F$.
Since $\Lambda(F)^2=[D]\cdot[D']=|D\cap D'|=0$ by \ASN{prp:C:b} and Claim~3,
the $\Rightarrow$~direction follows from Claims~2, 4 and 5.
The $\Leftarrow$~direction follows the second property at \ASN{prp:C:a}
and $\Lambda$ being well-defined.
\end{proof}

In this article,
$N(X)\cong\bas{\l_0,\l_1,\p_1,\p_2,\p_3,\p_4}_\Z$,
where the nonzero intersections between the generators are $\l_0\cdot\l_1=1$ and $\p_1^2=\p_2^2=\p_3^2=\p_4^2=-1$.
We use the following shorthand notation;
those are going to be elements in $B(X)\cup G(X)\cup E(X)$:
\begin{gather*}
\begin{array}{@{}l@{\hspace{1cm}}l@{\hspace{1cm}}l@{}}
b_1:=\p_1-\p_3, & b_{ij}:=\l_0-\p_i-\p_j   & b_0:=\l_0+\l_1-\p_1-\p_2-\p_3-\p_4, \\
b_2:=\p_2-\p_4, & b_{ij}':=\l_1-\p_i-\p_j, &
\end{array}
\\
\begin{array}{@{}l@{\hspace{1cm}}l@{\hspace{1cm}}l@{}}
g_0:=\l_0, & g_2:=2\l_0+ \l_1-\p_1-\p_2-\p_3-\p_4, & g_{ij}:=\l_0+\l_1-\p_i-\p_j,\\
g_1:=\l_1, & g_3:= \l_0+2\l_1-\p_1-\p_2-\p_3-\p_4, &
\end{array}
\\
e_i=\p_i,\qquad e_{ij}=\l_i-\p_j,\qquad e_i'=b_0+\p_i.
\end{gather*}
For convenience, we included at \cite{cyclides} a
\href{https://github.com/niels-lubbes/cyclides#initialization-of-classes-and-real-involutions-following-section-4}{table}
and
\href{https://github.com/niels-lubbes/cyclides#example-15-in-section-4}{graphs}
that
encode the pairwise intersection numbers of the above elements.

We consider the following
unimodular involutions $\sigma_*\c N(X)\to N(X)$ that are induced
by the real structure~$\sigma\c X\to X$:
\[
\begin{array}{rl}
2A_1 :& \sigma_*(\l_0)=\l_0,~ \sigma_*(\l_1)=\l_1,~ \sigma_*(\p_1)=\p_2,~ \sigma_*(\p_3)=\p_4, \\
3A_1 :& \sigma_*(\l_0)=\l_1,~ \sigma_*(\p_1)=\p_2,~ \sigma_*(\p_3)=\p_4,                      \\
D_4  :& \sigma_*(\l_0)=g_3,~ \sigma_*(\l_1)=\l_1,~ \sigma_*(\p_i)=\l_1-\p_i \text{ for } 1\leq i\leq 4. \\[2mm]
\end{array}
\]

Recall from \RP{C}{a} that the graph of a component $W\subset B(X)$
is a Dynkin graph of type $A_1$, $A_2$ or $A_3$.
The corresponding isolated double point is called a \df{node}, \df{cusp} or \df{tacnode}, \resp.
If $\sigma_*(W)\neq W$, then $\operatorname{type}(W)\in\{A_1,A_2,A_3\}$ denotes the type of~$W$.
If $\sigma_*(W)= W$, then $\operatorname{type}(W)\in\{\uAi,\uAii,\uAiii\}$ denotes the underlined type of~$W$.
If $\{W_1,\dots, W_n\}$ is the set of components in $B(X)$,
then we denote the \df{singular type} $\SingType X$ as a formal sum
$\operatorname{type}(W_1)+\ldots+\operatorname{type}(W_n)$.

A Darboux cyclide $X\subset\S^3$ is called a \df{S1 cyclide} or \df{S2 cyclide},
if $\bR(X)$ is smooth and homeomorphic to a sphere or the disjoint union of two spheres, \resp.

The following theorem follows from \citep[Theorem~4 and Corollary~5]{circle}.

\begin{theoremext}
\label{thm:B}
If $X\subset\S^3$ is a $\lambda$-circled Darboux cyclide \st $\lambda\geq 2$,
then
\[
N(X)\cong\bas{\l_0,\l_1,\p_1,\p_2,\p_3,\p_4}_\Z,
\]
and $\sigma_*$, $\SingType X$, $B(X)$, $E(X)$, $G(X)$, $\lambda$ are up to $\aut N(X)$
defined by a row in \Cref{tab:B,tab:B2}, together with the name of $X$.
\end{theoremext}

\begin{remark}
The $G(X)$, $E(X)$ and $\SingType X$ are, up to $\aut N(X)$,
uniquely determined by $B(X)$ together with $\sigma_*$ (see \citep[Theorem~4 and Corollary~5]{circle}).
Notice that the 3th, 5th, and 6th columns in \Cref{tab:B}
and the ordering/underlines in \Cref{tab:B2} are immediate corollaries of \Cref{thm:B}.
\END
\end{remark}

\begin{table}[!ht]
\centering
\caption{See \Cref{thm:B}.
A class is send by the unimodular involution $\sigma_*$ to
itself if underlined.}
\label{tab:B}
\csep{2mm}
\begin{tabular}{lclcccc}
cyclide  & $\sigma_*$ & components in $B(X)$                                                       & $\SingType X$    & $|E(X)|$ & $|G(X)|$ & $\lambda$\\\hline
Blum     & $2A_1$     & $\varnothing$                                                                & $\AN$            & $16$     & $10$     & $6$      \\
Perseus  & $2A_1$     & $\{b_1\},~\{b_2\}$                                                         & $2\Ai$           & $8$      & $7$      & $5$      \\
ring     & $2A_1$     & $\{b_{13}\},~\{b_{24}\},~\{b_{14}'\},~\{b_{23}'\}$                         & $4\Ai$           & $4$      & $4$      & $4$      \\
EH1      & $2A_1$     & $\{b_{12}\}$                                                               & $\uAi$           & $12$     & $8$      & $4$      \\
CH1      & $2A_1$     & $\{b_{13}\},~\{b_{24}\},~\{\underline{b_{12}'}\}$                          & $2\Ai+\uAi$      & $6$      & $5$      & $3$      \\
HP       & $2A_1$     & $\{\underline{b_{12}},\underline{b_{34}'}\}$                               & $\uAii$          & $8$      & $6$      & $2$      \\
EY       & $2A_1$     & $\{\underline{b_{12}},b_1,b_2\}$                                           & $\uAiii$         & $4$      & $5$      & $3$      \\
CY       & $2A_1$     & $\{\underline{b_{12}},b_1,b_2\},~\{b_{13}'\},~\{b_{24}'\}$                 & $\uAiii+2\Ai$    & $2$      & $2$      & $2$      \\
EO       & $2A_1$     & $\{\underline{b_{12}}\},~\{\underline{b_{34}}\}$                           & $2\uAi$          & $8$      & $7$      & $3$      \\
CO       & $2A_1$     & $\{\underline{b_{12}}\},~\{\underline{b_{34}}\},~\{b_{13}'\},~\{b_{24}'\}$ & $2\uAi+2\Ai$     & $4$      & $4$      & $2$      \\
\hdashline
EE/EH2   & $3A_1$     & $\{\underline{b_0}\}$                                                      & $\uAi$           & $12$     & $8$      & $2$      \\
EP       & $3A_1$     & $\{b_{13},b_{24}'\}$                                                       & $\uAii$          & $8$      & $6$      & $2$      \\
S1       & $3A_1$     & $\varnothing$                                                                & $\AN$            & $16$     & $10$     & $2$      \\
\hdashline
S2       & $D_4$      & $\varnothing$                                                                & $\AN$            & $16$     & $10$     & $2$      \\
\end{tabular}
\end{table}

\begin{table}
\centering
\caption{See \Cref{thm:B}.
A class is send by the unimodular involution $\sigma_*$ to
itself if underlined, and otherwise to
its left or right neighbor in the listing, if its position
is even or odd, \resp.
The dashed row dividers indicate when $\sigma_*$ is defined by $2A_1$, $3A_1$ or $D_4$.}
\label{tab:B2}
$
\begin{array}{ll}
\text{cyclide}   & B(X), E(X), G(X)
\\\hline
\text{Blum}      & \{\},                                                       \{e_1,e_2,e_3,e_4,e_{01},e_{02},e_{03},e_{04},e_{11},e_{12},e_{13},e_{14},e_4',e_3',e_2',e_1'\},
\\               & \{\underline{g_0},\underline{g_1},\underline{g_{12}},\underline{g_{34}},\underline{g_2},\underline{g_3},g_{13},g_{24},g_{14},g_{23}\}
\\\text{Perseus} & \{b_1,b_2\},                                                \{e_3,e_4,e_{01},e_{02},e_{11},e_{12},e_4',e_3'\},                                                  \{g_{13},g_{24},\underline{g_0},\underline{g_1},\underline{g_{12}},\underline{g_2},\underline{g_3}\}
\\\text{ring}    & \{b_{13},b_{24},b_{14}',b_{23}'\},                          \{e_1,e_2,e_3,e_4\},                                                                                \{\underline{g_0},\underline{g_1},\underline{g_{12}},\underline{g_{34}}\}
\\\text{EH1}     & \{\underline{b_{12}}\},                                     \{e_1,e_2,e_3,e_4,e_{03},e_{04},e_{11},e_{12},e_{13},e_{14},e_2',e_1'\},
\\               & \{\underline{g_0},\underline{g_1},\underline{g_{34}},\underline{g_3},g_{13},g_{24},g_{14},g_{23}\}
\\\text{CH1}     & \{b_{13},b_{24},\underline{b_{12}'}\},                      \{e_1,e_2,e_3,e_4,e_{13},e_{14}\},                                                                  \{g_{14},g_{23},\underline{g_0},\underline{g_1},\underline{g_{34}}\}
\\\text{HP}      & \{\underline{b_{12}},\underline{b_{34}'}\},                 \{e_1,e_2,e_3,e_4,e_{03},e_{04},e_{11},e_{12}\},                                                    \{\underline{g_0},\underline{g_1},g_{13},g_{24},g_{14},g_{23}\}
\\\text{EY}      & \{\underline{b_{12}},b_1,b_2\},                             \{e_3,e_4,e_{11},e_{12}\},                                                                          \{g_{13},g_{24},\underline{g_0},\underline{g_1},\underline{g_3}\}
\\\text{CY}      & \{\underline{b_{12}},b_1,b_2,b_{13}',b_{24}'\},             \{e_3,e_4\},                                                                                        \{\underline{g_0},\underline{g_1}\}
\\\text{EO}      & \{\underline{b_{12}},\underline{b_{34}}\},                  \{e_1,e_2,e_3,e_4,e_{11},e_{12},e_{13},e_{14}\},                                                    \{g_{13},g_{24},g_{14},g_{23},\underline{g_0},\underline{g_1},\underline{g_3}\}
\\\text{CO}      & \{\underline{b_{12}},\underline{b_{34}},b_{13}',b_{24}'\},  \{e_1,e_2,e_3,e_4\},                                                                                \{\underline{g_0},\underline{g_1},g_{14},g_{23}\}
\\\hdashline
\text{EE/EH2}    & \{\underline{b_0}\},                                        \{e_1,e_2,e_3,e_4,e_{01},e_{12},e_{02},e_{11},e_{03},e_{14},e_{04},e_{13}\},
\\               & \{g_0,g_1,\underline{g_{12}},\underline{g_{34}},g_{13},g_{24},g_{14},g_{23}\}
\\\text{EP}      & \{b_{13},b_{24}'\},                                         \{e_1,e_2,e_3,e_4,e_{02},e_{11},e_{04},e_{13}\},                                                    \{g_0,g_1,\underline{g_{12}},\underline{g_{34}},g_{14},g_{23}\}
\\\text{S1}      & \{\},                                                       \{e_1,e_2,e_3,e_4,e_{01},e_{12},e_{02},e_{11},e_{03},e_{14},e_{04},e_{13},e_4',e_3',e_1',e_2'\},
\\               & \{g_0,g_1,\underline{g_{12}},\underline{g_{34}},g_{13},g_{24},g_{14},g_{23},g_2,g_3\}
\\\hdashline
\text{S2}        & \{\},                                                       \{e_1,e_{11},e_2,e_{12},e_3,e_{13},e_4,e_{14},e_{01},e_1',e_{02},e_2',e_{03},e_3',e_{04},e_4'\},
\\               & \{g_0,g_3,g_{12},g_{34},g_{13},g_{24},g_{14},g_{23},\underline{g_1},\underline{g_2}\}
\end{array}
$
\end{table}

\clearpage

\begin{figure}[!ht]
\centering
\csep{2cm}
\begin{tabular}{cc}
\begin{tikzpicture}[scale=0.7] 
\draw[red] (-1.5,0) node[left, red] {$e_{11}$} -- (1.5,0);
\draw[red] (0,1.5) node[above, red] {$e_{12}$} -- (0,-1.5);
\draw[red] (-1,1) node[left, red] {$e_4$} -- (1,-1);
\draw[red] (-1,-1) node[left, red] {$e_3$} -- (1,1);
\draw[draw=black, fill=green!5]  (0,0) circle [radius=5mm] node[black] {$b_{34}$};
\node at (0,-3) {EY cyclide};
\end{tikzpicture}
&
\begin{tikzpicture}[scale=0.7] 
\draw[red] (-1,1) node[left, red] {$e_4$} -- (2,-2);
\draw[red] (-1,-1) node[left, red] {$e_3$} -- (2,2);
\draw[draw=black, fill=green!5]  (0,0) circle [radius=5mm] node[black] {$b_{12}$};
\draw[draw=black, fill=white, densely dotted]  (1,1) circle  [radius=5mm] node[black] {$b_{13}'$};
\draw[draw=black, fill=white, densely dotted]  (1,-1) circle [radius=5mm] node[black] {$b_{24}'$};
\node at (0,-3) {CY cyclide};
\end{tikzpicture}
\end{tabular}
\caption{Incidences between complex lines and isolated singularities on a Darboux cyclide~$X$ (see \Cref{exm:Y}).
Each complex line is represented as a line segment and labeled with its corresponding class in $E(X)$.
A real or non-real isolated singularity is represented as a disc with a solid and dashed
border, \resp.
Each singularity is labeled with the sum of classes in the corresponding component in $B(X)$.
}
\label{fig:Y}
\end{figure}

\begin{example}
\label{exm:Y}
The diagrams in \Cref{fig:Y,fig:pers,fig:O} show
the incidences between complex lines and isolated singularities
in EY cyclide, CY cyclide, Perseus cyclide, CH1 cyclide, ring cyclide, CO cyclide, and EO cyclide.
In case the Darboux cyclides are inversions of quadratic surfaces, it is straightforward to
compute such incidences via elementary computations. For example, an EY cyclide~$X\subset\S^3$
is for some $\lambda>0$, M\"obius equivalent to
\[
X':=\set{x\in\P^4}{x_1^2+\lambda^2\,x_2^2=(x_0-x_4)^2,~ x_1^2+x_2^2+x_3^2+x_4^2=x_0^2}.
\]
If $L\subset X'$ is a complex line, then it must be contained in $\U=\set{x\in\S^3}{x_0=x_4}$
and thus there exist $\alpha,\beta\in\{1,-1\}$ \st
\[
L=\set{x\in\P^4}{x_0=x_4,~ x_1=\alpha\,\ii\,\lambda\, x_2=\tfrac{\beta\,\lambda}{\sqrt{1-\lambda^2}}x_3}
.
\]
For Darboux cyclides $X\subset\S^3$ such that $X_\R$ is smooth this approach is less feasible.
Instead, we obtained the diagrams by applying \Cref{prp:C} to \Cref{thm:B} as follows.
Since the approach is the same for each case, we
suppose again that $X\subset \S^3$ is an EY cyclide.
We apply \Cref{thm:B} and find that
$\SingType X=\underline{A_3}$,
$B(X)=\{b_1,b_2,\underline{b_{12}}\}$ and $E(X)=\{e_{11},e_{12},e_3,e_4\}$
(see \Cref{tab:B,tab:B2}).
A component $W\subset B(X)$ corresponds to an isolated singularity
by \RP{C}{a}
and a class in $E(X)$ corresponds to a complex line in~$X$ by \RP{C}{c}.
If $\sigma_*(W)=W$, then the isolated singularity is real.
Each line segment in the diagram for the EY cyclide in~\Cref{fig:Y}
represents a complex line in $X$.
Two line segments intersect at a disc if and only if the
corresponding complex lines meet at a complex point~$p\in X$.
If $p$ is real or non-real, then the disc has solid and dashed border, \resp.
If $p$ is an isolated singularity of~$X$, then the disc is labeled
with the sum of the classes in the corresponding component $W\subset B(X)$.
In the EY cyclide case, $W=\{b_1,b_2,b_{12}\}$ and thus the label equals $b_1+b_2+b_{12}=b_{34}$.
We use \RPS{C}{a}{C}{d} to determine whether complex lines
and/or isolated singularities intersect. If $\sigma_*(a)=b$ and $a\cdot b=1$
for some $a,b\in E(X)$, then the corresponding complex conjugate lines
meet in a real point.
The diagrams for the remaining cyclides are obtained analogously,
and are automatically
\href{https://github.com/niels-lubbes/cyclides#example-15-in-section-4}{verified}
at \citep[\texttt{cyclides}]{cyclides}.
Notice that a non-real singularity
in an CY cyclide meets only one complex line.
\END
\end{example}

\begin{figure}[!ht]
\centering
\csep{0mm}
\begin{tabular}{ccc}
\begin{tikzpicture}[scale=0.4] 
\node at (0,6) {};\node at (0,-6) {};
\draw[red] (-3,2) -- (3.5,2) node[right] {$e_{12}$};
\draw[red] (-3.5,-2) node[left] {$e_{11}$} -- (3,-2);
\draw[red] (-2,3) -- (-2,-3.5) node[below] {$e_{01}$};
\draw[red] (2,3.5) node[above] {$e_{02}$} -- (2,-3);
\draw[red] (-5,4.3) node[left,red]  {$e_4$} -- (3,1.6);
\draw[red] (-4.3,5) node[above,red] {$e_3'$} -- (-1.6,-3);
\draw[red] (5,-4.3) node[right,red] {$e_3$} -- (-3,-1.6);
\draw[red] (4.3,-5) node[below,red] {$e_4'$} -- (1.6,3);
\draw[draw=black, fill=white, densely dotted] (-2,2) circle [radius=3mm];
\draw[draw=black, fill=white, densely dotted] (2,-2) circle [radius=3mm];
\draw[draw=black, fill=white, densely dotted] (-4,4) circle [radius=3mm];
\draw[draw=black, fill=white, densely dotted] (4,-4) circle [radius=3mm];
\draw[draw=black, fill=white, densely dotted] (-2,-2) circle [radius=6mm] node[black] {$b_1$};
\draw[draw=black, fill=white, densely dotted] (2,2)   circle [radius=6mm] node[black] {$b_2$};
\node at (-1,-5.5) {Perseus cyclide};
\end{tikzpicture}
&
\begin{tikzpicture}[scale=0.4] 
\node at (0,6) {};\node at (0,-6) {};
\draw[red] (-3,2)  -- (3.5,2) node[right] {$e_3$};
\draw[red] (-3.5,-2) node[left] {$e_4$} -- (3,-2);
\draw[red] (-2,3) -- (-2,-3.5) node[below] {$e_{13}$};
\draw[red] (2,3.5) node[above] {$e_{14}$} -- (2,-3);
\draw[red] (-5,4.3) node[left, red] {$e_1$} -- (3,1.6);
\draw[red] (-4.3,5) node[above, red] {$e_2$} -- (-1.6,-3);
\draw[draw=black, fill=white, densely dotted] (-2,2) circle [radius=3mm];
\draw[draw=black, fill=white, densely dotted] (2,-2) circle [radius=3mm];
\draw[draw=black, fill=white, densely dotted] (2,2) circle [radius=6mm] node[black] {\tiny$b_{13}$};
\draw[draw=black, fill=white, densely dotted] (-2,-2) circle [radius=6mm] node[black] {\tiny$b_{24}$};
\draw[draw=black, fill=green!5] (-4,4) circle [radius=6mm] node[black] {\tiny$b_{12}'$};
\node at (0,-5.5) {CH1 cyclide};
\end{tikzpicture}
&
\begin{tikzpicture}[scale=0.4] 
\node at (0,6) {};\node at (0,-6) {};
\draw[red] (-3,2)  -- (3.5,2) node[right] {$e_1$};
\draw[red] (-3.5,-2) node[left] {$e_2$} -- (3,-2);
\draw[red] (-2,3) -- (-2,-3.5) node[below] {$e_3$};
\draw[red] (2,4) node[above] {$e_4$} -- (2,-3);
\draw[draw=black, fill=white, densely dotted] (-2,2)  circle [radius=6mm] node[black] {\tiny$b_{13}$};
\draw[draw=black, fill=white, densely dotted] (2,2)   circle [radius=6mm] node[black] {\tiny$b_{14}'$};
\draw[draw=black, fill=white, densely dotted] (-2,-2) circle [radius=6mm] node[black] {\tiny$b_{23}'$};
\draw[draw=black, fill=white, densely dotted] (2,-2)  circle [radius=6mm] node[black] {\tiny$b_{24}$};
\node at (0,-5.5) {ring cyclide};
\end{tikzpicture}
\end{tabular}
\caption{Incidences between complex lines and isolated singularities (see the caption of \Cref{fig:Y}).}
\label{fig:pers}
\end{figure}

\begin{example}
\label{exm:vil}
Suppose that $Z\subset \R^3$ is a smooth torus of revolution.
The astronomer Yvon Villarceau observed
that $Z$ contains through a general point
a latitudinal circle, a longitudinal circle and two cospherical \df{Villarceau circles} \citep[1848]{vil}.
Let us identify the classes of these circles
and identify those circles that are members of a pencil with base points.
For this purpose, we suppose that
$\Gamma\c\cW(X)\to\Sing X$ and $\Lambda\c\cF(X)\to\cG(X)$
are the bijections defined at \Cref{prp:C}, where $X:=\bS(Z)$.
Let $F,G\in\cF(X)$ be the two pencils of Villarceau circles,
and let $F',G'\in\cF(X)$ be the pencils of latitudinal and longitudinal circles, \resp.
Since $X_\R$ is smooth and $X$ is covered by four pencils of circles,
it follows from \RPS{C}{a}{C}{b} that
$\sigma_*(W)\neq W$ for all $W\in\cW(X)$
and $|\cG(X)|=4$.
Thus \Cref{thm:B} implies
that the corresponding type entries in the $\SingType X$ column of \Cref{tab:B} are not underlined
and the corresponding value in the $\lambda$ columns is equal to $4$.
We find that $X$ is a ring cyclide so that
$\cG(X)=\{g_0,g_1,g_{12},g_{34}\}$ up to $\aut N(X)$ (see the corresponding row in \Cref{tab:B2}).
Notice that $\alpha\cdot\beta=2$ for $\alpha,\beta\in\cG(X)$
if and only if $\{\alpha,\beta\}=\{g_{12},g_{34}\}$.
It now follows from \RP{C}{b} that
\Wlog $\bigl(\Lambda(F),\Lambda(G),\Lambda(F'),\Lambda(G')\bigr)=(g_{12},g_{34},g_0,g_1)$.
Since
$g_0\cdot \{b_{14}'\}\succ 0$
and
$g_1\cdot \{b_{13} \}\succ 0$,
we deduce from \RPS{C}{a}{C}{e}
and that $F'$ and $G'$ each have complex conjugate base points in $\Sing X$.
\END
\end{example}

The following \Cref{lem:bp,lem:hyp} are needed in \SEC{c}.

\begin{lemma}
\label{lem:bp}
Suppose that $X\subset\S^3$ is a celestial Darboux cyclide.
\begin{claims}
\item\label{lem:bp:a}
$X$ is covered by a pencil of circles with complex conjugate base points
if and only if
$X$ is either a Perseus cyclide, ring cyclide, CH1 cyclide, CY cyclide or CO~cyclide.
\item\label{lem:bp:b}
If $X$ is covered by at least two pencils of circles
with complex conjugate base points, then $X$ is a ring cyclide.
\item\label{lem:bp:c}
If $X$ is a CY cyclide or CO cyclide,
then $X\cap \E$ does not contain complex lines.
\end{claims}
\end{lemma}

\begin{proof}
Suppose that $\Gamma\c\cW(X)\to \Sing X$
and $\Lambda\c\cF(X)\to\cG(X)$ are the bijections defined at~\Cref{prp:C}
and let
$\cM(X):=\set{W\in\cW(X)}{\sigma_*(W)\neq W }$.

{\bf Claim~1.} {\it
Each of the following items holds true up to $\aut N(X)$:
\begin{Mlist}
\item
If $X$ is a Perseus cyclide, then
$\cG(X)= \{g_0,g_1,g_{12},g_2,g_3\}$,
\[
\cM(X)=\{\{b_1\},\{b_2\}\},\qquad
g_{12}\cdot\{b_1\}\succ 0,\qquad
g_{12}\cdot \{b_2\}\succ 0,
\]
and $g\cdot b=0$ for all $g\in\{g_0,g_1,g_2,g_3\}$ and $b\in\{b_1,b_2\}$.

\item
If $X$ is a ring cyclide, then
$\cG(X)= \{g_0,g_1,g_{12},g_{34}\}$,
\begin{gather*}
\cM(X)=\{\{b_{13}\},\{b_{24}\},\{b'_{14}\},\{b'_{23}\}\},
\\
g_0\cdot\{b'_{14}\}\succ 0,\qquad
g_0\cdot \{b'_{23}\}\succ 0,\qquad
g_1\cdot\{b_{13}\}\succ 0,\qquad
g_1\cdot \{b_{24}\}\succ 0,
\end{gather*}
and $g\cdot b=0$ for all $g\in\{g_{12},g_{34}\}$ and $b\in\{b_{13},b_{24},b'_{14},b'_{23}\}$.

\item
If $X$ is a CH1 cyclide, then
$\cG(X)= \{g_0,g_1,g_{34}\}$,
\[
\cM(X)=\{\{b_{13}\},\{b_{24}\}\},\qquad
g_1\cdot\{b_{13}\}\succ 0,\qquad
g_1\cdot \{b_{24}\}\succ 0,
\]
and
$g\cdot b=0$ for all $g\in\{g_0,g_{34}\}$ and $b\in\{b_{13},b_{24}\}$.

\item
If $X$ is a CY cyclide or CO cyclide, then
$\cG(X)= \{g_0,g_1\}$,
\[
\cM(X)=\{\{b'_{13}\},\{b'_{24}\}\},\qquad
g_0\cdot\{b'_{13}\}\succ 0,\qquad
g_0\cdot \{b'_{24}\}\succ 0,
\]
and
$g_1\cdot b'_{13}=g_1\cdot b'_{24}=0$.
\end{Mlist}
}
This claim follows from \Cref{thm:B}
(see \Cref{tab:B,tab:B2} for $\cM(X)$ and $\cG(X)$, \resp).
See \citep[\texttt{cyclides}]{cyclides} for a
\href{https://github.com/niels-lubbes/cyclides#initialization-of-classes-and-real-involutions-following-section-4}{table}
with entries~$g\cdot b$
for all $g\in G(X)$ and $b\in B(X)$.

{\bf Claim~2} {\it
The pencil $F\in\cF(X)$ has complex conjugate base points
if and only if
$\Lambda(F)\cdot W\succ 0$ for some $W\in\cM(X)$.}
\\
This claim follows from \RPS{C}{a}{C}{b}.

{\bf Claim~3} {\it If $|\Sing X|-|\Sing X_\R|>0$,
then $X$ is either a Perseus cyclide, ring cyclide, CH1 cyclide, CY cyclide or CO~cyclide.}
\\
This claim follows from \Cref{thm:B} (see the $\SingType X$ column in \Cref{tab:B}).

\ref{lem:bp:a}
It follows from \RPS{C}{a}{C}{e}
that a non-real base point of a pencil~$F\in\cF(X)$ is contained in $\Sing X$.
Hence, the $\Rightarrow$~direction follows from Claim~3.
The $\Leftarrow$~direction follows from \RP{C}{b} together with Claims~1 and~2.

\ref{lem:bp:b}
Suppose that the pencil $F\in\cF(X)$ has complex conjugate base points.
It follows from Claims~1 and~2
that if $X$ is a Perseus~cyclide, CH1~cyclide, CY~cyclide or CO~cyclide,
then
$\Lambda(F)$ equals $g_{12}$, $g_1$, $g_0$ and~$g_0$, \resp.
Moreover, if $X$ is a ring cyclide, then $\Lambda(F)\in\{g_0,g_1\}$ (see also \Cref{exm:vil}).
Hence, this assertion follows from \ASN{lem:bp:a} and the injectivity of~$\Lambda$.

\ref{lem:bp:c}
Recall from \Cref{exm:Y} that the incidences between
complex lines and real isolated singularities in
a CY cyclide and CO cyclide are illustrated in the diagrams of \Cref{fig:Y} (right)
and \Cref{fig:O} (left), \resp.
We observe that each line segment in these two diagrams meet a disc with solid border.
Thus each complex line in $X$ meets some real isolated singularity.
Since $\E_\R=\varnothing$, we conclude that
$X\cap \E$ does not contain complex lines.
\end{proof}

\begin{lemma}
\label{lem:hyp}
Suppose that $X\subset\S^3$ is either
a Perseus cyclide, CH1 cyclide or ring cyclide,
and suppose that $C\subset X$ is a hyperplane section.
\begin{claims}

\item\label{lem:hyp:a}
If
$X$ is a ring cyclide and $\Sing X\subset C$,
then $C$ consists of four complex lines.

\item\label{lem:hyp:b}
If $R,\oR\subset C$ are complex conjugate lines \st $|R\cap\oR|=0$,
then there exist complex conjugate lines $L$ and $\oL$
\st $C=L\cup\oL\cup R\cup\oR$
and $|L\cap R|=|L\cap \oR|=|\oL\cap R|=|\oL\cap\oR|=1$.

\end{claims}
\end{lemma}

\begin{proof}
\ref{lem:hyp:a}
Recall from \Cref{exm:Y} that the incidences of complex lines and isolated singularities
in~$X$ are depicted in the rightmost diagram of \Cref{fig:pers}.
By assumption, the singular points are contained in the hyperplane section~$C$.
A complex line in~$X$ that meets a hyperplane in $\P^4$ in more than one complex point
must be contained in this hyperplane.
Hence, $C$ consists by B\'ezout's theorem
of four complex lines.

\ref{lem:hyp:b}
By \RP{C}{c}, two complex lines in $X$ are complex conjugate if
and only if their classes in $E(X)$ are related via $\sigma_*$.
By assumption, $X$ is either a Perseus cyclide, CH1 cyclide or ring cyclide.
Hence, $\sigma_*$ is of type $2A_1$ by \Cref{thm:B} (see the second column of \Cref{tab:B}).
In particular,
$\sigma_*(e_1)=e_2$,
$\sigma_*(e_3)=e_4$,
$\sigma_*(e_3')=e_4'$,
$\sigma_*(e_{01})=e_{02}$,
$\sigma_*(e_{11})=e_{12}$ and
$\sigma_*(e_{13})=e_{14}$.

First, suppose that $X$ is a Perseus cyclide.
Let us consider the incidences between the complex lines in $X$
as depicted in the leftmost diagram of \Cref{fig:pers}.
Notice that $R$ and $\oR$ are presented by line segments that are labeled with
$[R]$ and $\sigma_*([R])$, \resp.
Since $|R\cap\oR|=0$, we observe that $([R],[\oR])$ is equal to either
$(e_3,e_4)$,
$(e_3',e_4')$,
$(e_{11},e_{12})$, or
$(e_{01},e_{02})$.
If $([R],[\oR])=(e_3,e_4)$, then
the complex conjugate lines $L$ and $\oL$ \st $([L],[\oL])=(e_3',e_4')$,
each meet $R\cup \oR$ (and thus the hyperplane~$C$) in two complex points.
A complex line in~$X$ that meets a hyperplane in $\P^4$ in more than one complex point
must be contained in this hyperplane and thus~$L,\oL\subset C$.
Therefore,
$C=L\cup\oL\cup R\cup\oR$ by B\'ezout's theorem
and
$|L\cap R|=|L\cap \oR|=|\oL\cap R|=|\oL\cap\oR|=1$.
The remaining three cases for $([R],[\oR])$ are symmetric,
and thus \ASN{lem:hyp:b} holds for the Perseus cyclide.

If $X$ is a CH1 cyclide or ring cyclide, then
\ASN{lem:hyp:b} is shown analogously using the corresponding diagrams in \Cref{fig:pers}, except that
$([R],[L])\in\{(e_3,e_{13}),(e_{13},e_3)\}$
and
$([R],[L])\in\{(e_1,e_3),(e_3,e_1)\}$, \resp, where $[\oR]=\sigma_*([R])$ and $[\oL]=\sigma_*([L])$.
In particular, if $X$ is a CH1 cyclide,
then $([R],[\oR])\neq(e_1,e_2)$, since $|R\cap\oR|=0$ and the
line segments labeled with $e_1$ and $e_2$ in the middle diagram of~\Cref{fig:pers}
represent complex conjugate lines that meet at an isolated singularity.
\end{proof}

\section{Cliffordian Darboux cyclides}
\label{sec:c}

We develop a necessary condition for a Darboux cyclide $X$ to be Cliffordian
in terms of the sets $B(X)$, $E(X)$ and $G(X)$ in \Cref{tab:B2}.

Suppose that $X\subset \S^3$ is a Darboux cyclide. For $a,b\in N(X)$,
we set $a \odot b:=1$ if either $a\cdot b=1$ or if $a\neq b$ and there exists a component $W\subset B(X)$
\st both $a\cdot W\succ 0$ and $b\cdot W\succ 0$; in all other cases we set $a \odot b:=0$.
Notice that if $L,L'\subset X$ are different complex lines, then $|L\cap L'|=[L]\odot [L']$
by \RP{C}{d} and thus the operator~$\odot$ provides an algebraic criterion for complex lines to intersect.

A \df{Clifford quartet} is defined as a subset $\{a,b,c,d\}\subset E(X)$ of cardinality four \st
$\sigma_*(a)=b$,
$\sigma_*(c)=d$,
$a\odot b=c\odot d=0$ and
$a\odot c=c\odot b=b\odot d=d\odot a=1$.

\begin{example}
\label{exm:cq}
If $X$ is a ring cyclide, then $E(X)=\{e_1,e_2,e_3,e_4\}$ by \Cref{thm:B}.
Since
$|E(X)|=4$,
$\sigma_*(e_1)=e_2$,
$\sigma_*(e_3)=e_4$,
$e_1\odot e_2=e_3\odot e_4=0$ and
$e_1\odot e_3=e_3\odot e_2=e_2\odot e_4=1$,
we find that $E(X)$ forms a Clifford quartet.
Recall from \Cref{exm:Y} that
the diagram for the ring cyclide in \Cref{fig:pers}
represents each class~$a\in E(X)$ in terms of a line segment,
and $a\odot b=1$ for $a,b\in E(X)$
if and only if the two
corresponding line segments in the diagram meet at a disc.
Hence, we can use such diagrams, together with the specification
of $\sigma_*$, to recognize Clifford quartets.
\END
\end{example}

Recall from \Cref{rmk:P} that $\bP(S^3)=\S^3$ and thus $\bP(A\star B)\subset\S^3$ for all $A,B\subset S^3$.

\begin{lemma}
\label{lem:cq}
If $\bP(A\star B)$ is a Darboux cyclide for some circles $A,B\in S^3$,
then $\bP(A\star B)\cap\E$ consists of two left generators and two right generators
whose classes form a Clifford quartet.
\end{lemma}

\begin{proof}
Let $F\subset\bP(A\star B)\times\bP(A)$
and $G\subset\bP(A\star B)\times\bP(B)$ be the left and right associated pencils
of $A\star B$, \resp.
We set $X:=\bP(A\star B)$.

If neither $F$ nor $G$ has base points in $\E$,
then it follows from \Cref{lem:E} and B\'ezout's theorem
that $X\cap\E=\set{x\in X}{x_0=0}$ consist of two left generators
and two right generators. By \RP{C}{d}, the classes of these
generators form a Clifford quartet and thus the proof
is concluded for this case.

In the remainder of the proof we assume that $F$
has a base point in~$\E$.
Since $\E_\R=\varnothing$, we find that $F$ has two complex
conjugate base points in~$\E$.
Recall from \RPS{C}{a}{C}{e} that each base point corresponds
to a complex isolated singularity of $X$.

First suppose that $G$ has base points in $\E$ as well.
These base points must be complex conjugate
and thus~$X$ is a ring cyclide by \RL{bp}{b}.
It follows from \RL{hyp}{a} that the hyperplane section~$\E\cap X$ consists of
four complex lines.
Recall from \Cref{exm:cq} that $E(X)$ defines a Clifford quartet
and thus we concluded the proof.

Finally, suppose that $G$ does not have base points in~$\E$.
In this case, the hyperplane section~$X\cap\E$
contains two right generators $R$ and $\oR$ by \Cref{lem:E}.
We apply \RLS{bp}{a}{bp}{c} and find that $X$ is either
a Perseus cyclide, ring cyclide or CH1 cyclide.
The main assertion now follows from \RL{hyp}{b}.
\end{proof}

Now we introduce an algebraic necessary condition for a Darboux cyclide to be
Cliffordian.

\begin{definition}
\label{def:cc}
Suppose that $X\subset \S^3$ is a celestial Darboux cyclide.
We call $(A,a,g,U)$ a \df{Clifford data} if
\begin{Mlist}
\item $A=\{a,b,c,d\}$ is a Clifford quartet for $X$ with distinguished element~$a$,
\item $g\in G(X)$ \st $\sigma_*(g)=g$ and $g\cdot a\neq 0$, and
\item $U=\set{e\in E(X)}{e \cdot a=1 \text{ and } e \odot b=e \odot c=e \odot d=0}$.
\end{Mlist}
We call the Clifford data~$(A,a,g,U)$ a \df{certificate} if $g\cdot u\neq 0$ for all $u\in U$.
We say that $X$ satisfies the \df{Clifford criterion}
if there exist at least one certificate.
\END
\end{definition}

\begin{example}
\label{exm:cc}
It follows from \Cref{thm:B} that $(\{e_1,e_2,e_3,e_4\},e_1,g_{12},\varnothing)$ is a
certificate for a ring cyclide~$X\subset\S^3$,
and thus a ring cyclide satifies the Clifford criterion.
Similarly, a Perseus cyclide and CH1 cyclide satisfy
the Clifford criterion with certificates
$
(\{e_{01},e_{02},e_{11},e_{12}\},e_{01},g_1,\varnothing)
$
and
$
(\{e_3,e_4,e_{13},e_{14}\},e_3, g_{34},\varnothing),
$
\resp.
In contrast, the Clifford data
$(\{e_1,e_2,e_3',e_4'\},e_1,g_{12},\{e_{01},e_{11},e_2'\})$ for a Blum cyclide
is not a certificate, since $g_{12}\cdot e_{11}=0$ (see \Cref{fig:cc}).
In fact, the Blum cyclide does not satisfy the Clifford criterion.
See \citep[\texttt{cyclides}]{cyclides} for a
\href{https://github.com/niels-lubbes/cyclides#example-22-and-proposition-25-in-section-5}{software implementation}
that computes for each case in \Cref{tab:B2} all possible Clifford data,
and checks whether they are certificates.
For the Blum cyclide see alternatively \Cref{fig:greatblum}.
\END
\end{example}

\begin{figure}[!ht]
\centering
\begin{tikzpicture}[scale=0.5] 
\draw[red] (-3, 1) node[left]  {$e_1$}  -- ( 9, 1);
\draw[red] (-3,-1) node[left]  {$e_2$}  -- ( 9,-1);
\draw[red] (-1, 2) node[above] {$e_3'$} -- (-1,-2);
\draw[red] ( 1, 2) node[above] {$e_4'$} -- ( 1,-2);
\draw[blue] ( 3,  2) node[above] {$e_{01}$} -- ( 4, 0);
\draw[blue] ( 5,  2) node[above] {$e_{11}$} -- ( 6, 0);
\draw[blue] ( 7,  2) node[above] {$e_2'$}   -- ( 8, 0);
\draw[draw=black, fill=white, densely dotted] (3.5, 1)  circle [radius=2mm];
\draw[draw=black, fill=white, densely dotted] (5.5, 1)  circle [radius=2mm];
\draw[draw=black, fill=white, densely dotted] (7.5, 1)  circle [radius=2mm];
\draw[draw=black, fill=white, densely dotted] (1,1)   circle [radius=2mm];
\draw[draw=black, fill=white, densely dotted] (-1,-1) circle [radius=2mm];
\draw[draw=black, fill=white, densely dotted] (-1,1)  circle [radius=2mm];
\draw[draw=black, fill=white, densely dotted] (1,-1)  circle [radius=2mm];
\end{tikzpicture}
\caption{Incidences between 7 of the 16 lines in a Blum cyclide.}
\label{fig:cc}
\end{figure}

\begin{remark}
\label{rmk:clifford}
The following lemma can be seen as a generalization of \citep[7.94]{coxeter}.
Donald Coxeter refers to \citep[Chapter~X]{klein} and Felix Klein attributes
these insights to William Kingdon Clifford (1845--1879~CE).
Clifford passed away at an early age and
his theories in elliptic geometry were only partially published.
Klein saw it as a duty to workout and disseminate these theories \citep[page 238]{klein}.
To my mind, Clifford taught us that
by combining most elementary curves
we gain essential insights into the geometry of space.
\END
\end{remark}

\begin{lemma}
\label{lem:cc}
If $X=\bP(A\star B)\subset\S^3$ is a Darboux cyclide for some circles $A,B\in S^3$,
then $X$ satisfies the Clifford criterion.
\end{lemma}

\begin{proof}
Let $F\subset\bP(A\star B)\times\bP(A)$
and $G\subset\bP(A\star B)\times\bP(B)$ be the left and right associated pencils
of $A\star B$, \resp.

First suppose that $F$ or $G$ has base points on $\E$.
These base points must be complex conjugate as $\E_\R=\varnothing$.
Lemmas~\ref{lem:cq} and~\ref{lem:bp}\ref{lem:bp:c} imply
that $X$ is not a CY cyclide or CO cyclide.
Thus, it follows from \RL{bp}{a} that $X$ is either
a Perseus cyclide, ring cyclide or CH1 cyclide
and the main assertion
holds for these three cases by \Cref{exm:cc}.

In the remainder of the proof we assume that neither $F$
nor $G$ has base points on~$\E$.

Notice that $X\cap\E=\set{x\in X}{x_0=0}$ defines a hyperplane section.
Hence, the intersection $\bP(A)\cap\E$
consists by B\'ezout's theorem of the complex conjugate points~$\{\aa,\oa\}$.
We obtain for all complex $\alpha\in \bP(A)\setminus\{\aa,\oa\}$
the complex left Clifford translation
$\varphi_\alpha\in\lt\S^3$
\st
$\varphi_\alpha(x)=\alpha\,\hstar\,x$
for all $x\in \S^3\setminus \E$.
By definition, $\varphi_\alpha(\bP(B))$
is a member of the pencil $F$ for all $\alpha\in\bP(A)\setminus\{\aa,\oa\}$.
We know from \RP{E}{b} that $\varphi_\alpha\in\aut_\E\S^3$ and thus
$\varphi_\alpha(\bP(B))$ is an irreducible complex conic for all $\alpha\in\bP(A)\setminus\{\aa,\oa\}$.

We know from \Cref{lem:cq} that
$X\cap\E$ consists of two left generators
$L,\oL\subset\E$ and two right generators $R,\oR\subset\E$
intersecting in four complex points $\pp,\op,\qq,\oq\in\E$ (see \Cref{fig:LR}).
It follows from \RLS{L}{b}{E}{a} that $|\set{i\in \P^1}{p\in F_i}|=1$ for all $p\in L\cup\oL$.
Since $\E\subset\S^3$ is a hyperplane section and $|F_i\cap (L\cup\oL)|=2$,
we deduce from B\'ezout's theorem that $F_i\cap R=\varnothing$ for almost all $i\in \P^1$.
Thus, we know from \RL{L}{c} that the unique member of $F$ that contains $\pp\in L$ is a complex reducible conic
with $R$ as component.
Similarly, the unique member of $F$ that contains $\qq\in L$ has $\oR$ as component.
Since the complex left Clifford translations of~$\bP(B)$ are irreducible, it follows that
$\varphi_\alpha(\bP(B))\cap\{\pp,\qq\}=\varnothing$ for all $\alpha\in\bP(A)\setminus\{\aa,\oa\}$.
We established the following complex continuous map:
\[
\xi\c \bP(A)\setminus\{\aa,\oa\} \to L\setminus\{\pp,\qq\},
\quad
\alpha \mapsto \varphi_\alpha\bigl(\bP(B)\bigr)\cap L.
\]
In fact, $\xi$ is an complex isomorphism as each point on $L$ is reached by exactly
one member of $F$ by \RL{L}{b}.
In other words, the complex left Clifford translations of the circle $\bP(B)$ trace out $L\setminus\{\pp,\qq\}$.

\begin{figure}[!ht]
\centering
\begin{tikzpicture}[scale=0.35]
\draw[very thick, blue, densely dashed] (-8, 2) node[left] {$R$}   -- (10, 2);
\draw[very thick, blue, densely dashed] (-8,-2) node[left] {$\oR$} -- (10,-2);
\draw[very thick, red, densely dashed] (-2, 4) -- (-2,-4) node[below] {$L$};
\draw[very thick, red, densely dashed] ( 2, 4) -- ( 2,-4) node[below] {$\oL$};
\draw[draw=black, fill=white]  (-2, 2) circle [radius=0.9] node {$\pp$};
\draw[draw=black, fill=white]  ( 2, 2) circle [radius=0.9] node {$\oq$};
\draw[draw=black, fill=white]  (-2,-2) circle [radius=0.9] node {$\qq$};
\draw[draw=black, fill=white]  ( 2,-2) circle [radius=0.9] node {$\op$};
\draw[very thick, blue] (-10,-0.5) to [out=5, in=175] (10,-0.5) node[right] {$\varphi_\alpha(\bP(B))$};
\draw[draw=black, fill=orange]  ( -2, 0) circle [radius=0.4];
\draw[draw=black, fill=white]  (  2, 0) circle [radius=0.4];
\end{tikzpicture}
\caption{The incidences between $L$, $\oL$, $R$, $\oR$
and $\varphi_\alpha(\bP(B))$ for some $\alpha\in \bP(A)\setminus\{\aa,\oa\}$.}
\label{fig:LR}
\end{figure}

Suppose that $U\subset E(X)$ is the set of classes
of complex lines $M\subset X$ \st $[M]\cdot [L]=1$
and $[M]\odot [\oL]=[M]\odot [R]=[M]\odot [\oR]=0$.
By \RP{C}{d}, we have $M\cap L=\{m\}$ and $m\in L\setminus\{\pp,\qq\}$.
We established that there exists $\alpha\in\bP(A)$
\st $m\in F_\alpha$ and thus $F_\alpha\cap M\neq\varnothing$.
Suppose by contradiction that $[F_\alpha]\cdot[M]=0$.
By \RP{C}{d}, there exists a component~$W\subset B(X)$
\st $[F_\alpha]\cdot W\succ 0$ and $[M]\cdot W\succ 0$.
It follows from \RPS{C}{a}{C}{e} that $m$ is a base point of~$F$.
We arrived at a contradiction as $F$ does not have base points on~$\E$.
Therefore, we require that $[F_\alpha]\cdot [M]\neq 0$ for all~$[M]\in U$.
We conclude from Propositions~\ref{prp:C}\ref{prp:C:b}, \ref{prp:C}\ref{prp:C:c} and \ref{prp:C}\ref{prp:C:d} that
$\bigl(\left\{[L],[\oL],[R],[\oR]\right\},[L],[F_\alpha],U\bigr)$ is a certificate
and thus $X$ satisfies the Clifford criterion.
\end{proof}

\begin{proposition}
\label{prp:c}
A Cliffordian Darboux cyclide $X\subset\S^3$
is either a Perseus cyclide, ring cyclide or CH1 cyclide.
\end{proposition}

\begin{proof}
The Darboux cyclide~$X$ satisfies the Clifford criterion by \Cref{lem:cc}.
We apply \Cref{thm:B} and consider the 14 triples $\bigl(B(X),E(X),G(X)\bigr)$ in \Cref{tab:B2}.
For each such triple we go through all possible Clifford quartets in~$E(X)$.
For each such Clifford quartet~$A$ we consider all possible Clifford data~$(A,a,g,U)$.
We verify that only a Perseus cyclide, ring cyclide or CH1 cyclide
admits a Clifford data~$(A,a,g,U)$ that is a certificate.
We used \citep[\texttt{cyclides}]{cyclides} to do the
\href{https://github.com/niels-lubbes/cyclides#example-22-and-proposition-25-in-section-5}{verification}
automatically.
In particular, we find that a Clifford quartet exists only if $X$
is a
Blum cyclide,
Perseus cyclide,
ring cyclide,
EH1 cyclide,
CH1 cyclide,
HP cyclide, or
S1 cyclide.
\end{proof}

\begin{example}
\label{exm:cert}
We show that the surfaces $Z_{01}$, $Z_{23}$ and $Z_{45}$
defined at \SEC{intro} are a ring cyclide, Perseus cyclide and CH1 cyclide, \resp.
Moreover, we show that $Z_{06}$ and $Z_{78}$ are Cliffordian surfaces of degree 8.
The required computations are done
\href{https://github.com/niels-lubbes/cyclides#example-26-in-section-5}{automatically}
at \citep[\texttt{cyclides}]{cyclides}.
See \citep[\texttt{orbital}]{orbital} for an alternative implementation of these methods.
Suppose that $0\leq i\leq 8$.
Let $M_i:=${\tt Mi} be the corresponding $5\times 5$ matrix in \Cref{tab:M}.
\begin{table}[!ht]
\caption{$5\times 5$ matrices that represent elements in $\aut \S^3$.}
\centering
\label{tab:M}
\vspace{-5mm}
\begin{Verbatim}[fontsize=\scriptsize,baselinestretch=1.1,commandchars=\\\{\}]
 M0 = [( 1, 0,0, 0, 0), ( 0, 1,0,0,  0), (0,0, 1, 0,0), ( 0,0,0,1, 0), (0, 0,0, 0, 1)]
 M1 = [( 5, 0,0, 0, 0), ( 0, 5,0,0,  0), (0,0, 4,-3,0), ( 0,0,3,4, 0), (0, 0,0, 0, 5)]
 M2 = [( 3, 0,0,-2,-2), ( 0, 1,0,0,  0), (0,0, 1, 0,0), (-2,0,0,1, 2), (2, 0,0,-2,-1)]
 M3 = [( 3, 0,0, 2,-1), ( 0, 2,0,0,  0), (0,0, 2, 0,0), ( 2,0,0,2,-2), (1, 0,0, 2, 1)]
 M4 = [( 3,-2,0, 0,-1), (-2, 2,0,0,  2), (0,0, 2, 0,0), ( 0,0,0,2, 0), (1,-2,0, 0, 1)]
 M5 = [( 3, 2,0, 0,-1), ( 2, 2,0,0, -2), (0,0, 2, 0,0), ( 0,0,0,2, 0), (1, 2,0, 0, 1)]
 M6 = [(17,12,0, 0,-9), (12, 8,0,0,-12), (0,0, 8, 0,0), ( 0,0,0,8, 0), (9,12,0, 0,-1)]
 M7 = [( 3,-2,0, 0,-1), ( 2,-2,0,0, -2), (0,0,-2, 0,0), ( 0,0,0,2, 0), (1,-2,0, 0, 1)]
 M8 = [( 3, 2,0, 0,-1), ( 2, 2,0,0, -2), (0,0, 2, 0,0), ( 0,0,0,2, 0), (1, 2,0, 0, 1)]
\end{Verbatim}
\vspace{-5mm}
\end{table}

Let $J$ be the diagonal matrix with $(-1,1,1,1,1)$ on its diagonal.
We verify that
there exists $\lambda\in\Q\setminus\{0\}$
\st $M_i^\top\cdot J\cdot M_i=\lambda\, J$,
and thus $M_i$ defines a M\"obius transformation~$\varphi_i\c\S^3\to\S^3$.
The curve parametrization~$C_i(t)$ in \Cref{tab:C} is related to the matrix~$M_i$
as follows, where $\psi(t):=(1:\cos(t):\sin(t):0:0)$:
\[
\bP(\set{C_i(t)}{0\leq t\leq 2\pi})=
\set{(\varphi_i\circ\psi)(t)}{0\leq t\leq 2\pi}.
\]
Since $\psi(t)$ parametrizes a great circle,
we observe that $C_i(t)$ parametrizes a circle as well.
Let $\vec{c}:=(1,0,0,0,0)$ and notice that $\bR((1:0:0:0:0))$ is the center of $S^3$.
If $i\in\{0,1\}$, then we verify that there exists $\lambda\in\Q\setminus\{0\}$ \st
$M_i\cdot\vec{c}=\lambda\,\vec{c}$.
Hence, $C_0(t)$ and $C_1(t)$ parametrize great circles.
For all $(i,j)\in\{(0,1),(2,3),(4,5),(0,6),(7,8)\}$, we implicitize the surface
\[
X_{ij}:=\bP(\set{C_i(\alpha)\star C_j(\beta)}{0\leq\alpha,\beta<2\pi})\subset\S^3
\]
and find the following (we refer to \citep[\texttt{cyclides}]{cyclides} for the
\href{https://github.com/niels-lubbes/cyclides#example-26-in-section-5}{details}):
\begin{gather*}
(\deg X_{01},\deg X_{23},\deg X_{45},\deg X_{06},\deg X_{78})=(4,4,4,8,8)
\quad\text{and}\\
(|\Sing X_{01}|,|\Sing X_{23}|,|\Sing X_{45}|)=(4,2,3).
\end{gather*}
Since $X_{ij}=\bS(Z_{ij})$ by definition,
it follows from \Cref{prp:c} and \Cref{thm:B} (see the $\SingType X$ column in \Cref{tab:B}) that
$Z_{01}$, $Z_{23}$ and $Z_{45}$ are a ring cyclide, Perseus cyclide and CH1 cyclide, \resp.
As $\deg X_{06}=\deg X_{78}=8$, we find that $Z_{06}$ and $Z_{78}$ are Cliffordian surfaces of M\"obius degree 8.
Since $C_0(t)$ parametrizes a great circle, it follows that the surfaces $Z_{01}$ and $Z_{06}$ are great.
\END
\end{example}

\section{Bohemian Darboux cyclides}
\label{sec:b}

We show that a Bohemian Darboux cyclide is the pointwise sum of a line
and a circle in $\R^3$, namely a CY or EY as in \Cref{fig:4}, and thus not
the pointwise sum of two circles in~$\R^3$.

\begin{remark}
\label{rmk:boh}
Before we apply the methods established in the previous sections, let us
sketch an alternative proof strategy for the claim that a Darboux cyclide
is not the pointwise sum of two circles.
Indeed, if $Z\subset \R^3$ is the pointwise sum of two circles,
then there exist infinitely many
parallel plane sections of~$Z$ that consist of two circles that belong to the same pencil
(see for example \Cref{fig:8}).
As we continuously vary such a plane section, the
two circles can be deformed into a single circle.
The points along which the coplanar circles intersect trace out an arc in
in the singular locus of~$Z$.
A real singularity in a Darboux cyclide must be isolated
and thus $Z$ is not a Darboux cyclide.
Depending on the background of the reader this or a similar strategy
may be more appropriate.
\END
\end{remark}

\begin{lemma}
\label{lem:b}
Suppose that $X\subset\S^3$ is a Darboux cyclide that
contains two circles through a general point that do not meet in two points.
If some pair of complex conjugate lines in~$X$ intersect,
then these complex lines meet at an isolated singularity.
\end{lemma}

\begin{proof}
Let $\cG(X):=\set{g\in G(X)}{\sigma_*(g)=g}$
and suppose that $L,\oL\subset X$ are complex conjugate lines \st $L\cap\oL\neq\varnothing$.
We know from \RP{C}{b} that there exist
different $f,f'\in \cG(X)$ \st $f\cdot f'\neq 2$.
By \RP{C}{c}, both $[L]$ and $\sigma_*([L])=[\oL]$ belong to~$E(X)$.
We apply \Cref{thm:B} and verify for each of the 14 cases in \Cref{tab:B2}
that the following statement holds:
If there exist different
$f,f'\in \cG(X)$ \st $f\cdot f'\neq 2$,
then $e\cdot\sigma_*(e)=0$ for all $e\in E(X)$.
See \citep[\texttt{cyclides}]{cyclides} for an automatic
\href{https://github.com/niels-lubbes/cyclides#lemma-28-in-section-6}{verification}
of this statement.
It follows from \RP{C}{d}
that there exists a component $W\subset B(X)$ \st $[L]\cdot W\succ0$ and $\sigma_*([L])\cdot W\succ0$.
By \RP{C}{a} such a component~$W$ corresponds to an isolated singularity in~$L\cap\oL$
and thus we concluded the proof.
\end{proof}

\begin{proposition}
\label{prp:b}
If $A,B\subset\R^3$ are generalized circles \st $\bS(A+B)$ is a Darboux cyclide,
then $A+B$ is either a CY or EY.
\end{proposition}

\begin{proof}
We first consider the quadratic case:

{\bf Claim 1.} {\it
If $\deg(A+B)\leq 2$, then $A+B$ is either a CY or EY and
either $A$ or $B$ is a line.}
\\
Since $\deg\bS(A+B)=4$, we find that $\deg(A+B)=2$.
We go through the well-known
classification of quadratic surfaces up to
Euclidean similarity (see \citep[Proposition~4]{circle} and \Cref{fig:4})
and conclude that Claim~1 holds.

Now let us assume by contradiction that $\deg A=\deg B=2$.
Notice that $A+B$ is not a plane by Claim~1 and
thus the circles $A$ and $B$ are not coplanar.

Let $C_v:=\{v\}+B$ and $D_v:=A+\{v\}$ for $v\in \R^3$.

{\bf Claim 2.} {\it
$|C_a\cap D_b|=1$ for almost all $a\in A$ and $b\in B$.}
\\
Let $H_a\subset\R^3$ denote the spanning plane of $C_a$.
Since $a$ and $b$ are general,
there exists $q\in A\setminus\{a\}$ \st
the circle $D_b$ meets $H_a$ transversally at
the two points~$a+b\in A$ and $q+b$.
As we translate the circle $D_b$ along $C_a$,
the incidence points $a+b$ and $q+b$ trace out the coplanar circles~$C_a\subset H_a$
and~$C_q\subset H_a$, \resp.
Since $|C_a\cap C_q|<\infty$ and $b$ is general,
we deduce that $q+b\notin C_a$. It follows that $C_a\cap D_b=\{a+b\}$
and thus Claim~2 holds true.

{\bf Claim 3.} {\it
$|C_i\cap C_j|=0$ for almost all $i,j\in A$.}
\\
The spanning planes of $C_i$ and $C_j$
are parallel, but not equal.
This implies the assertion of Claim~3.

Let
$F\subset \bS(A+B)\times \bS(A)$ and $G\subset \bS(B+A)\times \bS(B)$
be the associated pencils of $A+B$ and $B+A$, \resp, where $A+B=B+A$.
It follows from Claim~2 that $|F_a\cap G_b|=1$
for almost all $a\in \bS(A)$ and $b\in \bS(B)$.
We deduce from Claim~3 that a base point of $F$ or $G$ must lie in $\U$.
We make a case distinction on the base points of $F$ and $G$.

First, we suppose that either $F$ or $G$ has no base points.
We know from \Cref{lem:U} that~$X$ contains complex conjugate
lines that meet at the vertex of $\U$.
General members of $F$ and $G$ meet in one point and thus $X$ has
an isolated singularity at this vertex by \Cref{lem:b}.
Hence, $\deg\pi(X)=\deg(A+B)=2$ with $\deg A=\deg B=2$.
We arrived at a contradiction with Claim~1.

Next, we suppose that $F$ has a real base point in $\U_\R$.
In this case $\deg\pi(X)=\deg(A+B)=2$ with $\deg A=\deg B=2$,
since $\U_\R\subset \Sing X$ by \RPS{C}{a}{C}{e}.
We arrived at a contradiction with Claim~1.

Finally, suppose by contradiction that both $F$ and $G$ have non-real base points in~$\U$.
Then $X$ must be a ring cyclide by \RL{bp}{b} and $\Sing X\subset\U$ by
\RPS{C}{a}{C}{e}.
It follows from \RL{hyp}{a}
that the hyperplane section $X\cap \U$ consists of four lines.
We arrived at a contradiction with the diagram in \Cref{fig:pers} (see \Cref{exm:Y}),
as all lines in $\U$ should be concurrent.

We arrived at a contradiction at all three cases
and thus we established either $A$ or $B$ is a line.
Therefore, $A+B$ is covered with lines, so that
either $F$ or $G$ has a base point in~$\U_\R$.
Hence, by \RPS{C}{a}{C}{e} the center of
stereographic projection~$\pi$ is in $\Sing X$ so that
$\deg(A+B)\leq 2$.
The proof is now concluded by Claim~1.
\end{proof}

\section{Great Darboux cyclides}
\label{sec:g}

In this section, we show that a Darboux cyclide is M\"obius equivalent to
a great celestial surface if and only if this cyclide is either
a Blum cyclide, Perseus cyclide, ring cyclide, EO cyclide or CO cyclide.
Moreover, we show that great ring cyclides are Cliffordian and that great Perseus cyclides
are not Cliffordian.
The hyperbolic and Euclidean analogues of great Darboux cyclides are considered as well.

\begin{remark}
\label{rmk:great}
We prove \Cref{thm:x}\ref{thm:x:c} by applying also methods from \cite{pot,tak,cool}
as formulated at \Cref{lem:cool,lem:ql,lem:h2} and \Cref{prp:hqf} below.
\Cref{rmk:class} below sketches how
these methods can be used to classify Darboux cyclides up to
M\"obius equivalence.
Since \Cref{thm:x}\ref{thm:x:c} can be read off from this classification,
this is an alternative proof strategy avoiding the use of \Cref{thm:B} and \Cref{prp:C}.
\END
\end{remark}

We call a Darboux cyclide $X\subset\S^3$
\df{elliptic} or \df{$\delta$-elliptic} for some $\delta\in\{1,2\}$,
if there exists a M\"obius transformation $\varphi\in\aut\S^3$  \st
$(\tau\circ\varphi)(X)$ is a ruled quadratic surface in~$\P^3$
that is covered by $\delta>0$ pencils of lines.
The $\delta$-\df{hyperbolic} and $\delta$-\df{Euclidean} Darboux cyclides
are defined analogously, but with the central projection~$\tau$ replaced with the vertical projection~$\nu$ and
stereographic projection~$\pi$, \resp~(recall \Cref{rmk:proj}).

Notice that an elliptic Darboux cyclide is M\"obius equivalent to
a great Darboux cyclide, since great circles are centrally projected to lines.

The following classical result is essentially \citep[Chapter~VII, Theorem~20, page~296]{cool},
but we followed the proof of \citep[Theorem~5.30 in the updated arXiv version]{sko}.

\begin{lemma}
\label{lem:cool}
If $X\subset\S^3$ is a surface \st $\deg X\neq 2$ and  $F,G\subset X\times\P^1$ are pencils of
circles \st $|F_i\cap G_j|=2$ for almost all $i,j\in\P^1_\R$,
then $X$ is a Darboux cyclide that is either elliptic, hyperbolic or Euclidean.
Moreover, if $F$ is base point free, then $X$ is either 2-elliptic or 2-hyperbolic.
\end{lemma}

\begin{proof}
Let $u,v\in\P^1_\R$ be general.
Two planes in $\P^4$ have non-empty intersection
and thus the planes spanned by the irreducible conics~$F_u$ and $F_v$ meet at a point $p\in \P^4$.
As $\deg X\neq 2$, we deduce that $p$ is the unique intersection point between these two planes.
Since $|F_u\cap G_j|=2$, the spanning planes of~$F_u$ and~$G_j$ intersect
along a line.
This implies that $p$ is the intersection point of
the two lines containing the pairs~$F_u\cap G_j$
and $F_v\cap G_j$, \resp.
We established that the spanning plane of $G_j$ contains~$p$.
Repeating the same argument with $F$ and $G$ interchanged shows that the
complex spanning plane of $F_i$ contains the point~$p$ as well.

We claim that there exists a M\"obius transformation~$\varphi\in\aut\S^3$
\st $\varphi(p)$ coincides with the projection center of either $\tau$, $\nu$
or $\pi$.
Indeed, $p\in\P^4$ corresponds via the hyperquadric~$\S^3$
uniquely to its polar hyperplane section~$H_p\subset\S^3$.
Since $H_p$ is M\"obius equivalent to either $\E$, $\Y$ or $\U$,
there exists $\varphi\in\aut\S^3$ \st $\varphi(p)$
is equal to the projection centers~$(1:0:0:0:0)$, $(0:0:0:0:1)$ and $(1:0:0:0:1)$, \resp.

First suppose that $\varphi(p)=(1:0:0:0:0)$.
In this case, the general members $F_i$ and $G_j$
are 2:1 projected to complex lines via the map $\tau\circ \varphi$.
Since $\deg X\neq 2$ by assumption, we deduce that $\deg (\tau\circ \varphi)(X)=2$
and thus $X$ is an elliptic Darboux cyclide.
If $F$ is base point free, then two general members of $F$ are disjoint.
This implies that the 2:1 projections of two general members
of~$F$ are disjoint lines in the quadric $(\tau\circ \varphi)(X)$.
Hence, the ruled quadric~$(\tau\circ \varphi)(X)$ must be smooth and thus doubly ruled.
This implies that $X$ is 2-elliptic.

If $\varphi(p)=(0:0:0:0:1)$, then $X$ must be a hyperbolic Darboux cyclide
and if $F$ is base point free, then $X$ is 2-hyperbolic.
The proof is analogous as in the elliptic case.

Finally, suppose that $\varphi(p)=(1:0:0:0:1)$.
In this case, $F$ and $G$ have a common base point at $p$.
Thus, the general members $F_i$ and $G_j$
are via the map $\pi\circ\varphi$ birationally projected to complex lines.
Since $\deg X\neq 2$, we find that $X$ must be an Euclidean Darboux cyclide.

We considered all three cases and thus the proof is concluded.
\end{proof}

We remark that a celestial Darboux cyclide
that is neither a CO cyclide nor a CY cyclide satisfies the hypothesis of \Cref{lem:cool}
by \Cref{thm:B} and \Cref{prp:C}.

Let $\aut_H\P^3:=\set{\varphi\in\aut\P^3}{\varphi(H)=H}$,
where $\aut\P^3$ denotes the real projective transformations of~$\P^3$
and $H\in\{\pi(\U),\tau(\E),\nu(\Y)\}$ (see \Cref{rmk:proj}).
Recall from \SEC{model} the related definition of~$\aut_W\S^3\cap \aut\S^3$
for $W\in\{\U,\E,\Y\}$.

\begin{remark}
\label{rmk:class}
We sketch how the methods in \cite{pot}
can recover the classification of Darboux cyclides up to $\aut\S^3$
as stated in \cite{tak}.
If $X\subset\S^3$ is a Darboux cyclide, then it follows from
\cite[Theorem~3]{pot} that that up to $\aut\S^3$ either $\pi(X)$,
$\tau(X)$ or $\nu(X)$ is a quadratic surface~$Y\subset \P^3$.
This implies that $X$ is the intersection of~$\S^3$ with a quadratic hypercone in~$\P^4$ over~$Y$.
The M\"obius transformations~$\aut_W\S^3\cap \aut\S^3$ for $W\in\{\U,\E,\Y\}$
induce the projective transformations~$\aut_H\P^3$ \st $H\in\{\pi(\U),\tau(\E),\nu(\Y)\}$, \resp.
The classification of $X$ up to $\aut\S^3$ can therefore be reduced to
the classification of~$Y\subset \P^3$ up to $\aut_H\P^3$
for $H\in\{\pi(\U),\tau(\E),\nu(\Y)\}$.
\begin{Mlist}
\item
If $Y=\pi(X)$, then $\aut_{\pi(\U)}\P^3$ is a projectivization of the Euclidean similarities of $\R^3$
and thus the classification of $Y$ is classically known (see for example \citep[Proposition~4]{circle}).
\item
If $Y=\tau(X)$, then an application of orthogonal diagonalization in linear algebra
shows that $Y$ is up to $\aut_{\tau(\E)}\P^3$ the zero set of a diagonal quadratic form.
\item
The case when $Y=\nu(\Y)$ is less known and treated by \Cref{prp:hqf} below.
Mikhail Skopenkov explained me how this proposition is an application of Uhlig's~\citep[Theorem~1]{uhlig}.
See \citep[\textsection6]{tak} by Takeuchi for a related approach.
\end{Mlist}
We can use \citep[Theorem~3]{pot}
to recover from $Y$ the number~$\lambda$ \st $X$ is $\lambda$-circled (see \Cref{lem:ql} below).
\END
\end{remark}

The \df{signature} of a real symmetric matrix~$M$ is denoted by $\sgn(M)=(p,n)$,
where $p$ and $n$ denote the number of positive and negative eigenvalues of~$M$, \resp.
Let $\diag(M_1,\ldots,M_r)$ denote the block diagonal matrix that has the square matrices $M_1,\ldots,M_r$
on its diagonal. If the square matrix $M_i$ is equal to the $1\times1$ matrix~$[m_i]$
for all $1\leq i\leq r$, then we write $\diag(m_1,\ldots,m_r)$ instead.

\begin{proposition}
\label{prp:hqf}
If $Y\subset\P^3$ is a quadratic surface, then there exist $\varphi\in\aut_{\nu(\Y)}\P^3$
and $a,b,c,d\in\R$ \st $Y$ is the zero set of
one of the following quadratic forms:
\begin{Mlist}
\item $q_1:=a\, x_0^2+b\, x_1^2 + c\, x_2^2 + d\, x_3^2$,
\item $q_2:=a\, x_0^2+ x_0\,x_1+(1-a)\, x_1^2 + b\, x_2^2 + c\, x_3^2$,
\item $q_3:=a\, x_0^2+ x_0\,x_1- a\, x_1^2 + b\, x_2^2 + c\, x_3^2$, or
\item $q_4:=a\, x_0^2 - a\, x_1^2 + x_0x_2 + x_1x_2 -a\,x_2^2 + b\,x_3^2$.
\end{Mlist}
\end{proposition}

\begin{proof}
Let $J$ and $M$ be the symmetric $4\times 4$ matrices
such that the quadratic surfaces~$\nu(\Y)$ and $Y$
are the zero set of quadratic forms~$\vec{x}^\top\cdot J\cdot \vec{x}$ and $\vec{x}^\top\cdot M\cdot \vec{x}$,
\resp, where $\vec{x}=(x_0,x_1,x_2,x_3)^\top$.
As a direct consequence of the definitions, we find that
the main assertion is equivalent to the following statement:
\textit{There exists a matrix $V$, coefficients~$a,b,c,d\in\R$ and $i\in\{1,2,3,4\}$ \st $V^\top\cdot J\cdot V=J$ and
\[
(V\cdot \vec{x})^\top\cdot M\cdot (V\cdot \vec{x})=\vec{x}^\top\cdot (V^\top\cdot M\cdot V)\cdot \vec{x}=q_i.
\]}%
The specialization of Uhlig's \citep[Theorem~1]{uhlig} to dimension~$4$
states that there exists a matrix~$R$
and $\epsilon_i\in\{1,-1\}$ and $\alpha_i,\beta_i\in\R$ for $1\leq i\leq 4$ and $1\leq r\leq 4$
\st
\[
K:=R^\top\cdot J\cdot R=\diag(\epsilon_1\cdot E_1,\ldots,\epsilon_r\cdot E_r)
~\text{and}~
N:=R^\top\cdot M\cdot R=\diag(\epsilon_1\cdot F_1,\ldots,\epsilon_r\cdot F_r),
\]
where $(E_i,F_i)$ is for all $1\leq i\leq r$ one of the following five pairs of matrices:
\[
\left(
\left[\begin{smallmatrix}\alpha_i \end{smallmatrix}\right],
\left[\begin{smallmatrix}1\end{smallmatrix}\right]
\right),\qquad
\left(
\left[\begin{smallmatrix}0&1\\1&0\end{smallmatrix}\right],
\left[\begin{smallmatrix}0&\alpha_i\\\alpha_i&1\end{smallmatrix}\right]
\right),\qquad
\left(
\left[\begin{smallmatrix}0&1\\1&0\end{smallmatrix}\right],
\left[\begin{smallmatrix}\beta_i&\alpha_i\\\alpha_i&-\beta_i\end{smallmatrix}\right]
\right),\qquad
\]\[
\left(
\left[\begin{smallmatrix}0&0&1\\0&1&0\\1&0&0\end{smallmatrix}\right],
\left[\begin{smallmatrix}0&0&\alpha_i\\0&\alpha_i&1\\\alpha_i&1&0\end{smallmatrix}\right]
\right),\qquad
\left(
\left[\begin{smallmatrix}0&0&0&1\\0&0&1&0\\0&1&0&0\\1&0&0&0\end{smallmatrix}\right],
\left[\begin{smallmatrix}0&0&0&\alpha_i\\0&0&\alpha_i&1\\0&\alpha_i&1&0\\\alpha_i&1&0&0\end{smallmatrix}\right]
\right).
\]
Moreover, we require that that if $\epsilon_i=-1$, then
\[
(E_i,F_i)\neq
\left(
\left[\begin{smallmatrix}0&1\\1&0\end{smallmatrix}\right],
\left[\begin{smallmatrix}\beta_i&\alpha_i\\\alpha_i&-\beta_i\end{smallmatrix}\right]
\right).
\]
We have $\sgn(K)=\sgn(J)$ as the signature is invariant under real congruence transformations and thus
\[
\sgn(K)=
\sgn(\epsilon_1\cdot E_1,\ldots,\epsilon_r\cdot E_r)=\sgn(\epsilon_1\cdot E_1)+\cdots+\sgn(\epsilon_r\cdot E_r)=(3,1),
\text{ where}
\]
\[
\sgn(E_i)\in\{(1,0),(1,1),(2,1),(2,2)\}
~~\text{and}~~
\sgn(-E_i)\in\{(0,1),(1,1),(1,2),(2,2)\}.
\]
It follows that the matrix pair $(K,N)$ is equal to one of the following five pairs:
\[
\left(
\left[\begin{smallmatrix}-1&0&0&0\\0&1&0&0\\0&0&1&0\\0&0&0&1\end{smallmatrix}\right],
\left[\begin{smallmatrix}\alpha_1&0&0&0\\0&\alpha_2&0&0\\0&0&\alpha_3&0\\0&0&0&\alpha_4\end{smallmatrix}\right]
\right),
\quad
\left(
\left[\begin{smallmatrix}1&0&0&0\\0&1&0&0\\0&0&0&1\\0&0&1&0\end{smallmatrix}\right],
\left[\begin{smallmatrix}\alpha_1&0&0&0\\0&\alpha_2&0&0\\0&0&0&\alpha_3\\0&0&\alpha_3&1\end{smallmatrix}\right]
\right),
\quad
\]
\[
\left(
\left[\begin{smallmatrix}1&0&0&0\\0&1&0&0\\0&0&0&1\\0&0&1&0\end{smallmatrix}\right],
\left[\begin{smallmatrix}\alpha_1&0&0&0\\0&\alpha_2&0&0\\0&0&\beta_3&\alpha_3\\0&0&\alpha_3&-\beta_3\end{smallmatrix}\right]
\right),
\quad
\left(
\left[\begin{smallmatrix}1&0&0&0\\0&1&0&0\\0&0&0&-1\\0&0&-1&0\end{smallmatrix}\right],
\left[\begin{smallmatrix}\alpha_1&0&0&0\\0&\alpha_2&0&0\\0&0&0&-\alpha_3\\0&0&-\alpha_3&-1\end{smallmatrix}\right]
\right),
\]
\[
\left(
\left[\begin{smallmatrix}1&0&0&0\\0&0&0&1\\0&0&1&0\\0&1&0&0\end{smallmatrix}\right],
\left[\begin{smallmatrix}\alpha_1&0&0&0\\0&0&0&\alpha_2\\0&0&\alpha_2&1\\0&\alpha_2&1&0\end{smallmatrix}\right]
\right).
\]
For each of these five cases for $(K,N)$, we do the following
\href{https://github.com/niels-lubbes/cyclides#proposition-33-and-lemma-35-in-section-7}{computations}
as implemented at \citep[\texttt{cyclides}]{cyclides}.
We compute the orthogonal diagonalization of $K$ and obtain an orthogonal matrix $U$ \st
$U^\top\cdot K\cdot U=J$. We verify that there exists $1\leq i\leq 4$ \st
$q_i=\vec{x}^\top\cdot (U^T\cdot N\cdot U)\cdot\vec{x}$.
It follows that if $V=R\cdot U$, then
$V^\top\cdot J\cdot V^\top=J$
and
$q_i=\vec{x}^\top\cdot (V^T\cdot M\cdot V)\cdot\vec{x}$
as was to be shown.
\end{proof}

We remark that if the ideal of~$X\subset\S^3$
is generated by the quadratic forms $q_3$ and $-x_0^2+x_1^2+x_2^2+x_3^2+x_4^2$,
then $X$ is M\"obius equivalent to the surface given by \citep[Equation~(3.2)]{tak}.

\begin{lemma}
\label{lem:ql}
Let $X:=\set{x\in\P^4}{\vec{x}^\top\cdot Q\cdot \vec{x}=\vec{x}^\top\cdot J\cdot \vec{x}=0}$,
where $J:=\diag(-1,1,1,1,1)$, $Q$ is a symmetric $5\times 5$ matrix and
$\vec{x}:=(x_0,\ldots,x_4)^\top$.
If $X\subset\S^3$ is a $\lambda$-circled Darboux cyclide
and $\set{t\in\R}{\det(Q-t\cdot J)=0}=\{t_1,\ldots,t_r\}$, then
\[
\lambda=\sum_{i=1}^r \chi\bigl(\sgn(Q-t_i\cdot J)\bigr),
\]
where the function $\chi\c \Z_{\geq 0}^2\to\{0,1,2\}$ is defined as follows:
\[
\chi(\alpha,\beta) :=
\begin{cases}
2 & \text{ if } \{\alpha,\beta\}=\{2\},\\
1 & \text{ if } \{\alpha,\beta\}=\{1,2\},\\
0 & \text{ otherwise}.
\end{cases}
\]
\end{lemma}

\begin{proof}
This assertion is a consequence of \citep[Theorem~3 and Remark~8]{pot}.
We outline the argument for convenience of the reader,
but refer to \citep[\textsection2.3 and \textsection2.4]{pot} for the details.
By assumption, $X$ is the base locus of the pencil of hyperquadrics in~$\P^4$
that are generated by the hyperquadrics~$\set{x\in\P^4}{\vec{x}^\top\cdot Q\cdot \vec{x}=0}$ and $\S^3$.
For each point $p\notin\S^3$ there exists a unique hyperquadric $H\subset \P^4$
\st $X=\S^3\cap H$ and $p\in H$ (see \citep[\textsection2.3]{pot}).
If we choose $p$ in a plane spanned by a circle~$C\subset X$, then $\{p\}\cup C \subset H$ and thus
$H$ must be a cone by B\'ezout's theorem.
This implies that $H$ is equal to the hypercone $H_i:=\set{x\in\P^4}{\vec{x}^\top\cdot(Q-t_i\cdot J)\cdot\vec{x}=0}$
for some $1\leq i\leq r$
and the 2-planes in $H_i$ intersect $X$ in conics.
Moreover, for each pencil of circles~$F$ on $X$ there exists a unique~$1\leq j\leq r$
and a unique pencil~$G$ of 2-planes on~$H_j$ \st
the circles that belong to~$F$ span infinitely many 2-planes that belong to~$G$.
As $X$ is a Darboux cyclide by assumption, we deduce that $\Sing H_i$ consists of a single vertex~$v_i$.
Suppose that $\rho_i\c\P^4\dto \P^3$ is a linear projection with center~$v_i$ and
let $\lambda_i$ denote the number of pencils of 2-planes on~$H_i$.
We deduce that $\rho_i(X)$ is a quadratic surface with signature~$\sgn(Q-t_i\cdot J)$
that is covered by exactly~$\lambda_i$ pencils of lines.
It is now straightforward to see that $\lambda_i=\chi\bigl(\sgn(Q-t_i\cdot J)\bigr)$
and thus $\lambda=\lambda_1+\cdots+\lambda_r$ as was to be shown.
\end{proof}

\begin{lemma}
\label{lem:h2}
If $X\subset\S^3$ is a $\lambda$-circled hyperbolic celestial Darboux cyclide
\st $\Sing X_\R=\varnothing$ and $X$ is not elliptic,
then $\lambda=2$.
\end{lemma}

\begin{proof}
Since $X$ is hyperbolic, we may assume \Wlog that $\nu(X)$ is a ruled quadric.
Notice that the hyperbolic transformation of~$\S^3$ commute via the vertical projection~$\nu$
with the projective transformations in $\aut_{\nu(\Y)}\P^3$.
Hence, it follows from \Cref{prp:hqf} that there exists $\varphi\in\aut\S^3$ and $a,b,c,d\in\R$
\st $(\nu\circ\varphi)(X)$ is the zero set of the quadratic form $q_i$ for some $1\leq i\leq 4$.
This implies that the ideal of $\varphi(X)$ is equal to $\bas{q_i,s}$, where
$s:=-x_0^2+x_1^2+x_2^2+x_3^2+x_4^2$.
We consider the following Jacobian matrix of the generators of the ideal $\bas{q_i,s}$:
\[
F_i(x):=
\begin{bmatrix}
\partial_0q_i&\partial_1q_i&\partial_2q_i&\partial_3q_i&\partial_4q_i\\
\partial_0s  &\partial_1s  &\partial_2s  &\partial_3s  &\partial_4s
\end{bmatrix}
=
\begin{bmatrix}
\partial_0q_i&\partial_1q_i&\partial_2q_i&\partial_3q_i&0\\
-2\,x_0&2\,x_1&2\,x_2&2\,x_3&2\,x_4
\end{bmatrix},
\]
where $\partial_j(\cdot)$ denotes the partial derivative \wrt
the variable $x_j$.
Notice that $p\in \Sing X$ if and only if $p\in X$ and the rank of the matrix~$F_i(p)$ is not equal to~$2$.
A direct
\href{https://github.com/niels-lubbes/cyclides#proposition-33-and-lemma-34-in-section-7}{computation}
shows that $F_k(x)$ has for all $k\in\{2,4\}$ rank one at the point $(1:-1:0:0:0)$.
We refer to \citep[\texttt{cyclides}]{cyclides} for the details of the matrix
\href{https://github.com/niels-lubbes/cyclides#proposition-33-and-lemma-35-in-section-7}{computations}
in this proof.
Notice that
\[
q_1+a\cdot s=(b+a)\,x_1^2+(c+a)\,x_2^2+(d+a)\,x_3^2+a\,x_4^2
\]
and thus $\bas{q_1,s}$ is the ideal of an elliptic Darboux cyclide.
Since $\Sing X_\R=\varnothing$ and $X$ is not elliptic by assumption,
we find that $i=3$ is the only remaining option.
Let $J:=\diag(-1,1,1,1,1)$ and let $Q$ be the symmetric matrix associated to~$q_3$.
In order to recover the number~$\lambda$ from the matrix~$Q$,
we apply \Cref{lem:ql}.
A direct
\href{https://github.com/niels-lubbes/cyclides#proposition-33-and-lemma-34-in-section-7}{computation}
shows that
\[
\set{t\in\R}{\det(Q-t\cdot J)=0}=\{0,b,c\}.
\]
Moreover, we
\href{https://github.com/niels-lubbes/cyclides#proposition-33-and-lemma-34-in-section-7}{find}
that the nonzero eigenvalues of $Q-t\cdot J$ for $t\in\{0,b,c\}$ are as follows,
where $u(t):=\tfrac{1}{2}\sqrt{1 + 4(a+t)^2}$:
\[
\{u(0), -u(0), b,  c\},\quad
\{u(b), -u(b),-b, -b + c\},\quad
\{u(c), -u(c),-c,  b - c\}.
\]
The next step is to use these eigenvalues to recover the following triple of signatures:
\[
\Lambda:=\bigl(\sgn(Q-0\cdot J),~\sgn(Q-b\cdot J),~\sgn(Q-c\cdot J)\bigr).
\]
We assume \Wlog that $b\geq c$ and observe that $u(t)>0$.
We make a case distinction:
\begin{Mlist}
\item If $b,c>0$ or $b,c<0$, then $\Lambda\in\{
\bigl((3,1),(1,3),(2,2)\bigr),~
\bigl((1,3),(2,2),(3,1)\bigr)
\}$.
\item If $b>0$ and $c\leq 0$, then $\Lambda\in\{
\bigl((2,2),(1,3),(3,1)\bigr),~
\bigl((2,1),(1,3),(2,1)\bigr)
\}$.
\item If $b=0$ and $c\leq0$, then $\Lambda\in\{
\bigl((1,2),(1,2),(3,1)\bigr),~
\bigl((1,1),(1,1),(1,1)\bigr)
\}$.
\item If $b=c\neq0$, then $\Lambda\in\{\bigl((3,1),(1,2),(1,2)\bigr),~\bigl((1,3),(2,1),(2,1)\bigr)\}$.
\end{Mlist}
It follows from \Cref{lem:ql}
and $|\{0,b,c\}|\leq 3$ that
$\lambda\leq\chi(\Lambda_1)+\chi(\Lambda_2)+\chi(\Lambda_3)$,
where $\Lambda=(\Lambda_1,\Lambda_2,\Lambda_3)$.
We conclude that $\lambda=2$ as was to be shown.
\end{proof}

\begin{lemma}
\label{lem:gg}
Suppose that $X\subset\S^3$ is a celestial Darboux cyclide.
\begin{claims}
\item\label{lem:gg:a}
If $X$ is elliptic,
then $|\Sing X_\R|\in\{0,2\}$ and $\Sing X_\R=(\Sing X)\setminus\E$.

\item\label{lem:gg:b}
If $X$ is 2-elliptic, then
complex conjugate lines in $X$ do not intersect.

\item\label{lem:gg:c}
The surface $X$ is 1-elliptic if and only if
$|\Sing X_\R|=2$ and
$X$ is covered by a pencil of circles with two real base points.
\end{claims}
\end{lemma}

\begin{proof}
\ref{lem:gg:a}
We may assume up to M\"obius equivalence that $\tau(X)$ is a ruled quadric.
Notice that the central projection~$\tau$
is a 2:1 covering that
defines, \wrt the complex analytic topology, locally a complex isomorphism outside
the ramification locus~$\E$.
This implies that $\tau((\Sing X)\setminus\E)\subset\Sing\tau(X)$.
If $|\Sing X_\R|=0$, then $\tau(X)$ must be smooth and thus $\Sing X\subset \E$
so that $\Sing X_\R=(\Sing X)\setminus\E=\varnothing$.
Now suppose that $|\Sing X_\R|>0$. In this case, $\tau(X)$ must be singular
and thus $(\Sing X)\setminus \E$ consist of two real antipodal points
that are send via $\tau$ to the vertex of the quadratic cone~$\tau(X)$.
It follows that $|\Sing X_\R|=2$ and $\Sing X_\R=(\Sing X)\setminus\E$ as asserted.

\ref{lem:gg:b}
We may assume up to M\"obius equivalence that $\tau(X)$ is a doubly ruled quadric.
Suppose by contradiction that
there exist complex conjugate lines $L,\oL\subset X$
that intersect at some point~$\pp$.
Complex conjugate lines in the doubly ruled quadric~$\tau(X)$ do not intersect and thus
$\tau(L)=\tau(\oL)$ so that $\pp$ is contained in the ramification locus~$\E$.
We arrived at a contradiction since $\pp\in X_\R$ and $\E_\R=\varnothing$.

\ref{lem:gg:c}
First, we show the $\Rightarrow$ direction.
We may assume \Wlog that $\tau(X)$ is a quadratic cone and thus $\Sing \tau(X)_\R=\{v\}$
for some point~$v$.
We deduce that there exist antipodal points $\pp, \qq\in \Sing X_\R$ \st $\tau(\pp)=\tau(\qq)=v$.
Moreover, there exists a pencil of circles $F\subset X\times\P^1$
with base points $\pp$ and $\qq$, and its members are 2:1 centrally projected to
lines in $\tau(X)$ that meet at~$v$. We know from \ASN{lem:gg:a} that $|\Sing X_\R|=2$.

Next, we show the $\Leftarrow$ direction.
We may assume up to M\"obius equivalence that the base points $\pp, \qq\in X_\R$ are antipodal
so that the circles in the pencil are 2:1 projected by $\tau$ to lines in $\tau(X)$
that pass through the point~$\tau(\pp)=\tau(\qq)$. Hence, $\tau(X)$ is a quadratic cone
with vertex~$\tau(\pp)$, which implies that $X$ is 1-elliptic.
\end{proof}

\begin{proposition}
\label{prp:g}
A celestial Darboux cyclide $X\subset\S^3$ is elliptic
if and only if
$X$ is either a
Blum cyclide,
Perseus cyclide,
ring cyclide,
EO cyclide, or
CO cyclide.
\end{proposition}

\begin{proof}
We first show the $\Rightarrow$ direction.

If $\Sing X_\R\neq\varnothing$,
then $|\Sing X_\R|=2$ by \RL{gg}{a}
and thus $X$ is either a EO cyclide or CO cyclide by \Cref{thm:B} (see \Cref{tab:B}).

Now suppose that $\Sing X_\R=\varnothing$.
By \Cref{thm:B} (see \Cref{tab:B}), $X$ is either a S1, S2, Blum, Perseus or ring cyclide.
If $X$ is a S1 cyclide, then we know from \Cref{thm:B} (see \Cref{tab:B2}) that
$e_{01},e_{12}\in E(X)$, $\sigma_*(e_{01})=e_{12}$ and $e_{01}\cdot e_{12}=1$.
Similarly, if $X$ is a S2 cyclide, then
$e_{1},e_{11}\in E(X)$, $\sigma_*(e_{1})=e_{11}$ and $e_{1}\cdot e_{11}=1$.
We apply \RPS{C}{c}{C}{d} and find that
S1 cyclides and S2 cyclides contain two complex conjugate lines that intersect.
Such cyclides are not elliptic by \RLS{gg}{b}{gg}{c}.
Thus if $X_\R$ is smooth, then $X$ must be either a Blum, Perseus or ring cyclide.

We proceed by showing the $\Leftarrow$ direction.

First, suppose that $X$ is either an EO cyclide or CO cyclide.
We know from \Cref{thm:B}
that $|\Sing X_\R|=2$,
$b_{12},b_{34}\in B(X)$,
$g_1\in G(X)$,
$\sigma_*(g_1)=g_1$,
$\sigma_*(b_{12})=b_{12}$,
$\sigma_*(b_{34})=b_{34}$
and $g_1\cdot b_{12}=g_1\cdot b_{34}=1$.
Hence, $X$ is by \RPS{C}{b}{C}{e} covered by a
pencil of circles with two real base points.
It follows from \RL{gg}{c} that $X$ is elliptic.

Finally, we assume that $X$ is either a Blum, Perseus or ring cyclide.
By definition, $\Sing X_\R=\varnothing$ and $X$ is $\lambda$-circled with $\lambda\geq 4$.
It follows that $\deg \pi(X)>2$ and thus $X$ is not Euclidean.
Hence, $X$ is by \Cref{lem:cool} either hyperbolic, elliptic or both.
We conclude from \Cref{lem:h2} that $X$ must be elliptic.
\end{proof}

\begin{remark}
\label{rmk:ring}
I do not know how to show that all the Blum, Perseus and ring cyclides are elliptic
without using \Cref{prp:hqf} and \Cref{lem:ql,lem:h2}.
We include an alternative proof for the ring cyclide case
as it provides additional geometric insight and does not rely
on this proposition and two lemmas.

\textit{A ring cyclide $X\subset\S^3$ is 2-elliptic and not 2-hyperbolic.}

This proof goes as follows.
By \Cref{thm:B} (see \Cref{tab:B2})
there exist
\[
g_{12},g_{34}\in \set{g\in G(X)}{\sigma_*(g)=g}
\]
\st
$g_{12}\cdot g_{34}=2$ and $g_{12}\cdot b=0$
for all $b\in B(X)$, where $B(X)=\{b_{13},b_{24},b'_{14},b'_{23}\}$.
Hence, we know from \RPS{C}{b}{C}{e}
that there exist pencils $F,G\subset X\times\P^1$ \st
$|F_i\cap G_j|=2$ for general $i,j\in\P^1$ and $F$ is base point free.
It follows from \Cref{lem:cool} that $X$ is either 2-elliptic, 2-hyperbolic or both.

Now suppose by contradiction that $X$ is 2-hyperbolic.
We may assume \Wlog that $\nu(X)\subset\P^3$ is a doubly ruled quadric.
The restriction of the 2:1 covering~$\nu$ to $X$
defines outside the ramification locus~$\Y\subset\S^3$ a local complex analytic isomorphism
on each of the two sheets.
Since $\nu(X)$ is smooth, it follows that $\Sing X\subset \Y$.
Recall from \Cref{exm:Y} (see the rightmost diagram of \Cref{fig:pers})
that there exist skew complex lines $L,\oL\subset X$ \st $[L]=e_1$, $[\oL]=e_2$,
$|L\cap\oL|=0$ and $|L\cap \Sing X|=|\oL\cap \Sing X|=2$.
Since $\Y\cap X$ is a hyperplane section of~$X$ and $|L\cap\Y|=|\oL\cap\Y|=2$,
it follows from B\'ezout's theorem that $L,\oL\subset\Y$.
Notice that $L$ and $\oL$ are complex conjugate lines as $\sigma_*([L])=[\oL]$.
We arrived at a contradiction, because
$\Y_\R\cong S^2$ and thus complex conjugate lines in $\Y$ are not disjoint (see \Cref{rmk:gens}).
This concludes the alternative proof.
\END
\end{remark}

\begin{proposition}
\label{prp:ring}
If $X\subset\S^3$ is a great ring cyclide, then
there exist great circles $A,B\subset S^3$ \st $X=\bP(A\star B)$.
\end{proposition}

\begin{proof}
We fix a point $e\in X_\R$.
First, suppose that $\bR(e)$ equals the identity quaternion in~$S^3$.

We know from \Cref{thm:B} that $|\Sing X_\R|=0$ and $|\Sing X|=4$.
It follows from \RL{gg}{c} and \RL{gg}{a} that
$\tau(X)$ is a doubly ruled quadric and $\Sing X\subset\E$.
We apply \RL{hyp}{a} and find that $X\cap\E$ consist of two left generators
$L$, $\oL$ and two right generators~$R$, $\oR$.
As $\tau(X)$ is doubly ruled, there exist
two great circles $A,B\subset \bR(X)$ \st $e\in \bP(A)\cap\bP(B)$.
We may assume \Wlog that the line $\tau(\bP(A))$ does not belong
to the pencil of lines containing $\tau(R)$ and~$\tau(\oR)$.
Hence, each circle in the pencil containing $\bP(A)$
meets both right generators $R,\oR\subset \E$.

We assume by contradiction that $X\neq \bP(A\star B)$.
Let $F\subset \bP(A\star B)\times \bP(B)$ and $G\subset \bP(A\star B)\times \bP(A)$
be the right and left associated pencils of $A\star B$.
Notice that $\bR(c)$ with $c:=(1:0:0:0:0)$ is the center of $S^3$
and that $p\,\hat{\star}\,c=c\,\hat{\star}\,p$ for all~$p\in\S^3$.
It follows that the left or right Clifford translation of a great circle in~$\S^3$ is again great.
Thus, both $F$ and~$G$ have infinitely many members that are great circles
on~$\bP(A\star B)$.
These great circles are centrally projected to lines on
the doubly ruled quadric~$\tau(\bP(A\star B))$.
This implies that both $F$ and $G$ are base point free.
It now follows from \RL{E}{b} that $|F_b\cap R'|=|F_b\cap\overline{R'}|=1$
for almost all~$b\in B$ and some complex conjugate right generators $R',\overline{R'}\subset \E$.
Since $F_e=\bP(A)$ and $|\bP(A)\cap R|=|\bP(A)\cap\oR|=1$,
we deduce from \RL{L}{b} that $R'=R$ and $\overline{R'}=\oR$.
Notice that $F_e=\bP(A)$ and $G_e=\bP(B)$
and thus $\bP(A),\bP(B)\subset X\cap\bP(A\star B)$.
We fix some general point~$b\in B$.
Let $C\subset X$ be the great circle that passes through~$b$
and belongs to the same pencil on~$X$ as $\bP(A)$.
Recall that $|C\cap R|=|C\cap\oR|=1$ as is illustrated in~\Cref{fig:converse}.
By assumption, $F$ does not cover $X$
and thus $C$ is not a member of $F$.
We arrived at a contradiction as the lines~$\tau(F_b)$ and~$\tau(C)$ span a plane
so that the complex lines~$\tau(R)$ and $\tau(\oR)$ cannot be skew.
We establised that $X=\bP(A\star B)$ for great circles $A,B\subset \bR(X)$.

\begin{figure}[!ht]
\centering
\begin{tikzpicture}[scale=0.8]
%
\draw[thick, red] (-2,1) -- (2,1) node[right] {$R$};
\draw[thick, red] (-2,-1) -- (2,-1) node[right] {$\oR$};
\draw[blue] (-1,2) to [out=260, in=100] (-1,-2) node[below] {$\bP(A)$};
\draw[blue] (1,2) to [out=260, in=100] (1,-2) node[below] {$F_b$};
\draw[densely dotted, blue] (0,2) to [out=290, in=125] (2,-2) node[below right] {$C$};
\draw[black!20!green] (-3,-0.2) to [out=10, in=170] (3,-0.2) node[right] {$\bP(B)$};
\draw[fill=white] (0.8,0.05) circle [radius=0.8mm] node[black,above right] {$b$};
\draw[fill=white] (-1.2,0.05) circle [radius=0.8mm] node[black,above right] {$e$};
\draw[fill=white] (-1.15,1) circle   [radius=0.8mm];
\draw[fill=white] (-1.15,-1) circle  [radius=0.8mm];
\draw[fill=white] (1-0.15,1) circle  [radius=0.8mm];
\draw[fill=white] (1-0.15,-1) circle [radius=0.8mm];
\draw[fill=white] (0.35,1) circle    [radius=0.8mm];
\draw[fill=white] (1.35,-1) circle   [radius=0.8mm];
\end{tikzpicture}
\caption{See the proof of \Cref{prp:ring}.
The incidences between the great circles $\bP(A)$, $\bP(B)$, $C$, $F_b$ and the
right generators $R,\oR\subset\E$, under the assumption that
$\bR(e)$ is the identity quaternion in~$S^3$.}
\label{fig:converse}
\end{figure}

Finally, suppose that $\bR(e)$ is not equal to the identity quaternion in~$S^3$.
There exists a right Clifford translation $\varphi\in\rt\S^3$ \st
$\bR(\varphi(e))$ is equal to the identity quaternion.
Recall from \RP{E}{d} that $\varphi$ leaves the right generators of~$\E$ invariant.
The right Clifford translation of a great circle is again great
and thus $\varphi(X)=\bP(A\star B)$ for some great circles $A,B\subset \varphi(X)$.
The unit quaternions~$S^3$ form a group and thus
there exists $r\in S^3$ \st
$\bR(\varphi(p))=\bR(p)\star r$ for all $p\in\S^3_\R$.
Therefore, $X=\bP(A\star B')$, where $B':=\set{b\star r^{-1}}{b\in B}$.
This concludes the proof as $X=\bP(A\star B')$ for some great circles $A,B'\subset S^3$.
\end{proof}

\begin{proposition}
\label{prp:perseus}
If $X\subset \S^3$ is a great Perseus cyclide,
then $X$ is not Cliffordian,
$\tau(X)$ is a doubly ruled quadric,
$\Sing X\subset \E$ and $X\cap\E$ does not contain lines.
\end{proposition}

\begin{proof}
Recall from \Cref{exm:Y} and the leftmost diagram in~\Cref{fig:pers}
that
the incidences between all the complex lines
and complex conjugate isolated singularities
are as depicted in \Cref{fig:gp}, where
$L,\oL,R,\oR,M,\oM,T,\oT\subset X$
denote the complex lines and
$\pp,\op \in \Sing X$ are
the complex conjugate isolated singularities.
We know from \RLS{gg}{c}{gg}{a} that
$\tau(X)$ is a doubly ruled quadric and $\Sing X=\{\pp,\op\}\subset \E$.

\begin{figure}[!ht]
\vspace{-5mm}
\centering
\begin{tikzpicture}[scale=0.5] 
\node at (0,6) {};\node at (0,-6) {};
\draw[red] (-3,2)    -- (3.5,2);
\draw[magenta] (-3.5,-2) -- (3,-2);
\draw[black!30!green] (-2,3)    -- (-2,-3.5);
\draw[blue] (2,3.5)   -- (2,-3);
\draw[red] (-5,4.3)  -- (3,1.6);
\draw[black!30!green] (-4.3,5)  -- (-1.6,-3);
\draw[magenta] (5,-4.3)  -- (-3,-1.6);
\draw[blue] (4.3,-5)  -- (1.6,3);
\draw[draw=black, fill=white] (-2,2) circle [radius=3mm];
\draw[draw=black, fill=white] (2,-2) circle [radius=3mm];
\draw[draw=black, fill=white] (-4,4) circle [radius=3mm];
\draw[draw=black, fill=white] (4,-4) circle [radius=3mm];
\draw[draw=black, fill=orange!20] (-2,-2) circle [radius=6mm] node[black] {$\pp$};
\draw[draw=black, fill=orange!20] (2,2)   circle [radius=6mm] node[black] {$\op$};
\node[magenta] at (0,-1.3) {$L$};
\node[magenta] at (1,-3.8) {$M$};
\node[red] at (0,1.3) {$\oL$};
\node[red] at (-1,3.7) {$\oM$};
\node[black!30!green] at (-1.4,0) {$R$};
\node[black!30!green] at (-3.5,1) {$T$};
\node[blue] at (1.4,0) {$\oR$};
\node[blue] at (3.8,-1) {$\oT$};
\end{tikzpicture}
\vspace{-10mm}
\caption{Incidences between complex conjugate lines and isolated singularities in a great Perseus cyclide,
where $\pp,\op\in \E$.}
\label{fig:gp}
\end{figure}

We claim that none of the complex lines in $X$ are contained in $\E$.
First suppose by contradiction that $L\subset \E$.
In this
case the complex conjugate line $\oL$ must also be contained in $\E$.
Therefore, $L,\oL,R,\oR\subset\E$ by B\'ezout's theorem.
It follows again from B\'ezout's theorem that $T,\oT,M,\oM\nsubseteq\E$.
We arrived at a contradiction, since $\tau(X)$ contains
three complex lines $\tau(R)$, $\tau(L)$ and $\tau(T)$
through the complex point~$\tau(\pp)$ instead of two.
We established that $L\nsubseteq \E$, and by using the same arguments
we find that $\oL,R,\oR,M,\oM,T,\oT\nsubseteq \E$ as well.
Since $X\cap \E$ does not contain complex lines,
we conclude from \Cref{lem:cq} that $X$ is not Cliffordian.
\end{proof}

We proceed in \Cref{exm:EOCO,exm:RPB}
to provide implicit equations for some great celestial Darboux cyclides.
This section is concluded with \Cref{rmk:EOCO,rmk:greatblum,rmk:greatperseus},
namely an analysis of
the geometries of great celestial Darboux cyclides by using the introduced methods.
The reader may opt to jump directly to \SEC{proof} at this point.

\begin{example}[great EO/CO cylides]
\label{exm:EOCO}
We consider the following surface
\[
X:=
\tau^{-1}(\set{y\in\P^3}{\alpha\,y_0^2 + \beta\,y_1^2 - y_2^2=0})
=
\set{x\in\S^3}{ \alpha\,x_1^2+\beta\,x_2^2-x_3^2=0},
\]
for some $\alpha,\beta\in\R_{>0}$.
Notice that $\set{y\in\P^3}{\alpha\,y_0^2 + \beta\,y_1^2 - y_2^2=0}$
is a ruled quadric and thus $X$ is great.
Suppose that $X'\subset\S^3$ is a CO cyclide or EO cyclide.
We claim that there exist $\alpha,\beta\in\R_{>0}$
\st $X'$ is M\"obius equivalent to $X$
and $\alpha=\beta$ if and only if $X$ is a CO cyclide.
We may assume up to M\"obius equivalence that the center~$p$ of~$\pi$ lies in $\Sing X_\R$.
The M\"obius transformations that leave $p$ invariant correspond via $\bR(\pi(\_))$
to Euclidean similarities in $\R^3$.
Hence, it follows from the classical result~\citep[Proposition~4]{circle}
that there exists $\varphi\in\aut\S^3$ and $\alpha,\beta\in\R_{>0}$ \st
$\varphi(p)=p$ and $(\pi\circ\varphi)(X')=\set{y\in\P^3}{\alpha\,y_1^2+\beta\,y_2^2-y_3^2}$.
Moreover, $\alpha=\beta$ if and only if $X'$ is a CO cyclide.
Since $\pi(x)=(x_0-x_4:x_1:x_2:x_3)$, we deduce that
$\varphi(X')=\set{x\in\S^3}{ \alpha\,x_1^2+\beta\,x_2^2-x_3^2=0}$ as was to be shown.
\END
\end{example}

\begin{figure}[!ht]
\centering
\fig{3}{5}{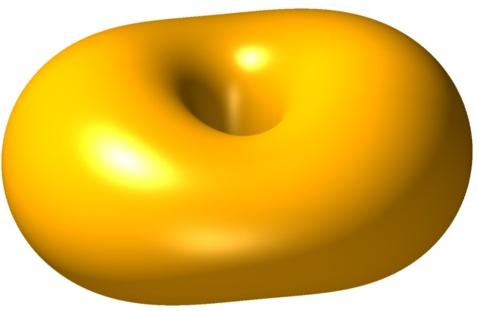}
\caption{Stereographic projection of a great Blum cyclide (see \Cref{exm:RPB}).}
\label{fig:blum}
\end{figure}

\begin{example}[great ring/Persues/Blum cyclide]
\label{exm:RPB}
We consider the following surface
\[
X:=\tau^{-1}(\set{y\in\P^3}{\alpha\,y_0^2 + y_1^2 - y_2^2 - \beta\,y_3^2=0})
=\set{x\in \S^3}{ \alpha\,x_1^2 + x_2^2 - x_3^2 - \beta\,x_4^2=0 },
\]
for some $\alpha,\beta\in\R_{>0}$.
We claim that $X\subset\S^3$ is a great ring cyclide, great Persues cyclide
or great Blum cyclide if $(\alpha,\beta)$ is equal to $(1,1)$, $(1,2)$ and $(2,2)$, \resp.
Since $\tau(X)$ is a doubly ruled quadric it follows that $X$ is great.
To identify $X$, let us first compute~$\Sing X$.
The Jacobian matrix of the generators of the ideal
$\bas{\alpha\,x_1^2 + x_2^2 - x_3^2 - \beta\,x_4^2,-x_0^2+x_1^2+x_2^2+x_3^2+x_4^2}$
is up to scaling of the rows as follows:
\[
\begin{bmatrix}
   0 & \alpha\,x_1 & x_2 & -x_3 & -\beta\,x_4\\
-x_0 &         x_1 & x_2 &  x_3 &         x_4
\end{bmatrix}.
\]
If $(\alpha,\beta)=(1,1)$,
then the Jacobian matrix has rank one at
the four complex points $(0:1:\pm \ii:0:0)$ and $(0:0:0:1:\pm \ii)$ in $X\cap \E$.
Hence, $X$ is a ring cyclide by \Cref{thm:B} (see the $\SingType X$ column in \Cref{tab:B}).
If $(\alpha,\beta)=(1,2)$,
then the Jacobian matrix has rank one at the complex points $(0:1:\pm \ii:0:0)$ in $X\cap \E$.
Hence, $X$ is a Perseus cyclide by \Cref{thm:B}.
If $(\alpha,\beta)=(2,2)$,
then the Jacobian matrix has rank two
at all complex points in~$X$ so that $\Sing X=\varnothing$.
Hence, $X$ is a Blum cyclide by \Cref{thm:B} and \Cref{prp:g}.
In \Cref{fig:blum} we depicted its stereographic projection
$
\bR(\pi(X))=\set{z\in\R^3}{
(x^2 + y^2 + z^2 )^2 - 6\,x^2 - 4\,y^2 + 1
}.
$
\END
\end{example}

\begin{remark}[great CO cyclide]
\label{rmk:EOCO}
Suppose that $X\subset \S^3$ is a great CO cyclide.
Recall from \Cref{exm:Y} that we encoded the corresponding
row in \Cref{tab:B2} in terms of the diagram in \Cref{fig:O} (left),
where $G(X)=\{\underline{g_0},\underline{g_1},g_{14},g_{23}\}$.
By \RP{C}{a}, the components $\{b_{34}\}$ and $\{b_{12}\}$
correspond to real antipodal singularities of $X$
and are centrally projected to the vertex of the quadratic cone~$\tau(X)$.
The great circles in~$X$ have class~$g_1$ and form a pencil with real antipodal base points in $\Sing X$.
The components $\{b_{13}'\}$ and $\{b_{24}'\}$
correspond to the complex isolated singularities that lie
in the ramification locus $\E$. The central projection of these
complex singular points are smooth complex branching points in~$\tau(X)$.
\END
\end{remark}

\begin{figure}[!ht]
\centering
\csep{1cm}
\begin{tabular}{cc}
\begin{tikzpicture}[scale=0.7] 
\draw[blue] (-2,1)  -- (2,1) node[right] {$e_3$};
\draw[blue] (-1,2) -- (-1,-2) node[below] {$e_2$};
\draw[red] (-2,-1) node[left] {$e_1$} -- (2,-1);
\draw[red] (1,2) node[above] {$e_4$} -- (1,-2);
\draw[draw=black, fill=green!5] (1,1)   circle [radius=5mm] node[black] {$b_{34}$};
\draw[draw=black, fill=green!5] (-1,-1)  circle [radius=5mm] node[black] {$b_{12}$};
\draw[draw=black, fill=white, densely dotted] (-1,1)  circle [radius=5mm] node[black] {$b_{13}'$};
\draw[draw=black, fill=white, densely dotted] (1,-1) circle [radius=5mm] node[black] {$b_{24}'$};
\node at (0,-3.5) {CO cyclide};
\end{tikzpicture}
&
\begin{tikzpicture}[scale=0.7] 
\draw[red, thick] (0.4,2.6) -- (-2,-1) node[left] {$e_{1}$};
\draw[red, thick] (-0.4,2.6) -- (2,-1) node[right] {$e_{2}$};
\draw[red] ( 0.2,3) node[above left] {$e_{13}$~~}  -- (-0.5,-0.6);
\draw[red] (-0.2,3) node[above right] {~$e_{14}$} -- ( 0.5,-0.6);

\draw[blue, thick] (0.4,-2.6) -- (-2,1) node[left] {$e_{12}$};
\draw[blue, thick] (-0.4,-2.6) -- (2,1) node[right] {$e_{11}$};
\draw[blue] ( 0.2,-3) node[below left] {$e_{3}$~~}  -- (-0.5,0.6);
\draw[blue] (-0.2,-3) node[below right] {~$e_{4}$} -- ( 0.5,0.6);

\draw[draw=black, fill=green!5] (0,2) circle [radius=5mm] node[black] {$b_{12}$};
\draw[draw=black, fill=green!5] (0,-2) circle [radius=5mm] node[black] {$b_{34}$};
\draw[draw=black, fill=white, densely dotted] (-1.3,0) circle [radius=2mm];
\draw[draw=black, fill=white, densely dotted] ( 1.3,0) circle [radius=2mm];
\draw[draw=black, fill=white, densely dotted] (-0.4,0) circle [radius=2mm];
\draw[draw=black, fill=white, densely dotted] ( 0.4,0) circle [radius=2mm];

\node at (-3.5,-3.5) {EO cyclide};
\end{tikzpicture}
\end{tabular}
\caption{Incidences between complex lines and isolated singularities on Darboux cyclides
(see \Cref{exm:Y} or the caption of \Cref{fig:Y}).
If the Darboux cyclide is great, then a red line is centrally projected to a blue line (see \Cref{rmk:EOCO}).
}
\label{fig:O}
\end{figure}

\begin{remark}[great Blum cyclide]
\label{rmk:greatblum}
Suppose that $X$ is the great Blum cyclide in \Cref{exm:RPB} with $(\alpha,\beta)=(2,2)$.
Our goal is to identify, up to $\aut N(X)$, the classes
of great and small circles, and pairs $([L],[L'])$ of classes
\st $L,L'\subset X$ are complex lines \st $\tau(L)=\tau(L')$.
As a byproduct, we recover the compact diagram in~\Cref{fig:greatblum}
from which we can read off the Clifford quartets
and the incidences between the complex lines in Blum cyclides.
By \Cref{thm:B}, we may assume up to $\aut N(X)$ that $B(X)$, $E(X)$ and $G(X)$ are as in \Cref{tab:B2}.
Since $\tau(X)$ is a doubly ruled quadric by \RL{gg}{c},
there exist great circles $C,C'\subset X$ \st $|C\cap C'|=2$.
By \RP{C}{b}, we may assume \Wlog that $([C],[C'])=(g_{12},g_{34})$ as
$g_{12},g_{34}\in\set{g\in G(X)}{\sigma_*(g)=g}$ and $g_{12}\cdot g_{34}=2$.
If $D,D'\subset X$ are antipodal small circles,
then $|D\cap D'\cap\E|=2$
and thus
$([D],[D'])\in\{(g_0,g_3),(g_1,g_2)\}$ by \RP{C}{b}.
Since the ramification locus~$X\cap\E$ of the central projection~$\tau$ does not contain complex lines,
it follows that for all complex lines $L\subset X$, there
exists a complex line $L'\subset X$ \st $\tau(L)=\tau(L')$.
The ramification locus~$\E$ is a hyperplane section,
which implies that $|L\cap L'\cap\E|=1$ and thus $|L\cap L'|=[L]\cdot [L']=1$ by \RP{C}{d}.
Let $a:=g_{12}\cdot ([L]+[L'])$ and $b:=g_{34}\cdot ([L]+[L'])$.
Notice that $\tau(L)$ belongs to one of the two rulings of $\tau(X)$
and thus $(a,b)$ is equal to either $(0,2)$ or $(2,0)$.
If $(a,b)=(0,2)$, then
$
([L],[L'])\in\{(e_1,e_2'),(e_2,e_1'),(e_{13},e_{04}),(e_{14},e_{03})\}
$
and if $(a,b)=(2,0)$, then
$
([L],[L'])\in\{(e_{01},e_{12}),(e_{02},e_{11}),(e_3',e_4),(e_4',e_3)\}.
$
See \Cref{fig:greatblum} for a diagrammatic representation of these pairs
\st
two line segments represent complex lines in the same pencil
on the doubly ruled quadric $\tau(X)$ if and only if the line segments
are both horizontal or both vertical.
The details are left to the reader.
\END
\end{remark}

\begin{figure}[!ht]
\centering
\vspace{-4mm}
\begin{tikzpicture}[xscale=1.8,yscale=0.7]
\draw[red ] (-1, 3) -- (4, 3) node[right] {$(e_{01},e_{12})$};
\draw[red ] (-1, 2) -- (4, 2) node[right] {$(e_{02},e_{11})$};
\draw[blue] (-1, 1) -- (4, 1) node[right] {$(e_{3}',e_{4})$};
\draw[blue] (-1, 0) -- (4, 0) node[right] {$(e_{4}',e_{3})$};
\draw[blue] (0, -1) -- (0, 4) node[above] {$(e_{1},e_{2}')$};
\draw[blue] (1, -1) -- (1, 4) node[above] {$(e_{2},e_{1}')$};
\draw[red ] (2, -1) -- (2, 4) node[above] {$(e_{13},e_{04})$};
\draw[red ] (3, -1) -- (3, 4) node[above] {$(e_{14},e_{03})$};
\foreach \x in {0,1,2,3}
   \foreach \y in {0,1,2,3}
      \draw[draw=black,fill=green] (\x,\y) circle [radius=0.6mm];
\end{tikzpicture}
\caption{
Each line segment is labeled with $([L],[L'])$, where $L$ and $L'$
are complex lines in a great Darboux cyclide $X\subset\S^3$
\st $\tau(L)=\tau(L')$. For each such label, we have $|L\cap L'|=1$.
Two line segments labeled with $([L],[L'])$ and $([M],[M'])$ meet at a green disc
if and only if $\{[L]\cdot[M],[L]\cdot[M']\}=\{0,1\}$ and $\{[L']\cdot[M],[L']\cdot[M']\}=\{0,1\}$.
The four Clifford quartets
are up to $\aut N(X)$ given by
\textcolor{blue}{$\{e_1 , e_2 , e'_3, e'_4\}$},
\textcolor{blue}{$\{e'_2, e'_1, e_3 , e_4\}$},
\textcolor{red}{$\{e_{13}, e_{14}, e_{01}, e_{02}\}$} and
\textcolor{red}{$\{e_{04}, e_{03}, e_{12}, e_{11}\}$}.
}
\label{fig:greatblum}
\end{figure}

\begin{remark}[great Perseus cyclide]
\label{rmk:greatperseus}
Let us describe the geometry of a great Perseus cyclide~$X\subset \S^3$,
by identifying the classes of great circles, small circles and complex lines,
and the components corresponding to base points.
By \Cref{thm:B}, we may assume up to $\aut N(X)$ that
\[
\set{g\in G(X)}{\sigma_*(g)=g}=\{g_0,g_1,g_{12},g_2,g_3\}.
\]
By \RP{C}{c}, there exist complex lines
$L,\oL,R,\oR,M,\oM,T,\oT\subset X$
\st
\[
[L]=e_{11},~
[\oL]=e_{12},~
[R]=e_{01},~
[\oR]=e_{02},~
[M]=e_3,~
[\oM]=e_4,~
[T]=e_3',~
[\oT]=e_4'.
\]
By \RP{C}{a}, the complex conjugate isolated singularities $\pp$ and~$\op$ in~$\Sing X$
correspond to the components $\{b_1\}$ and $\{b_2\}$, \resp.
The incidences between the complex lines and isolated singularities
are illustrated in \Cref{fig:gp} (see also \Cref{fig:pers}).

We know from \Cref{prp:perseus} that
$\tau(X)$ is a doubly ruled quadric, $\pp,\op\in\E$
and $L,\oL,R,\oR,M,\oM,T,\oT\nsubseteq \E$.
It follows that either $\tau(L)=\tau(R)$, $\tau(L)=\tau(M)$ or $\tau(L)=\tau(T)$.

First, suppose by contradiction that $\tau(L)=\tau(R)$.
In this case $\tau(\oL)=\tau(\oR)$.
This is a contradiction, since $\tau(L)$
intersects its complex conjugate line $\tau(\oL)$
and thus $\tau(X)$ must be an ellipsoid instead of being a doubly ruled quadric.

Next, we suppose that $\tau(L)=\tau(M)$.
In this case, the complex conjugate lines $\tau(L)=\tau(M)$ and $\tau(\oL)=\tau(\oM)$
belong to the first pencil of lines on the doubly ruled quadric $\tau(X)$.
The complex conjugate lines, $\tau(R)=\tau(T)$ and $\tau(\oR)=\tau(\oT)$
belong to the second pencil of lines on $\tau(X)$.
By Proposition~\ref{prp:C}, a circle with class~$g_1$ meets each of the lines
in $\{L,\oL,M,\oM\}$ and belongs to a base point free pencil.
Similarly, a circle with class $g_2$ meets each of the lines
in $\{R,\oR,T,\oT\}$ and belongs to a base point free pencil.
From this we establish that $g_1$ and $g_2$ correspond to pencils of great circles
that are centrally projected the first and second ruling of the quadric~$\tau(X)$, \resp.
A circle $C''\subset X$ \st $[C'']$ equals either $g_0$ or~$g_3$ meets each of the lines in
$\{L,\oL,T,\oT\}$ and $\{R,\oR,M,\oM\}$, \resp.
Therefore, each circle $C\subset X$ \st $[C]=g_0$
is a small circle whose antipodal points
form a small circle~$C'\subset X$ \st $[C']=g_3$ and $\tau(C)=\tau(C')$.
A circle with class $g_{12}$ passes through the complex conjugate isolated singularities~$\pp$ and~$\op$.
Hence, each circle $C\subset X$ \st $[C]=g_{12}$
is a small circle whose antipodal points
form a small circle $C'\subset X$ \st $[C']=g_{12}$ and $\tau(C)=\tau(C')$.

The case $\tau(L)=\tau(T)$
is analoguous to the previous case: $\{g_0,g_3\}$ are classes of great circles,
$\{g_1,g_2\}$ are the classes of antipodal little circles, and $g_{12}$ is the class of
a small circle that meets $\pp$ and $\op$.
The details are left to the reader.
\END
\end{remark}

\section{Combining the results}
\label{sec:proof}

In order to prove \Cref{thm:x}
we use the following theorem from \citep[Theorem~1]{circle}.

\begin{theoremext}
\label{thm:8}
If $X\subset\S^3$ is a $\lambda$-circled surface of degree $d$ \st $\lambda\geq 2$,
then either $X$ is either a Darboux cyclide or
$(\lambda,d)\in\{(\infty,2),(2,8)\}$.
\end{theoremext}

\begin{lemma}
\label{lem:sphere}
If $X\subset \S^3$ is a surface \st $\bR(X)$ is a 2-dimensional sphere, then $X$ is not Cliffordian.
\end{lemma}

\begin{proof}
Suppose by contradiction that $X=\bP(A\star B)$ for some circles $A,B\subset S^3$.
The left Clifford translations correspond to isoclinic rotations of~$S^3$.
Thus the infinitesimal left Clifford translations of points on the circle $B$
define a nowhere vanishing vector field on~$\bR(X)\subset S^3$.
We arrived at a contradiction, since a 2-dimensional sphere does not admit such a vector field
by the hairy ball theorem.
\end{proof}

\begin{proof}[Proof of \Cref{thm:x}.]
Since $Z\subset\R^3$ is $\lambda$-circled and of M\"obius degree $d$
\st $(d,\lambda)\neq (8,2)$, it follows from \Cref{thm:8} that $d\in\{2,4\}$.
If $d=2$, then $Z$ is either a plane or a 2-dimensional sphere.

\ref{thm:x:a}
A CY or EY is always Bohemian as it can be obtained
by translating a circle along a line (see \Cref{fig:4}).
A plane is the translation of a line along a line.
Hence the proof for this assertion is concluded by \Cref{prp:b}.

\ref{thm:x:b}
By \Cref{lem:sphere} we have $d=4$
and thus the first part follows from \Cref{prp:c},
where $X=\bS(Z)$.
It follows from \Cref{prp:g,prp:ring}
that ring cyclides are M\"obius equivalent to Cliffordian surfaces.

\ref{thm:x:c}
Direct consequence of \Cref{prp:g}, where $X=\bS(Z)$.
\end{proof}

\Cref{cor:moebius,cor:bohemian} are direct consequences of \Cref{thm:x}.

\begin{proof}[Proof of \Cref{cor:g2}.]
The central projection of a surface~$Z\subseteq S^n$
that is covered by two pencils of great circles
is a doubly ruled quadric.
Therefore, $n=3$ and $Z$ has no real singularities.
In particular, the stereographic projection~$\mu(Z)$
is not a CO cyclide or EO cyclide.
Thus the proof is concluded by \Cref{thm:x}\ref{thm:x:c}.
\end{proof}

\begin{proof}[Proof of \Cref{cor:great}.]
Direct consequence of \Cref{prp:ring,prp:perseus}.
\end{proof}

\section{Acknowledgements}

I would like to thank Mikhail~Skopenkov for
the proof idea of \Cref{prp:hqf}
and many other insightful and detailed remarks and corrections.
The surface figures were generated using \citep[Sage]{sage}.
This work was supported by the Austrian Science Fund (FWF) project P33003.

\bibliography{product-4}

\paragraph{address of author:}
Johann Radon Institute for Computational and Applied
Mathematics (RICAM), Austrian Academy of Sciences
\\
\textbf{email:} info@nielslubbes.com

\end{document}